\documentclass[aap]{imsart}

\RequirePackage{amsthm,amsmath,amsfonts,amssymb}
\RequirePackage[numbers,sort&compress]{natbib}
\RequirePackage[colorlinks,citecolor=blue,urlcolor=blue]{hyperref}
\RequirePackage{graphicx}

\usepackage{tikz}
\usetikzlibrary{positioning}
\usepackage{eucal}
\usepackage{amssymb,epsfig,latexsym,bm}
\usepackage{psfrag}
\usepackage{xcolor,comment,bm}
\newcommand{\dd}{\mathrm d}

\usepackage{amsthm}

\newtheoremstyle{italictheorem}
  {\topsep}   
  {\topsep}   
  {\itshape}  
  {0pt}       
  {\bfseries} 
  {.}         
  {5pt plus 1pt minus 1pt} 
  {}          

\startlocaldefs
\numberwithin{equation}{section}

\theoremstyle{remark}

\theoremstyle{italictheorem} 

\theoremstyle{italictheorem} 
\newtheorem{Theorem}{Theorem}[section]
\newtheorem{Proposition}[Theorem]{Proposition}

\newtheorem{Definition}[Theorem]{Definition}
\newtheorem{Corollary}[Theorem]{Corollary}
\newtheorem{Lemma}[Theorem]{Lemma}
\newtheorem{Remark}[Theorem]{Remark}

\def\ind{{\rm 1\hspace{-0.90ex}1}}

\newcommand{\supp}{\text{supp}}

\newcommand{\primotheokappab}{h_B}
\newcommand{\secondotheokappazero}{h_0}
\newcommand{\secondotheokappauno}{h_1}
\newcommand{\secondotheokappazerob}{h_B^{(0)}}
\newcommand{\secondotheokappaunob}{h_B^{(1)}}

\newcommand{\terzotheokappazero}{h_0}
\newcommand{\terzotheokappazerob}{h_B^{(0)}}
\newcommand{\nsupra}{non-suprathreshold }
\newcommand{\supra}{suprathreshold}

\newcommand{\Rsupra}{$R$-suprathreshold}
\newcommand{\Bsupra}{$B$-suprathreshold}
\newcommand{\Ssupra}{$S$-suprathreshold}

\newcommand{\FH}{\zeta}
\newcommand{\sidebyside}[4]{
	\begin{figure}[phtb]
		\begin{minipage}[phtb]{6.5cm}
			\includegraphics[width=\textwidth]{#1.pdf}
			\caption{#2\label{fig:#1}}
		\end{minipage}
		\hfill \hspace{0.3cm}
		\begin{minipage}[phtb]{6.5cm}
			\includegraphics[width=\textwidth]{#3.pdf}
			\caption{#4\label{fig:#3}}
		\end{minipage}
		\hfill
		\vspace{-0mm}
	\end{figure}
}

\usepackage{ulem,cancel}

\usepackage{changepage}
\usepackage{ifthen}
\endlocaldefs

\begin{document}

\begin{frontmatter}
\title{Competing bootstrap processes \\ on the random graph $G(\lowercase{n},\lowercase{p})$}
\runtitle{Competing bootstrap processes}

\begin{aug}

\author[A]{\fnms{Michele}~\snm{Garetto}\ead[label=e1]{michele.garetto@unito.it}}, 
\author[B]{\fnms{Emilio}~\snm{Leonardi}\ead[label=e2]{emilio.leonardi@polito.it}}
\and 
\author[C]{\fnms{Giovanni Luca}~\snm{Torrisi}\ead[label=e3]{giovanniluca.torrisi@cnr.it}}

\address[A]{Universit\`a di Torino, Torino, Italy.
\printead[presep={,\ }]{e1}}
	
\address[B]{Politecnico di Torino, Torino, Italy.
\printead[presep={,\ }]{e2}}

\address[C]{Consiglio Nazionale delle Ricerche,  Roma, Italy.
	\printead[presep={,\ }]{e3}}

\end{aug}

\begin{abstract}
		We introduce and analyze a competitive extension of classical bootstrap
		percolation on the Erd\H{o}s--R\'{e}nyi random graph $G(n,p_n)$.
		Nodes can be red, black, or white, with red and black representing two
		competing active states. Starting from two sets of initially active seeds,
		white nodes activate according to independent Poisson clocks and become
		permanently red or black whenever the number of active neighbors of one
		color exceeds that of the competing color by at least a fixed threshold
		$r\geq2$.
		
		We characterize the asymptotic dynamics and final sizes of the two
		competing activation processes over all relevant seed-density scales.
		Our results reveal a sharp qualitative dichotomy. In the sub-critical
		regime, competition is asymptotically negligible at first order: each
		process reaches the same normalized final size as it would in the
		absence of the competing process. In the super-critical regime, instead,
		the initial advantage of the process with the larger seed density is
		amplified: the dominant process activates $n-o(n)$ nodes, while the
		competing process is suppressed and remains confined to a much smaller
		scale.
		
		Beyond final-size asymptotics, we derive a fluid-limit description of
		the joint activation trajectories and characterize the relevant
		activation time scales. The analysis combines uniform concentration
		estimates, deterministic limiting Cauchy problems, stochastic coupling,
		and multiscale arguments.
\end{abstract}



\end{frontmatter}

\section{Introduction}

Bootstrap percolation is an activation process on a graph that begins with a set of initially active nodes (seeds). The process unfolds in discrete rounds: any inactive node with at least $r\ge 2$ active neighbors becomes active, and once active, remains so permanently (the process is irreversible). In each round, all eligible nodes activate simultaneously, and the process terminates when no further activations are possible.

Like many percolation processes, bootstrap percolation exhibits \lq\lq all-or-nothing" behavior: either the activation spreads to nearly all nodes in the graph, or it quickly ceases, resulting in a final number of active vertices only slightly larger than the initial seeds. The process is said to almost percolate if the final number of active nodes is $n-o(n)$ as $n\to\infty$ (rigorous definitions of the asymptotic notation used in this paper can be found in the next section).

Historically, bootstrap percolation was first introduced on a Bethe lattice \cite{chalupa} and later explored on regular grids \cite{kogut, enter, aizenman,
gravner, morris} and trees \cite{BPP, gunderson}. More recently, its study has expanded to various random graphs, driven by growing interest in large-scale complex systems such as technological, biological, and social networks.

A key contribution in this direction comes from Janson et al. \cite{JLTV}, who provided a detailed analysis of the bootstrap percolation process on the Erd\H{o}s--R\'{e}nyi random graph $G(n,p_n)$. Their work identifies a critical size $a_c^{(n)}$ for the initial number of seeds: if the number of seeds asymptotically exceeds $a_c^{(n)}$, the bootstrap percolation process spreads throughout nearly the entire graph; otherwise, the process largely ceases to develop. We note that the analysis in \cite{JLTV} considers seeds selected uniformly at random. However, subsequent studies have shown that the critical threshold for percolation can be considerably reduced when seed selection is optimized through the formation of "contagious sets" \cite{feige,GMM}.

Related to our work is the study in \cite{inhibition}, which explores a variant of classical bootstrap percolation on the random graph $G(n,p_n)$. In their model, nodes are classified as either excitatory or inhibitory, and activation spreads to nodes where the number of active excitatory neighbors sufficiently outweighs the number of active inhibitory neighbors. Interestingly, when more than half the nodes are inhibitory, they observe non-monotonic effects on the final active size in the traditional round-based model. These effects disappear in a continuous-time setting that incorporates exponential transmission delays on edges. While we also utilize a continuous-time framework, our exponential delays are placed on nodes rather than edges. Furthermore, our model differs significantly from \cite{inhibition} because we investigate the competition between two opposing activation processes. Another related variant is majority bootstrap percolation \cite{cecilia17}, where a node activates if at least half of its neighbors are active.

Large deviations in classical bootstrap percolation on $G(n,p_n)$ have also been studied. In \cite{largedevsub}, the authors calculate the rate function for the event where a small (subcritical) set of initially active nodes unexpectedly activates a large number of vertices, also identifying the most probable "least-cost" trajectory for such events. Large deviations in the supercritical regime were fully characterized in our previous work \cite{LDP-nostro}.

Bootstrap percolation has been analyzed on various other graph types, including random regular graphs \cite{pittel}, random graphs with given vertex degrees \cite{amini1}, and random geometric graphs \cite{rgg}. It has also been explored on Chung--Lu random graphs \cite{amini2,fountolakis}, which are particularly useful for modeling power-law node degree distributions, as well as on small-world \cite{kozma16,turova15} and Barabasi--Albert random graphs \cite{fountolakis2}. In \cite{Bernoulli}, we examined bootstrap percolation on the stochastic block model, an extension of the Erd\H{o}s--R\'{e}nyi random graph that captures the community structures prevalent in many real networks.

The main contribution of this paper is not merely to extend bootstrap
percolation to a different graph structure, but to introduce and analyze
a genuinely competitive activation mechanism. Our analysis reveals two
qualitatively different behaviors. In the sub-critical regime, the two
processes are asymptotically decoupled at first order, whereas in the
super-critical regime the initial advantage of one process is amplified,
leading to almost complete percolation of the dominant color and to the
early suppression of the competing process. Beyond final-size results,
we characterize the activation trajectories and the relevant physical
time scales through fluid-limit and multiscale arguments.


The remainder of the paper is organized as follows.
Section~\ref{Sect:model} introduces the model and states the main results, together
with a numerical illustration. Section~\ref{Sect:PA} develops the preliminary
probabilistic representation of the process. Sections~\ref{Sect:q-timescalle-res} and~\ref{sec:V}  analyze the initial evolution of the activation process and provide
the corresponding proofs.
Sections~\ref{sub:scaling} and~\ref{sect-theo->q} 
extend the analysis to later stages of the
activation process and establish the asymptotic behavior of the
system when the activation propagates to a macroscopic fraction of
the graph.
Technical and auxiliary proofs
are deferred to the appendices.

\section{Model and Main Results}
\label{Sect:model}
\subsection{Notation}
	
Throughout the paper, all unspecified limits are as $n\to\infty$. 
We will use the following standard asymptotic notation. Given two numerical sequences 
$\{f(n)\}_{n\in\mathbb N}$ and $\{g(n)\}_{n\in\mathbb N}$, $\mathbb N:=\{1,2,\ldots\}$,
we write: $f(n)\ll g(n)$ if $f(n)=o(g(n))$, i.e., $f(n)/g(n)\to 0$; $f(n)=O(g(n))$ if $\limsup\Big|\frac{f(n)}{g(n)}\Big|<\infty$; $f(n)=\Theta(g(n))$ if $f(n)=O(g(n))$ and $g(n)=O(f(n))$;
$f(n)\sim g(n)$ if $f(n)/g(n)\to 1$.  Unless otherwise stated, all random quantities considered in this paper are defined on an underlying probability space $(\Omega,\mathcal{F},\mathbb P)$.  Let $\{X_n\}_{n\in\mathbb N}$ be a sequence of real-valued random variables.  We write $X_n=o_{\mathrm{a.s.}}(f(n))$ if $\mathbb{P}\left(\lim\Big|\frac{X_n}{f(n)}\Big|=0\right)=1$; $X_n=O_{\mathrm{a.s.}}(g(n))$
if $\mathbb{P}\left(\limsup\Big|\frac{X_n}{g(n)}\Big|<\infty\right)=1$;
$X_n=\Theta_{\mathrm{a.s.}}(g(n))$ if 
\[
\mathbb{P}\left(\limsup\Big|\frac{g(n)}{X_n}\Big|<\infty\right)=\mathbb{P}\left(\limsup\Big|\frac{X_n}{g(n)}\Big|<\infty\right)=1.
\]

We denote by $\|\cdot\|$ the Euclidean norm on $\mathbb R^d$ for some $d\in\mathbb N$, and by $\lfloor\cdot\rfloor$ and $\lceil\cdot\rceil$ the floor and ceiling functions, respectively.
Given a set $\mathcal A$, we denote by ${\mathcal A}^c$ its complement and by  
$|\mathcal{A} |$ its cardinality. 
Let $X$ and $Y$ denote two real-valued random variables.
We denote by $X \leq_{\mathrm{st}} Y$ the usual stochastic order, i.e., we write $X \leq_{\mathrm{st}} Y$ if $\mathbb P(X>z)\leq \mathbb P(Y>z)$ for $z\in\mathbb R$. Hereafter, the symbols ${\rm Be}(u)$,
$\text{Bin}(m,\theta)$, $\text{Po}(\lambda)$ and $\text{Exp}(\lambda)$ 
denote random variables
distributed according to the Bernoulli law with mean $u\in [0,1]$, the binomial law with parameters $(m,\theta)$, the Poisson law, and the exponential law, both with parameter $\lambda>0$, respectively. The symbol $\overset{\mathrm L}{=}$ denotes equality in law.

Finally, throughout this paper we will use the function
\begin{equation}\label{eq:H}
	\FH(x):=1-x+x\log x,\quad x>0, \qquad \FH(0):=1.
\end{equation}
In the following,  with occasional exceptions dictated by standard usage (e.g.  $\Omega$ for the sample space),    we adopt the following conventions: $(i)$  capital letters denote random variables; $(ii)$ lowercase letters indicate deterministic quantities, including constants, parameters and  functions; $(iii)$ capital calligraphic letters denote set-valued random variables,  events and sigma-algebras;
	$(iv)$ boldface  symbols indicate vectors; $(v)$ blackboard bold capital letters denote  sets of points or numbers and probability measures.

\subsection{Model description}\label{sect:mod-descr}
We consider a generalization of the bootstrap percolation process
on the Erdős--Rényi random graph $G(n,p_n)=(\mathcal{V}^{(n)},\mathcal{E}^{(n)})$, $n\in\mathbb N$, introduced in \cite{JLTV}.
The graph consists of a node set $\mathcal{V}^{(n)}:=\{1,\ldots,n\}$ and an edge set $\mathcal{E}^{(n)}$, where each potential edge between two distinct nodes is included independently with probability $p_n\in (0,1)$.
Our model is defined as follows:

\begin{adjustwidth}{1cm}{0cm} 
\begin{description}
\item[{\bf Node states:}] Nodes can be in one of three states: red ($R$), black ($B$), or white ($W$). We refer to $R$ and $B$ nodes as active nodes, and to $W$ nodes as inactive nodes.

\item[{\bf Initial condition:}] At time $0$, an arbitrary number $a_R^{(n)}$ of nodes are chosen uniformly at random and set to $R$. Subsequently, an arbitrary number $a_B^{(n)}$ of nodes are selected uniformly at random  from the remaining $n-a_R^{(n)}$ nodes and set to $B$. These nodes, active at time $0$, are referred to as seeds.\footnote{Since the seeds are selected uniformly at random in $G(n,p_n)$, the order in which the two sets of seeds are generated is not relevant, i.e., it has no impact on the evolution of the bootstrap percolation processes.} All other nodes are initially set to $W$.

\item[{\bf Activation mechanism:}] Each white node has an independent Poisson process (with intensity 1) attached to it, which dictates when the node "wakes up".
When a white node wakes up, it assesses its neighbor states to decide whether to change its color to either $R$ or $B$.
A $W$ node changes its state to $S\in \{R,B \}$ if the number of its neighbors with color $S$ exceeds the number of its neighbors with the opposite color $\overline{S}$ (if $S$ is red, $\overline{S}$ is black, and vice versa) by at least $r\in \mathbb{N}\setminus \{1\}$. 
Throughout this paper we refer to this condition as the "threshold condition with respect to $S$"
{and to nodes satisfying it as  \Ssupra{}  nodes}.
Otherwise, the node stays white.

\item[{\bf Irreversibility:}] Once active (either $R$ or $B$), a node remains so indefinitely, meaning that it cannot deactivate or change its color. This ensures that the total number of $R$ and $B$ nodes in the system is non-decreasing.

\item[{\bf Termination condition:}] The process stops when no more nodes can be activated, i.e., no $W$ node satisfies the "threshold condition with respect to either $R$ or $B$".
\end{description}
\end{adjustwidth}

\begin{Remark} \label{remark_Poisson}
	Unlike the bootstrap percolation process in \cite{JLTV}, where the activation order does not affect the final number of active nodes (as noted in Proposition 4.1 of \cite{Bernoulli}), in our model the activation order is crucial, as toy examples demonstrate. To circumvent this problem, we have introduced Poisson clocks on the nodes. Essentially, this allows us to model a system where, at any given time, the next node to activate is chosen uniformly at random from those satisfying the threshold condition with respect to either $R$ or $B$.
\end{Remark}

The aim of this paper is to study the asymptotic behavior of the final number 
$A_R^{*(n)}$ ($A_B^{*(n)}$) of nodes $R$ ($B$) as $n$ grows large.
Following a common practice in the theory of large random graphs, we will omit the dependence on $n$ of the various mathematical objects or quantities, writing e.g.  $p$ in place of $p_n$,  $G(n,p)$ or simply $G$ in place of $G(n,p_n)$,  $a_S$ in place of $a_S^{(n)}$, $A_S^*$ in place of $A_S^{*(n)}$, $S\in\{R,B\}$, and so on.
We will specify such dependence only when necessary to avoid confusion.
{We remark that the threshold $r$ is supposed to be constant, i.e., not depending on $n$.}

\begin{Remark}\label{remark-bootstrap}
	
	When $a_{\overline{S}}= 0$,
	our process reduces to an asynchronous variant of the 
	classical bootstrap percolation model studied in  \cite{JLTV},
where the next node to activate is chosen uniformly at random among nodes   satisfying the threshold condition. 
Consequently,  as noted in  Remark \ref{remark_Poisson},
	$A_S^*$ matches the final count of active nodes in a classical bootstrap 
	percolation process on $G$ with $r\geq 2$ and $a_S$ seeds.
	
\end{Remark}

Throughout this paper we assume that 
\begin{equation}\label{eq:bootcond}
	\frac{1}{n}\ll p\ll\frac{1}{n^{1/r}\log n}. 
\end{equation}

{Condition \eqref{eq:bootcond} strengthens the assumption
	$p\ll n^{-1/r}$ used in \cite{JLTV}. Indeed, it implies
	$g/\log n\to\infty$, making the exponential error bounds arising in
	our proofs summable. The Borel--Cantelli lemma then yields almost sure
	convergence, whereas the weaker assumption is sufficient only for
	convergence in probability.}

Our model of competing bootstrap percolation gives rise to different regimes
depending on how $a_R$ and $a_B$ scale with $n$.
As in \cite{JLTV}, we first define the critical seed-set size of standard bootstrap percolation in $G$ (the meaning of $g$ is explained in Remark \ref{BPRemark}):
\begin{equation}\label{eq:gcrit}
g:=\left(1-\frac{1}{r}\right)\left(\frac{(r-1)!}{n p^r}\right)^{\frac{1}{r-1}} \qquad \text{(note that $pg\to 0$).}
\end{equation}
We consider the following different choices of sequence $\{q_n\}$ (hereinafter written simply as $q$, and also referred to as the system \lq\lq time-scale\rq\rq ):
\begin{equation}
	\begin{split} \label{eq:grelq}
		&\text{$(i)$} \quad q=g; \quad 	
		\text{$(ii)$}\quad g\ll q \ll p^{-1}; \quad   
		\text{$(iii)$}\quad q=p^{-1}; \quad \text{$(iv)$} \quad p^{-1}\ll q\ll n. 
	\end{split}
\end{equation}
and we assume that:
\begin{equation}\label{eq:trivial}
	a_R/q \to \alpha_R,\quad a_B/q\to\alpha_B,\quad\text{for some positive constants $\alpha_R,\alpha_B>0$.}
\end{equation}

Without loss of generality, 
we will always assume $\alpha_R>\alpha_B$, 
deferring the analysis of the case $\alpha_R=\alpha_B$ to future studies.
{
	\begin{Remark}
		We do not analyze the regime $q\ll g$ further.
		If $a_S/q\to\alpha_S>0$, then
		\[
		\frac{A_S^*-a_S}{q}\xrightarrow{\mathbb P}0,
		\qquad\text{and hence}\qquad
		\frac{A_S^*}{q}\xrightarrow{\mathbb P}\alpha_S,
		\qquad S\in\{R,B\}.
		\]
		
		Indeed, let $\overline A_S^*$ denote the final number
		of active vertices in the corresponding single-color
		bootstrap percolation process with $a_S$ seeds.
		The classical results in \cite{JLTV,Bernoulli} give
		\[
		\frac{\overline A_S^*}{q}
		\xrightarrow{\mathbb P}\alpha_S.
		\]
		The claim then follows from
		Proposition~\ref{lemma-compare}, which gives
		$A_S^*\leq_{\mathrm{st}}\overline A_S^*$,
		together with the deterministic bound $A_S^*\geq a_S$.
		
		Unlike the regimes treated in our main theorems,
		$q\ll g$ does not by itself provide summable
		deviation probabilities at scale $q$.
		We therefore make no almost sure convergence claim
		in this regime without additional assumptions.
	\end{Remark}
}
\begin{Remark}
	We do not explore the $q\ll g$ scenario since it yields straightforward results. The analysis from \cite{JLTV, Bernoulli} indicates that classical bootstrap percolation barely evolves under this condition, meaning 
	$A_R^*= \alpha_R q+o_{\mathrm{a.s.}}(q)$. This behavior extends to our model with two competing bootstrap processes, a claim directly supported by Proposition \ref{lemma-compare}.
\end{Remark}
\begin{Remark} \label{BPRemark}
Under the condition  $\frac{1}{n}\ll p\ll\frac{1}{n^{1/r}}$, the main results from \cite{JLTV} provide the asymptotic behavior of $A_R^*$ when
$a_R/g\to \alpha_R$ and $a_B=0$. 
Specifically, $A_R^*/g\to z_R+\alpha_R$ 
in probability if $\alpha_R<1$, while $A_R^*/n\to 1$ in probability if $\alpha_R>1$ (a precise definition of $z_R$ will be provided in Remark \ref{re:29Ott}). This implies the existence of a critical threshold for the number of seeds: below it, the bootstrap percolation process remains largely unchanged, but above it, the bootstrap percolation process percolates almost the entire graph.
\end{Remark}
The well-known behavior of classical bootstrap percolation described above motivates the following terminology for the model introduced in this paper:
$(i)$ We say that the system is in the sub-critical regime when $q=g$ 
and $\alpha_{R}<1$;
$(ii)$ We say that the system is in the super-critical regime 
if either $q=g$ and $\alpha_{R}>1$, or $g\ll q\ll n$.

\subsection{Main results}

To state our results we need to introduce the following function $\beta_S:[0,\infty)^2\to \mathbb R$, $S\in\{R,B\}$:
\begin{align}\label{eq:beta0}
	\beta_S(x_R,x_{B}):=
	\left\{ 
	\begin{array}{cc}
		r^{-1}(1-r^{-1})^{r-1}(x_S+\alpha_S)^r-x_S & \text{if $q=g$} \\
		\frac{1}{r!}(x_S+\alpha_S)^r &\text{if $g\ll q\ll p^{-1}$} \\
		\sum_{r'=r}^\infty\sum_{r''=0}^{r'-r} \frac{(x_S+\alpha_S)^{r'}}{r'!}
		\frac{(x_{\overline{S}}+\alpha_{\overline{S}})^{r''}}{r''!} {\mathrm e}^{-(x_R+x_{B}+\alpha_R+\alpha_{B} )}&
		\text{if $q=p^{-1}$}\\
		\ind_{[0,\infty)}\left(\frac{x_S+\alpha_S}{ x_R +\alpha_R+ x_B+\alpha_{B}} -\frac{1}{2}\right) 
		& \text{if $p^{-1}\ll q \ll n$,}\\
	\end{array} \right.
\end{align}
where $\ind_{\mathcal A}(\cdot)$ denotes the indicator function of the set $\mathcal A$.
{
Along the physical trajectories and on the domains
used below, $\beta_S$ describes the leading-order
activation rate after the appropriate normalization,
as made precise in Proposition~\ref{cor:BC}.
Here $x_Sq$ represents the number of non-seed
$S$-activations, while the initial seeds contribute
through $a_S\sim\alpha_Sq$.
In particular, when $q=g$, the term $-x_S$ accounts
for vertices that have already activated and are
therefore no longer available for activation.

In the regime $p^{-1}\ll q\ll n$, the values assigned
by \eqref{eq:beta0} on the hyperplane
\[
x_R+\alpha_R=x_B+\alpha_B
\]
are merely conventional. They are not asserted
to be the limiting probabilities of satisfying
the two threshold conditions on that hyperplane.
The $q$-scale uniform approximation in this regime
is used only on compact regions that remain
a positive distance from the hyperplane.

This restriction is sufficient for our analysis:
under $\alpha_R>\alpha_B$, the relevant limiting
trajectory is
\[
\bm f(x)=(x,0),
\]
as stated in Proposition~\ref{le:coupledphi}.
It therefore remains strictly on the red-dominated
side of the hyperplane, with a positive separation
on every compact interval.}


\begin{Remark}\label{re:29Ott}				
For $q=g$, the sign of $\beta_S$ is determined by $\alpha_S$: it is strictly positive for any $x_S\ge 0$ when $\alpha_S>1$; when $\alpha_S=1$, $\beta_S$ has one strictly positive zero, say $z_S$; when $\alpha_S<1$, it has two strictly positive zeros; letting $z_S$ denote the smaller one, it turns out that $\beta_S$ is strictly decreasing in the interval $(0,z_S)$.
If either $g\ll q\ll p^{-1}$ or $q=p^{-1}$, then $\beta_S$ is strictly positive in the whole domain. $\beta_S$ is non-negative if $p^{-1}\ll q\ll n$.
\end{Remark}

\begin{Remark}\label{remark-BPthreshold}			
We have excluded the case $r=1$ from the analysis of the competing bootstrap processes since when $r=1$, the classical bootstrap percolation itself has a qualitatively different behavior. Indeed, a single seed that lies in the giant component
is enough to trigger
an almost complete graph percolation (see Remark 5.9 in \cite{JLTV}). 
This phenomenon fundamentally removes the sub-critical phase and the existence of a critical threshold.
As a consequence, the analysis of competition between two bootstrap processes with $r=1$ requires substantially different techniques, as it necessitates considering finite seed sets (i.e., those that don't scale with $n$). The exploration of the $r=1$ case in our model is reserved for future studies.
\end{Remark}	

Consider the system evolution within the sub-critical regime. One might intuitively expect that competition would lead to smaller asymptotic final sizes for $S$-active nodes ($S\in \{R,B\}$) compared to scenarios without competition (i.e., when $a_{\overline{S}}=0$). However, the following theorem shows that this is not the case.

\begin{table}[t]
	\centering
	\tiny
	\renewcommand{\arraystretch}{1.25}
	\begin{tabular}{|p{0.22\textwidth}|p{0.18\textwidth}|p{0.23\textwidth}|p{0.27\textwidth}|}
		\hline
		\textbf{Regime}
		&
		\textbf{$R$-final size}
		&
		\textbf{$B$-final size}
		&
		\textbf{Qualitative behavior}
		\\
		\hline
		
		$q=g$, $\alpha_R<1$
		&
		$\displaystyle
		\frac{A_R^*}{q}
		\to z_R+\alpha_R$
		&
		$\displaystyle
		\frac{A_B^*}{q}
		\to z_B+\alpha_B$
		&
		The two processes remain at scale $q$ and are asymptotically
		decoupled at first order.
		\\
		\hline
		
		$q=g$, $\alpha_R>1$
		&
		$\displaystyle
		\frac{A_R^*}{n}\to1$
		&
		$\displaystyle
		\frac{A_B^*}{q}
		\to g_B(\kappa_{\bm g})+\alpha_B$
		&
		The initial advantage of the $R$-process is amplified:
		$R$ almost percolates, while the $B$-process remains confined
		to scale $q$.
		\\
		\hline
		
		$g\ll q\ll p^{-1}$
		&
		$\displaystyle
		\frac{A_R^*}{n}\to1$
		&
		$\displaystyle
		\frac{A_B^*}{n}\to0$
		&
		The processes are weakly coupled at the initial scale $q$,
		but the $R$-process eventually dominates and almost percolates.
		\\
		\hline
		
		$q=p^{-1}$
		&
		$\displaystyle
		\frac{A_R^*}{n}\to1$
		&
		$\displaystyle
		\frac{A_B^*}{n}\to0$
		&
		Competition is already visible at scale $q$; the $R$-process
		eventually dominates.
		\\
		\hline
		
		$p^{-1}\ll q\ll n$
		&
		$\displaystyle
		\frac{A_R^*}{n}\to1$
		&
		$\displaystyle
		\frac{A_B^*}{n}\to0$
		&
		Strong interaction occurs already at scale $q$, and the
		$R$-process dominates the evolution.
		\\
		\hline
	\end{tabular}
	
	\caption{Summary of the asymptotic behavior of the two competing
		activation processes. Throughout the table, $\alpha_R>\alpha_B$.
		All convergences hold almost surely.}
	\label{tab:phase-diagram}
\end{table}

\begin{Theorem}[{\textbf{sub-critical regime}}]\label{thm:subcriticostoc}
{\it Assume $q=g$ 
	with $\alpha_{R}<1$. Then
	\[
	\frac{A^*_R}{q}\to z_R+\alpha_R\quad\text{and}\quad\frac{A^*_B}{q}\to z_{B}+\alpha_B,\qquad {\text{a.s.}}
	\]
	where $z_S$ is the smallest zero of $\beta_S$ $($see Remark \ref{re:29Ott}$)$.}
\end{Theorem}

Theorem \ref{thm:subcriticostoc} states that, in the sub-critical regime, the two competing processes essentially do not interact with each other.
Indeed, $A^*_S/q$
converges exactly to the same value it would converge to if
$a_{\overline{S}}=0$ (see Remark \ref{BPRemark}).

Next, we consider 
the more interesting super-critical regime.

\begin{Theorem}[{\textbf{super-critical regime}}]\label{prop:supercritical}
{\it The following statements hold:\\
	\noindent $(i)$ Assume $q=g$ and $\alpha_{R}>1$. Then
	\begin{equation}\label{eq:14072022}
		\text{$\frac{A^*_R}{n}\to 1$\quad and\quad $\frac{A^*_B}{q}\to g_B(\kappa_{\bm g})+\alpha_B$,\qquad a.s.}
	\end{equation}
	\noindent$(ii)$ Assume $g\ll q\ll n$. Then, for any $\alpha_R>0$, 
	\[
	\text{$\frac{A^*_R}{n}\to 1$\quad and\quad $\frac{A^*_B}{n}\to 0$,\qquad a.s.}
	\]}
\end{Theorem}

The quantities $\kappa_{\bm g}$ and $g_B(\kappa_{\bm g}):=\lim_{y\uparrow \kappa_{\bm g}}g_B(y)$ are defined as follows.
\begin{Definition}\label{def:CPsimp}{\it (Cauchy problem).	
We denote by $\bm g(y) = (g_R(y),g_B(y))$ the maximal solution of the Cauchy problem:
\begin{equation}\label{simplified:CPcoupled}
	\bm g'(y)=\bm{\beta}(\bm g(y)),\quad\text{$y\in [0,\kappa_{\bm g})$,\quad $\bm g(0)=(0,0)$}
\end{equation}
where $\bm{\beta}:=(\beta_R,\beta_B)$.
}
{
It is worth noting that, as an immediate consequence of the celebrated Cauchy-Lipschitz theorem, Cauchy problem \eqref{simplified:CPcoupled} has a unique local solution. This is guaranteed because $\bm{\beta}(\cdot,\cdot)$ is Lipschitz on an open set containing $(0,0)$. This unique local solution can then be extended to its maximal domain.}
\end{Definition}
A more explicit characterization of $\kappa_{\bm g}$ and $g_B(\kappa_{\bm g})$ will be provided in Proposition \ref{le:coupledphi}.
{ In simple terms, Theorem \ref{prop:supercritical} indicates that, in the super-critical regime, the $R$-activation process spreads across nearly the entire graph. This effectively causes an "early stop" of the competing $B$-activation process. 
Specifically,
when $q = g$ and $\alpha_B\le 1<\alpha_R$, the value of $g_B(\kappa_{\bm g})$ is strictly less than $z_B$, which implies that,
in this case $\frac{A^*_B}{q}$ tends to 
a value strictly smaller than its counterpart without competition (as detailed in Remark \ref{BPRemark}).

Furthermore $\frac{A^*_B}{q}$ remains finite even when $\alpha_B>1$,  which is  particularly noteworthy  because  in the absence of competition the $B$-activation process would percolate almost the whole graph (again, see Remark \ref{BPRemark}).
Finally, when $g\ll q\ll n$, the final number of black nodes is of smaller order than $n$ for every value of $\alpha_B$.}

Note that, while in the sub-critical regime the activation process stops when $O_{\mathrm{a.s.}}(q)$ nodes are active, in the super-critical regime
almost all nodes become $R$-active (i.e., the final size of $R$-active nodes is $n-o_{\mathrm{a.s.}}(n)$).
The qualitative picture emerging from Theorems
\ref{thm:subcriticostoc} and \ref{prop:supercritical} is summarized
in Table~\ref{tab:phase-diagram}.
{
\begin{Remark}
To make the main ideas of the proofs more transparent, we separate the
core probabilistic arguments from the more technical auxiliary results.
Each major proof is preceded by a brief overview of its conceptual steps,
while lengthy auxiliary derivations are deferred to the appendices.	
\end{Remark}	
}

\subsection{Numerical illustration of the results} \label{sect:numerical}

For the purpose of numerical illustration of our results, we consider the case $r=2$, which allows closed-form solutions of the main quantities of interest.
We focus on the super-critical regime with $q=g$.
In this case, using results reported in Proposition \ref{le:coupledphi},~\footnote{For $q\ll p^{-1}$, $\beta_S$ depends on $x_S$ only and we write $\beta_S(x_S)$.}
$\kappa_{\bm g}:=\int_0^\infty \frac{\dd  x}{\beta_R(x)} < \infty$.
Specifically,  with $r=2$,  from 
\eqref{eq:beta0} we have 
\[
\beta_S(x_R,x_B)=\frac{(x_S+\alpha_S)^2}{4}-x_S, 
\]
and we get the closed-form expression:
\begin{equation}\label{eq:kappaS1r2}
	\kappa_{\bm g}:=\int_{0}^{\infty}\frac{\dd x}
	{\frac{(x+\alpha_R)^2}{4}-x} = \frac{2}{\sqrt{\alpha_R-1}}
	\left( \frac{\pi}{2}-\arctan\left(\frac{\alpha_R -2}{2 \sqrt{\alpha_R-1}}\right)\right).
\end{equation}
{ The constant $\kappa_{\bm g}$ is the finite explosion time
		of the limiting red trajectory at scale $q$.
		Since $\eta=1$ when $q=g$, as specified in \eqref{eta},
		this limiting trajectory is parameterized by unscaled
		physical time. Its explosion indicates the breakdown
		of the description at a fixed multiple of scale $q$;
		it does not mean that $n-o(n)$ vertices have already
		activated at physical time $\kappa_{\bm g}$.
		
		Indeed, for any fixed $t<\infty$, let $W_0(t)$ denote
		the number of initially white vertices whose clocks
		have not rung by time $t$. Then
		\[
		W_0(t)\overset{\mathrm L}{=}
		\operatorname{Bin}(n-a_R-a_B,e^{-t}),
		\]
		and binomial concentration, together with
		$a_R+a_B=o(n)$, gives
		\[
		\frac{W_0(t)}{n}\longrightarrow e^{-t},
		\qquad\text{a.s.}
		\]
		All these vertices are still white at time $t$.
		Consequently, almost percolation cannot occur at any
		fixed finite physical time.
		
		As $\alpha_R\downarrow1$, we have
		$\kappa_{\bm g}\to\infty$, reflecting the slowdown
		of the initial red activation dynamics near the
		critical threshold.}

As shown in Appendix  \ref{Appendix-CP}  (see \eqref{simplified:CP} , \eqref{eq:identityCP} and \eqref{sep-var}),  $g_B(\kappa_{\bm g})$ in the considered case satisfies the following equation: 
\begin{equation}\label{eq:mic}
	\int_0^{g_B(\kappa_{\bm g})} \frac{1}{\beta_B(y)} \,\dd  y = \kappa_{\bm g}.
\end{equation}
We distinguish two cases, depending on whether $\alpha_B$ is smaller or 
greater than 1.

\paragraph*{Case $\alpha_B < 1$}
In this case
$\beta_B(x)=\frac{(x+\alpha_B)^2}{4}-x$ 
has two zeros. The smallest one is at $z_B = 2-\alpha_B - 2\sqrt{1-\alpha_B}$
and the other one is at $w_B = 2-\alpha_B + 2\sqrt{1-\alpha_B}$. Note also that
$z_B \cdot w_B = \alpha_B^2$.
Conveniently, the integral in \eqref{eq:mic} admits the following closed-form expression:
\begin{equation}\label{eq:kappaS1r2bis}
	\int_{0}^{g_B(\kappa_{\bm g})}\frac{\dd y}{\beta_B(y)} =\int_{0}^{g_B(\kappa_{\bm g})}\frac{\dd y}
	{\frac{(y+\alpha_B)^2}{4}-y} = \frac{1}{\sqrt{1-\alpha_B}}
	\log\left|\frac{z_B(w_B-g_B(\kappa_{\bm g}))}{w_B(z_B-g_B(\kappa_{\bm g}))}\right|
\end{equation}
From \eqref{eq:mic} and \eqref{eq:kappaS1r2bis} one can compute $g_B(\kappa_{\bm g})$ explicitly.   Theorem \ref{prop:supercritical} then provides the asymptotic 
behavior of the (normalized) final number of $B$-active nodes in terms of  $g_B(\kappa_{\bm g})$:
\begin{equation*}
	\frac{A_B^*}{q} \rightarrow g_B(\kappa_{\bm g})+\alpha_B   = \frac{\alpha_B^2(\xi-1)}
	{(2-\alpha_B)(\xi-1)+2\sqrt{1-\alpha_B}(\xi+1)}+\alpha_B, \quad{\text{a.s.}}
\end{equation*}
where $$\xi=\xi(\alpha_R,\alpha_B):=e^{\kappa_{\bm g}\,\sqrt{1-\alpha_B}} $$

Note that the above quantity is strictly smaller than $\alpha_B + z_B$ for any 
$\xi$. As $\alpha_R \downarrow 1$, $\xi$ diverges to $\infty$, and 
we recover the well-known result of classical subcritical bootstrap percolation process with $r=2$, where the (normalized) final number of active nodes 
converges to $\alpha_B + z_B = 2-2\sqrt{1-\alpha_B}$.
Numerical results for different choices of $\alpha_R > 1 > \alpha_B$ are reported in
Fig. \ref{fig:fig1}.

\paragraph*{Case $\alpha_B > 1$}
In this case
\begin{equation}\label{eq:kappaS1r2tris}
	\int_0^{g_B(\kappa_{\bm g})} \frac{1}{\beta_B(y)}\,  \dd  y = 
	\frac{2}{\sqrt{\alpha_B-1}}
	\left( \arctan\left(\frac{g_B(\kappa_{\bm g}) + \alpha_B - 2}{2 \sqrt{\alpha_B-1}}\right)-\arctan \left(\frac{\alpha_B -2}{2 \sqrt{\alpha_B-1}}\right)\right).
\end{equation}

From \eqref{eq:mic} and \eqref{eq:kappaS1r2tris} one can compute $g_B(\kappa_{\bm g})$ explicitly also in this case. The (normalized) final number of $B$-active nodes is asymptotically estimated by
\begin{equation*}
	\frac{A_B^*}{q} \rightarrow 2+2\sqrt{\alpha_B-1} \tan(\xi'),\quad {\text{a.s.}}
\end{equation*}
where  
$$\xi'=\xi'(\alpha_R,\alpha_B):=
\arctan \left(\frac{\alpha_B -2}{2 \sqrt{\alpha_B-1}} \right)
+ \sqrt{\frac{\alpha_B-1}{\alpha_R-1}}
\left( \frac{\pi}{2}-\arctan\left(\frac{\alpha_R-2}{2 \sqrt{\alpha_R-1}}\right)\right).
$$
{
	The behavior near $(\alpha_R,\alpha_B)=(1,1)$ depends on how
	this point is approached. For example, along the paths
	\[
	\alpha_R=1+c_R\varepsilon,
	\qquad
	\alpha_B=1+c_B\varepsilon,
	\qquad
	c_R>c_B>0,
	\]
	with $c_R$ and $c_B$ fixed, the deterministic limiting
	final size satisfies
	\[
	g_B(\kappa_{\bm g})+\alpha_B\longrightarrow 2,
	\qquad \varepsilon\downarrow0.
	\]
	However, this conclusion does not hold for every approach
	to $(1,1)$ satisfying $\alpha_R>\alpha_B>1$. For instance,
	along the path
	\[
	\alpha_B=1+\varepsilon,
	\qquad
	\alpha_R=1+\varepsilon+\varepsilon^2,
	\]
	the explicit formula above gives
	\[
	g_B(\kappa_{\bm g})+\alpha_B
	\sim \frac{4}{\pi\sqrt{\varepsilon}}
	\longrightarrow+\infty.
	\]
	Here the limit $\varepsilon\downarrow0$ is taken after
	$n\to\infty$ for each fixed pair $(\alpha_R,\alpha_B)$.
	Thus, the limiting final-size formula does not have a
	path-independent joint limit as
	$\alpha_R,\alpha_B\downarrow1$.

}
{
One can also easily check that, for increasing values of $\alpha_R$, 
$A_B^*/q$ approaches $\alpha_B$,  meaning that the infection of $B$ nodes 
essentially does not evolve, being immediately stopped by the infection of $R$ nodes.
Instead, as $\alpha_R \downarrow \alpha_B$, $A_B^*/q$ diverges (note indeed
that in this case $\xi' \uparrow \frac{\pi}{2}$).  
Numerical results for different choices of $\alpha_R > \alpha_B > 1$ are reported in
Fig. \ref{fig:fig2}.

\sidebyside{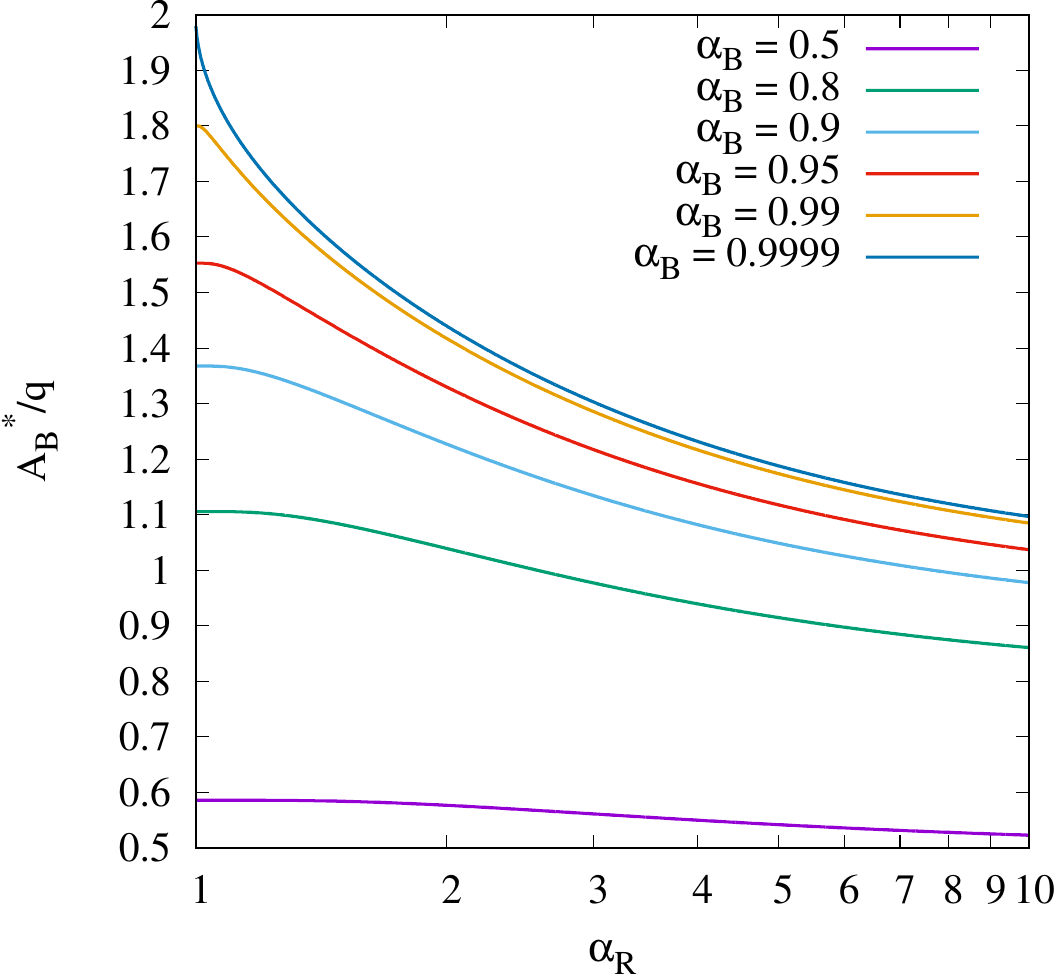}{ Case $q=g$, $\alpha_R > 1 > \alpha_B$: $A_B^*/q$ as function of $\alpha_R$, for different values 
	of $\alpha_B$.}{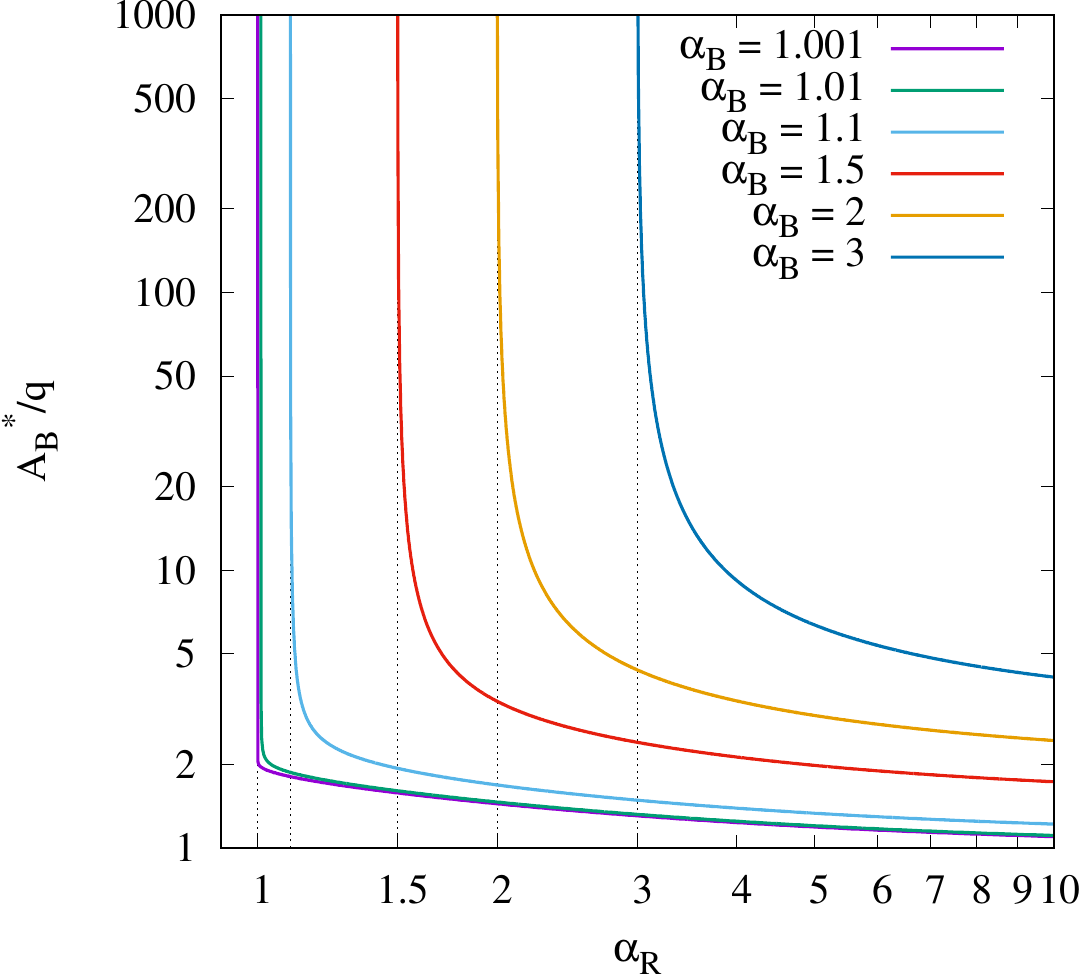}{Case $q=g$, $\alpha_R > \alpha_B > 1$: $A_B^*/q$ as function of $\alpha_R$, for different values of $\alpha_B $.}

\section{Process Representation and Preliminary Results}\label{Sect:PA}
\subsection{Definition of main variables and sets}\label{sect:def-var}

In this subsection we introduce the random quantities in terms of which we will
describe the dynamics of the competing bootstrap percolation processes on $G$. 

Let $\mathcal{V}_W\subset\mathcal{V}$ be the set of non-seed nodes and 
let $n_W:=|\mathcal V_W|=n-(a_R+a_B)$. To each node $v\in\mathcal{V}_W$,
we attach an independent, unit-rate Poisson process
(called Poisson clock),   
whose ordered points represent the successive wake-up times of node $v$. 
More formally, we define a collection $\{N_v'\}_{v\in\mathcal{V}_W}$ of independent Poisson processes on $[0,\infty)\times\mathcal{V}_W$ where each process $N_v'$ has mean measure $\dd t\delta_v(\dd \ell)$, where $\delta_v(\cdot)$ is the Dirac measure on $\mathcal{V}_W$ concentrated at $v\in\mathcal{V}_W$.
As it is well known,  the superposition
$$N':=\sum_{v\in\mathcal{V}_W}N_v'$$
is still a Poisson process on $[0,\infty)\times\mathcal{V}_W$
with mean measure $n_W\dd t\mathbb{U}(\dd v)$, where 
$\mathbb{U}$ is the uniform law on $\mathcal{V}_W$. We denote by $\{(T'_k,V'_k)\}_{k\in\mathbb N}$ the points of $N'$,  ordered by increasing time coordinate.  Here,  $T'_k$  represents the time of the $k$-th wake-up event and $V'_k$ the corresponding node.
For each $S\in\{R,B\}$,  we consider the $S$-activation point process $N_S$ on $[0,\infty)\times\mathcal{V}_W$: for any  $t>0$ and any $L\subseteq\mathcal{V}_W$, $N_S([0,t]\times L)$ 
counts the number of $S$-active nodes in $L\subseteq\mathcal{V}_W$ 
at time $t$. {In the following we refer to $t$ as physical time.}

Let $(T_k^S,V_k^S)$, $k\in\mathbb N$, denote the $k$-th point of $N_S$.
By construction, $T_k^S$ 
is the "activation time" of node $V_k^S$, 
i.e., the physical time at which node $V_k^S$ becomes 
$S$-active by taking color $S$.
A node $V'$ becomes $S$-active upon waking up at time $T'$ if and only if it is still white and fulfills the "threshold condition with respect to $S$".
Therefore, point process
$N_S$ can be constructed by thinning $\{(T'_k,V'_k)\}_{k\in\mathbb N}$ as follows: we retain only those couples
$(T'_k,V'_k)$, $k\in\mathbb N$, for which, at time $(T'_k)^{-}$, the node $V'_k$ is still white and
satisfies the "threshold condition with respect to $S$".

We set $N:=N_R+N_B$ and denote by $(T_k,V_k)$, $k\in\mathbb N$, the points of $N$.
Throughout this paper we refer to $N$ as the (global) activation process.
In the following we will use $N_S(t)$ and $N(t)$ 
as shorthand notation for $N_S([0,t]\times\mathcal{V}_W)$ and $N([0,t]\times\mathcal{V}_W)$, respectively.
Hereafter, we denote by $\mathcal{V}_S(t)\subset\mathcal{V}_W$, $t\geq 0$, the set of 
non-seed nodes which are $S$-active at {physical} time $t$, i.e.,
\[
\mathcal{V}_S(t):= \{V_k^S\}_{k:\,T_k^S\in [0,t]} \quad \text{with} \quad  \mathcal{V}_S(0):= \emptyset,
\]
and we denote by $\mathcal{V}_W(t)\subseteq\mathcal{V}_W$, $t\geq 0$, the set of 
nodes which are still $W$ at time $t$, i.e.,
\[
\mathcal{V}_W(t):=\mathcal{V}_W\setminus (\mathcal{V}_R(t)\cup \mathcal{V}_B(t)).
\]	

Let $\{E_i^{R,(v)}\}_{i\in\mathbb N}$, $\{E_i^{B,(v)}\}_{i\in\mathbb N}$, $v\in\mathcal{V}_W$, be two independent sequences
of independent and identically distributed random variables with Bernoulli distribution with mean $p$,  independent of $\{(T'_k,V'_k)\}_{k\in\mathbb N}$. The event 
$\{E_i^{S,(v)}=1\}$ ($\{E_i^{S,(v)}=0\}$) indicates the presence (absence) of an edge between node $v\in\mathcal V_W$ and the $i$-th  $S$-active node. 
We will often refer to random variables $E_i^{S,(v)}$ as $S$-marks.
We define the quantities:
\begin{equation}\label{eq:D}
D_{R}^{(v)}(t):=\sum_{i=1}^{N_R(t)+a_R}E_i^{R,(v)}\quad\text{and}\quad
D_{B}^{(v)}(t):=\sum_{i=1}^{N_B(t)+a_B}E_{i}^{B,(v)},\quad v\in\mathcal{V}_W.
\end{equation}
Specifically, $D_{S}^{(v)}(t)$ denotes the number of node $v$ neighbors with color $S$ at {physical} time $t$.
The sets of  \Ssupra{} and \supra{} nodes, at time $t$, 
are defined by
\begin{equation}\label{eq:Susceptible}
\mathcal{S}_{S}(t):=\{v\in\mathcal{V}_W:\,\,D_{S}^{(v)}(t)-D_{\overline{S}}^{(v)}(t)\geq r\}\quad\text{and}\quad\mathcal{S}(t):=\mathcal{S}_{R}(t)\cup\mathcal{S}_{B}(t),
\end{equation}	
respectively.
Note that all previously introduced variables and sets can be defined at {physical} time $t^-$ by replacing $[0,t]$ with $[0,t)$.	
The final number of active nodes is given by
\[ 
A^*:=A^*_R+A^*_B,\quad\text{where}\quad A_S^*:=N_S([0,\infty))+a_S.
\]
Recalling that the epidemic process naturally stops as soon as no more jointly \supra{} and white nodes can be found, we can define the random time-index at which the process stops, $K^*$, as:
\[
{K^*}:=\min\{k\in\mathbb N\cup\{0\}:\,\,\mathcal{S}(T_{k})\cap\mathcal{V}_W(T_{k})=\emptyset\},\quad\text{$T_0:=0$.}
\]
Consequently,  by construction,  the global activation process ceases at time $T_{K^*}$,  and we have
\begin{equation}\label{Kappa-Astar}
A^*={K^*}+a_R+a_B.
\end{equation}
For the moment, we conventionally set
$T_{K^*+1}:=\infty$, and note that, on the event $\{ K^*<k\}$, we have $T_{k}=\infty$. 
It is worth mentioning that, for technical reasons, in Section \ref{sect:extensionbeyondT} we will artificially extend the activation process $N$ beyond $T_{{K^*}}$, by redefining the times $T_k$ on the event $\{ K^*<k\}$.
We emphasize that this extension does not alter  dynamics of the process before time $T_{{K^*}}$.

Finally, we remark that, without loss of generality, throughout this paper, we assume that the random graphs $G(n,p_n)$ and the dynamical   processes evolving on them are independent for different values of $n$.

\begin{Remark}
During the evolution of the activation process, an edge $\{v, w\}\in\mathcal E$ is unveiled potentially twice (i.e., when $v$ becomes active and when $w$ becomes active). As in the classical bootstrap percolation process studied in \cite{JLTV}, this
has no effect on the dynamics of the competing bootstrap percolation processes.	
Indeed, if $v$ activates before $w$, then any mark potentially added to $v$ when $w$ activates has no impact on the system evolution.
\end{Remark}

\subsection{Discrete time notation}\label{discrete-time_notation}

To study the evolution of the system, 
it is convenient to introduce some discrete time notation.
For time-index $k \in \mathbb{N}\cup\{0\}$,
we set			
\begin{align*}
\widehat{T}_k&:=\min\{T_{k}, T_{K^*}\},\quad N_S[k]:=N_S(\widehat{T}_k), \quad  \mathcal{S}_S[k]:=\mathcal{S}_S(\widehat{T}_k), \\
\mathcal{V}_W[k]&:=\mathcal{V}_W(\widehat{T}_k)\quad\text{and}\quad D^{(v)}_S[k]:={D^{(v)}_S(\widehat{T}_k)}.
\end{align*}
Moreover, we define:
\begin{equation}\label{eq:defU}
U^R_{k+1} :=\frac{|\mathcal{V}_W[k]\cap\mathcal{S}_{R}[k]|} {|\mathcal{V}_W[k]\cap\mathcal{S}[k]|},
\quad
U^B_{k+1}=1-U^R_{k+1},\quad k\in\mathbb{N}\cup\{0\}
\end{equation}
where we conventionally set $0/0:=1/2$. 

Given $\mathcal{V}_W[k]$,  $\mathcal{S}_{R}[k]$ and $\mathcal{S}_{B}[k]$, $U_{k+1}^R$ ($U_{k+1}^B$) is defined as 
the probability that a node, taken uniformly at random in $\mathcal{V}_W[k]\cap\mathcal{S}[k]$, is \Rsupra{}  (\Bsupra). Building on the properties of Poisson processes (see Remark \ref{remark_Poisson}), this can be understood as the conditional probability that the $(k+1)$-th node (excluding seeds) 
to activate is assigned color $R$ ($B$). Proposition \ref{M-U} will clarify this point further.

Finally, we note that by construction it holds:
\[
\mathcal{V}_W[k]=\mathcal{V}_W\setminus(\mathcal{V}_R[k]\cup \mathcal{V}_B[k]), \qquad  
| \mathcal{S}_{S}[k] \cap \mathcal{V}_S[k] |=| \mathcal{V}_S[k]| - |\mathcal{V}_S[k]\cap ( \mathcal{V}_W\setminus  \mathcal{S}_{S}[k]) |
\]
and 
\begin{equation*}
\mathcal{V}_{S}[k]\subseteq \{  v:  D_S^{(v)}[k]\geq r  \},\label{03052025}
\end{equation*}
and so,
for all $S\in\{R,B\}$ and $k\in\mathbb N\cup\{0\}$, we have
\begin{align}
|\mathcal{V}_W[k]\cap\mathcal{S}_{S}[k]| =  |\mathcal{S}_S[k]|-N_S[k]
&+|(\mathcal{V}_W\setminus\mathcal{S}_S[k])\cap \mathcal{V}_S[k]\cap\{ v:  D_S^{(v)}[k]\geq r   \}|  \nonumber\\
&-|\mathcal{S}_S[k]\cap \mathcal{V}_{\overline{S}}[k]\cap\{ v:  D_{\overline{S}}^{(v)}[k]\geq r   \}|:=Q^S_{k+1}. 
 \label{eq:QSnew} 
\end{align}
{
\begin{Remark}\label{anticipated}
Even though the expression for $Q_{k+1}^S$ looks complex,
its asymptotic behavior is tractable. This 
because the final two terms in 
\eqref{eq:QSnew} are negligible compared to the first two, as shown by Proposition \ref{cor:BC}.    
This tractability will also prove useful when we extend the process beyond $T_{K^*}$
in the next section.
\end{Remark}
}
\subsection{Prolonging process $N$ beyond $T_{{K^*}}$}\label{sect:extensionbeyondT}
Even though no activation events occur after time  $T_{{K^*}}$,  i.e., $N({(T_{K^*}},\infty)\times\mathcal{V}_W)=0$, 
it  is convenient, for analytical purposes, to extend the point process $N$
beyond $T_{{K^*}}$  by allowing the activation 
of nodes that are not \supra{}. 
This extension facilitates
 the analysis without altering the process dynamics up to time $T_{{K^*}}$.
 {
 To avoid ambiguity, before introducing the artificial extension, we
 reserve the notation $N^\circ_S$ and $N^\circ=N^\circ_R+N^\circ_B$
 for the original physical processes. More precisely, for
 $S\in\{R,B\}$, we set
 \[
 N^\circ_S(t):=N_S(t\wedge T_{K^*}),
 \qquad
 N^\circ(t):=N^\circ_R(t)+N^\circ_B(t).
 \]
 We also define
 \[
 T^{S,\circ}_j
 :=
 \inf\{t\geq0:N^\circ_S(t)\geq j\},
 \qquad
 T^\circ_j
 :=
 \inf\{t\geq0:N^\circ(t)\geq j\},
 \]
 with the convention $\inf\emptyset:=+\infty$.
 The quantities $K^*$ and $A_S^*$ always refer to the original
 physical process.
From this point onward, unless the superscript $\circ$ is explicitly
displayed, $N_S$ and $N$ denote the processes artificially extended
beyond $T_{K^*}$. We retain the notation
$\{(T_k^S,V_k^S)\}_{k\geq1}$ and
$\{(T_k,V_k)\}_{k\geq1}$ for the supports of these extended
processes.
By construction, the original and extended processes coincide up to
the original stopping time:
\[
N_S(t)=N^\circ_S(t),
\qquad 0\leq t\leq T_{K^*},\quad S\in\{R,B\}.
\]
Equivalently, all the corresponding discrete-time variables coincide
on the event $\{k\leq K^*\}$.

On the event $\{K^*<k\}$, the variables introduced in
Section~\ref{discrete-time_notation}, including
$N_S[k]$, $\mathcal V_W[k]$, $\mathcal S_S[k]$,
$D_S^{(v)}[k]$, $Q^S_{k+1}$, and $U^S_{k+1}$,
refer to the extended process and may differ from their original
physical counterparts.
}
Points $\{(T_{k},V_{k})\}_{k>K^*}$ are obtained
by thinning the point process 
$\{(T'_{k'},V'_{k'}) \}_{k'>K^*}$
retaining only those couples $(T'_{k'},V'_{k'})$ such that $V'_{k'}$ is still $W$. 
Then we determine $\{(T_{k}^R,V_{k}^R)\}_{k>K^*}$ ($\{(T_{k}^B,V_{k}^B)\}_{k>K^*}$) by randomly assigning
color $R$ ($B$) to each node $V_{k}$, on $\{ K^*< k\}$.
More precisely,  setting

\begin{equation}
U^S_{K^*+1}:= \frac{1}{2}, \qquad  U^S_{k+1}:= \frac{|Q^S_{k+1}|}{|Q^R_{k+1}|+|Q^B_{k+1}|}  \quad \text {on } \{ K^*<k\},
\label{eq:redefineU}
\end{equation}
where, i)  consistently with \eqref{eq:defU}, we adopt the convention
$0/0:=1/2$; ii)  $Q^S_{k+1}$ is still defined as in \eqref{eq:QSnew}, then
for every $u\in \supp(U^R_{k})$,  
conditionally
on $\{U^R_{k}=u, K^*< k\}$,  color $R$ is assigned to $V_k$ with probability $u$. This can be achieved,  as
explained in more detail in Section \ref{sect-markov},  sampling a uniformly distributed random variable with support $(0,1)$ and comparing it with $U^R_{k}$.

We wish to emphasize that this extension implies a redefinition of all the random variables in Section \ref{discrete-time_notation} on the set $\{ K^*< k\}$.
Notably, on $\{ K^*<k\}$ the activation of  \nsupra nodes invalidates the identity in \eqref{eq:QSnew}, potentially resulting in negative values for $Q^S_{k+1}$.
This necessitates the absolute values in the rightmost term of \eqref{eq:redefineU}.
A simple computation shows that
definition \eqref{eq:redefineU} aligns with \eqref{eq:defU} for $\{K^* \ge k\}$.
Note that the extended process naturally stops at $k=n_W$; for $k> n_W$, we set $T_k:=\infty$, $N[k]:=N[n_W]$, $\mathcal{S}[k]:=\mathcal{S}[n_W]$, etc.

Moreover, since the sum between the second and the third term in \eqref{eq:QSnew} is non-positive and the absolute value of the fourth term in
\eqref{eq:QSnew} is at most $N_{\overline{S}}[k]$, the number of \Ssupra{} nodes satisfies:
\begin{equation}\label{180523-1}
|\mathcal{S}_{S}[k]| -k\leq Q^S_{k+1} \le |\mathcal{S}_{S}[k]|,\qquad    k \in \mathbb{N}\cup \{0\}.
\end{equation}

We conclude this section observing that by \eqref{eq:Susceptible} we have
\begin{equation} \label{eq:susceptible2}
|\mathcal{S}_S[k]|=\sum_{v\in\mathcal{V}_W} \ind_{\{D^{(v)}_S[k]-D_{\overline{S}}^{(v)} [k]\geq r\}}, \qquad  k \in \mathbb{N}\cup\{0\},
\end{equation}

{
For later use, it is convenient to introduce the following auxiliary variables.
Given a  vector $\bm{k}=(k_R,k_B)\in(\mathbb{N}\cup\{0\})^2$ and
$S\in\{R,B\}$, define
\begin{equation}
	\widehat D^{(v)}_S(\bm{k})
	:=
	\sum_{i=1}^{k_S+a_S} E_i^{S,(v)},
	\qquad v\in \mathcal V_W,
\end{equation}
and
\begin{equation}
	\widehat{\mathcal S}_S(\bm{k})	:=
	\left\{v\in \mathcal V_W:\widehat D^{(v)}_S(\bm{k})-
	\widehat D^{(v)}_{\overline S}(\bm{k})	\ge r	\right\}
\end{equation}
in the following we will refer to variables $	\widehat{\mathcal S}_S(\bm{k})$ and $\widehat D^{(v)}_S(\bm{k})$ as 
fixed-configuration \Ssupra{} set and fixed-configuration $S$ mark count for node $v$,  respectively.
The introduction of these auxiliary variables  allows us to establish concentration bounds uniformly over deterministic configurations and then evaluate them along the random trajectory $\bm N[k]$.

Define, moreover
\begin{equation}
	\label{pi-static}
	\pi_S(\bm{k})	:=	\mathbb{P}\left(  \text{Bin}(k_S+a_S,p)	-
	\text{Bin}(k_{\overline S}+a_{\overline S},p)	\ge r	\right),
\end{equation}
where the two binomial random variables appearing in the
right-hand side are independent.
With a slight abuse of notation, for
$\bm{k}=(k_R,k_B)$ we also write
$\pi_S(k_R,k_B):=\pi_S(\bm{k})$.
\begin{Lemma}
	\label{lemma:deterministic-marks}
	For every 
	$\bm{k}=(k_R,k_B)\in(\mathbb{N}\cup\{0\})^2$ and
	$S\in\{R,B\}$,
	\begin{equation}
		\label{static-binomial}
		\left|\widehat{\mathcal S}_S(\bm{k})\right|
		\overset{\mathrm L}{=}  \text{Bin}	\left(n_W,\pi_S(\bm{k})\right).
	\end{equation}
	Moreover, for every $k\le n_W$,
	\begin{equation}
		\label{static-representation}
		\mathcal S_S[k]
		=
		\widehat{\mathcal S}_S\bigl(\bm N[k]\bigr)
		\qquad\text{a.s.},
	\end{equation}
	where
	$\bm N[k]=(N_R[k],N_B[k])$.
	
	Finally, the family	$\widehat{\mathcal S}_S(\bm{k})$ is monotone in the
	following pathwise sense. If
	$\bm{k}=(k_S,k_{\overline S})$ and
	$\bm{\ell}=(\ell_S,\ell_{\overline S})$ satisfy
	\begin{equation}
		k_S\ge \ell_S,
		\qquad
		k_{{\overline S}}\le \ell_{{\overline S}},
	\end{equation}
	then
	\begin{equation}
		\label{static-monotonicity}
		\widehat{\mathcal S}_S(\bm{k})
		\supseteq
		\widehat{\mathcal S}_S(\bm{\ell})
		\qquad\text{a.s.}
	\end{equation}
\end{Lemma}
A brief  proof of Lemma \ref{lemma:deterministic-marks} is provided in Appendix \ref{App-new}.

\begin{Remark}
	\label{rem:static-uniform}
	Lemma~\ref{lemma:deterministic-marks} allows us to avoid conditioning
	on the random vector $\bm N[k]$.
	
	Indeed, let ${\mathbb A}_n$ be any deterministic finite subset of
	$(\mathbb N\cup\{0\})^2$. For any $x\ge 0$, 
	the union bound gives
	\begin{align}
		&\mathbb P\left(
		\exists\,\bm k\in{\mathbb A}_n: 	\left||\widehat{\mathcal S}_S(\bm k)|	-		n_W\pi_S(\bm k)\right|	>x 
			\right)		\nonumber\\
	&\qquad\le 	\sum_{\bm k\in\mathbb{A}_n}
		\mathbb P\left(
		\left|
		\text{Bin}(n_W,\pi_S(\bm k))
		-
		n_W\pi_S(\bm k)
		\right|
		>x 
		\right).
		\label{uniform-static-bound}
	\end{align}
	Consequently, any concentration estimate proved uniformly over
	$\bm k\in{\mathbb A}_n$ can subsequently be evaluated at the
	random vector $\bm N[k]$, using
	\eqref{static-representation}, without requiring any independence
	between $\bm N[k]$ and the underlying mark process.
\end{Remark}
}

\begin{Remark}  
We have the freedom to choose any form for the quantity $U_{k}^S$, for $k>K^*$.   
Indeed, it has no impact on the process dynamics up to $K^*$.
The selected form of $U_{k}^S$ for $k>K^*$ simplifies the analysis considerably, even if it appears somewhat artificial. This comes from the fact that the asymptotic behavior of
$Q_{k+1}^S$ is actually easy to characterize, as anticipated in Remark \ref{anticipated}.
\end{Remark}

\subsection{Markovianity of the system}\label{sect-markov}

The next proposition states the Markovianity of the system. Its proof is based on standard computations and therefore it is omitted. We refer the reader to
\cite{B2, norris} for any unexplained notion concerning Markov chains.

\begin{Proposition}\label{MarkovChain}
{\it
The stochastic process
\[
\bm Z=	\{{\bm Z}(t)\}_{t\geq 0}:= \{((\ind_{\{v\in \mathcal{V}_R(t) \}}, \ind_{\{ v\in  \mathcal{V}_B(t)\} },  D_{R}^{(v)}(t) , D_{B}^{(v)}(t))_{v\in\mathcal{V}_W}, \ind_{\{T_{K^*}\le t\}}) \}_{t\geq 0}	
\]
is a regular-jump, continuous time, homogeneous Markov chain, i.e., a continuous time homogeneous Markov chain such that, for almost all $\omega$, $|\text{Disc}(\omega)\cap [0,c]|<\infty$, for any $c\geq 0$. Here $\text{Disc}(\omega)$ denotes the set of discontinuity points of the mapping $t\mapsto\bm{Z}(t,\omega)$.

Let 
$$\mathbb{S}\subseteq (\{0,1\}\times\{0,1\}\times\{0,\ldots,n_W\}\times\{0,\ldots,n_W\} )^{|\mathcal{V}_W|}
\times\{0,1\}$$ denote the state space of $\bm Z$ and	 
\begin{equation}\label{201023-2}
{R(\bm z)}:=\lim_{h\to 0}\frac{1-\mathbb{P}(\bm{Z}(h)=\bm{z}\,|\,\bm{Z}(0)=\bm{z} )} {h},\quad\bm z\in\mathbb{S}
\end{equation}
the diagonal elements of the transition-rate matrix.
}
\end{Proposition}
Since for any $t\ge 0$, $\mathcal{V}_W(t)$, $\mathcal{S}_S(t)$ and $N(t)$ are $\sigma\{\bm Z(t)\}$-measurable random variables,
with a slight abuse of notation, we conveniently denote
them with the symbols $\mathcal{V}_W(\bm Z(t))$, $\mathcal{S}_S(\bm Z(t))$ and $N(\bm Z(t))$, respectively.
From the properties of Poisson processes and their thinnings it immediately follows that 
\begin{equation}\label{201023-3}
	R({\bm z})=| \mathcal{S}(\bm z) \cap \mathcal{V}_W(\bm z)| \ind_{\{\bm z^{(E)}=0\}}+(n_W- N(\bm z) )  \ind_{\{\bm z^{(E)}=1\}},
\end{equation}
where $\bm z^{(E)}$ denotes the last component of vector $\bm z$ (which is equal to $0$ or $1$).
Note, indeed, that at time $t$ only the jointly \supra{}
and white nodes,  i.e. ,  nodes in $\mathcal{S}(\bm Z(t)) \cap \mathcal{V}_W(\bm Z(t))$,  are enabled for activation if $t<T_{K^*}$, while the entire set of white nodes, whose cardinality is $n_W- N(\bm Z(t))$,  
is enabled for activation if $t\ge T_{K^*}$.

The sequence of 
transition times of $\bm Z$ coincides with the sequence of activation times $\{T_k\}_{k\geq 0}$, $T_0:=0$, of the nodes.
Let $\mathcal{F}_t^{\bm Z}:=\sigma\{\bm{Z}(s):\,\,s\leq t\}$ be the natural filtration of the Markov chain $\bm{Z}$ and let $\{\bm{Z}_k\}_{k\in\mathbb N\cup\{0\}}$ be the embedded chain defined by $\bm Z_k=\bm{Z}(T_k)$. 
We have
\[
\bm{Z}_{n_W}\in\mathbb{S}_0:=\{\bm{z}\in\mathbb{S}:\,\,R({\bm z})=0\},
\]
and $\{\bm{Z}_k\}_{0\leq k< n_W}\in \mathbb{S}\setminus\mathbb{S}_0=\{\bm{z}\in\mathbb{S}:\,\,R({\bm z})>0\}$. 
Moreover, given $\{\bm{Z}_k\}_{0\leq k<n_W}$, sojourn times $\{W_k\}_{0\leq k<n_W}$,
$W_k:=T_{k+1}-T_k$ are independent and $W_k$ is exponentially distributed with mean $\frac{1}{R({\bm Z_k})}$ (see  \eqref{131123}).  

Since all the random variables defined in Section \ref{discrete-time_notation}, i.e.
\[
\text{$N_S[k]$, $\mathcal{V}_W[k]$, $\mathcal{S}[k]$, $\mathcal{S}_S[k]$, $U^S_{k+1}$, and $Q_{k+1}^S$} 
\]
are $\sigma(\bm Z_k)$-measurable,  with a slight abuse of notation,  they will be conveniently denoted by 
\[
\text{$N_S(\bm Z_k)$, $\mathcal{V}_W(\bm Z_k)$, $\mathcal{S}(\bm Z_k)$, $\mathcal{S}_S(\bm Z_k)$, $U^S(\bm Z_k)$ and $Q^S(\bm Z_k)$,}
\]
respectively.

We define binary random variables:
\begin{equation}\label{26102022}
M_{k+1}^S:=N_S[k+1]-N_S[k],\quad\text{$k\geq 0$, $S\in\{R,B\}$.}
\end{equation}
$M_{k+1}^S$ indicates whether $V_{k+1}$ gets color $S$.
Clearly, $M_{k+1}^S\in \sigma\{\bm Z_k, \bm Z_{k+1}\}$. 
Moreover, recalling that on $\{K^*> k\}$, $V_{k+1}$ receives color $R$ if and only if 
$V_{k+1} \in\mathcal{V}_W[k]\cap \mathcal{S}_R[k]$, while on $\{ K^*\leq k \}$, a color is
randomly assigned to $V_{k+1}$ as detailed in Section \ref{sect:extensionbeyondT}, 
we can write:
\begin{equation}\label{eq:12062025due}
M_{k+1}^R= \ind_{\{V_{k+1}\in\mathcal{V}_W[k]\cap \mathcal{S}_R[k]\}}\ind_{\{K^*>k\}}+ 
\ind_{\{ K^*\le k\}}\ind_{\{ L_{k+1}<U^R_{k+1}\}}\quad\text{and}\quad M_{k+1}^B=1-M_{k+1}^R, 
\end{equation}
where $\{L_{k+1}\}_{k\geq 0}$ is a sequence of random variables uniformly distributed on $(0,1)$, and
such that $L_{k+1}$ is independent of $\mathcal{H}_k:=\sigma\{\bm{Z}_h:\,\,0 \le h \le k\}$, $k\in\mathbb{N}\cup\{0\}$.
Observe indeed that 
for any $u \in \supp(U^R_{k+1})$ we have $\ind_{\{ L_{k+1}<U^R_{k+1}\}} \mid
\{ U^R_{k+1}=u \}\overset{\mathrm L}{=} \text{Be}(u)$.

\begin{Proposition}\label{M-U}
For $S\in\{R,B\}$ and $k < n_W$, we have
\begin{align}
\mathbb{P}(M_{k+1}^S=1\mid \mathcal{H}_{k}) &= \mathbb{P}(M_{k+1}^S=1\mid \bm Z_{k} )
=\mathbb{E}[M_{k+1}^S\mid \bm Z_{k}] \nonumber \\
&=\mathbb{E}[\ind_{\{V_{k+1}\in\mathcal{V}_W[k]\cap \mathcal{S}_S[k]\}} \mid \bm Z_k]
\ind_{ \{ K^*> k \}  } + 
\mathbb{E}[ \ind_{\{ L_{k+1}<U^S_{k+1}\} }\mid \bm Z_k ] \ind_{\{ K^*\le k  \} } \nonumber \\ 
&=U_{k+1}^S.\label{eq:12062025uno}
\end{align}		
\end{Proposition}
The first equality in \eqref{eq:12062025uno} is a consequence of the Markovianity of $\{\bm Z_k\}$,  
the third equality is a consequence of \eqref{eq:12062025due} and the fact that $\{K^*\le  k  \} \in \sigma(\bm Z_k)$,
while the last equality follows from
the fact that, given $\bm Z_{k}$, provided that $K^*>k$,
$V_{k+1}$ is uniformly selected from $\mathcal{V}_W({\bm Z}_k)\cap \mathcal{S}({\bm Z}_k) \supseteq 
\mathcal{V}_W({\bm Z}_k)\cap \mathcal{S}_S({\bm Z}_k)= \mathcal{V}_W[k]\cap \mathcal{S}_S[k]$
(see Remark \ref{remark_Poisson}).
	
Proposition \ref{M-U} formally states that $U^S_{k+1}=U^S(\bm Z_k)$ can be interpreted as the conditional probability, given $\bm Z_k$, that color $S$ is assigned to $V_{k+1}$.
Finally, to avoid interrupting the main paper flow, we have moved two standard consequences of Markovianity, along with their straightforward proofs, to Appendix \ref{appendix-MC}. These results, which will be referenced in Theorem \ref{propVk} and Theorem \ref{thm:Lemma-new-3}, are best read when specifically invoked.

\subsection{Brief overview of the proofs of our main results}

Theorems \ref{thm:subcriticostoc} and \ref{prop:supercritical}, our main results, follow directly from intermediate findings detailed in Sections \ref{sub:Cauchy} and \ref{sub:scaling}.
As a guide to the reader, we briefly describe, at a high level, the strategy of the proofs. First, we analyze the activation process on time-scale $q$, i.e., we analyze the asymptotic behavior of $\bm N[\lfloor x q \rfloor ]/q$ for bounded values of $x$ (see Section \ref{sub:Cauchy}).
The main result on time-scale $q$
is provided by Theorem \ref{prop:equadiff}, which shows that a suitably regularized version of the trajectories $x\mapsto \bm N[\lfloor x q \rfloor ]/q$ 
converges almost surely to the (deterministic) solution of the Cauchy problem stated in Definition \ref{def:CP}.
To prove the convergence of such trajectories, we proceed as follows.
Exploiting the Ascoli-Arzelà theorem, we show that a subsequence of trajectories converges uniformly to a limiting function, almost surely.
Then, we provide sufficiently tight upper and lower bounds for the incremental ratio of the trajectories within a neighborhood of a fixed point. By doing so we show that the limiting trajectory is differentiable and indeed that it is the solution of the Cauchy problem formulated in Definition \ref{def:CP}. 
As a consequence, the uniqueness of the solution to the Cauchy problem
allows us to show that the whole sequence of trajectories converges pointwise to the limiting trajectory, almost surely.
Finally, the regularity of the trajectories enables us to upgrade the pointwise convergence to uniform convergence. Theorem \ref{propVk} complements this result by showing that normalized versions of 
both $T_{\lfloor x q \rfloor}$ and $T^S_{\lfloor x q \rfloor}$ (with a suitable $x$) converge almost surely to deterministic quantities. {This is accomplished by constructing upper and lower bounds for $T_{\lfloor x q \rfloor}$ and $T^S_{\lfloor x q \rfloor}$ through two appropriately defined sums of independent and exponentially distributed random variables. We subsequently show that these sums exhibit sufficient concentration around their means, and ultimately, that the means of these bounds are arbitrarily close.}

When the activation processes of the nodes do not stop at time-scale $q$, (i.e., in the super-critical regime)
we extend our analysis  to time-scales larger than $q$ (see Section \ref{sub:scaling}).   
In this case, an analysis of the properties of the solutions of the Cauchy problem \eqref{eq:CP} 
reveals that the ratio 
$N_B[\lfloor xq \rfloor]/N[\lfloor xq \rfloor]$ becomes arbitrarily small as $x$ grows large. 
The analysis at time-scales $q' \gg q$
hinges on the observation that the number $|\mathcal S_S(t)|$ of \Ssupra{} nodes
is sufficiently concentrated around its average. This average,  in turn, 
depends superlinearly on the number of active nodes $N_S(t)$. As a result, and as demonstrated in 
Theorems \ref{lemma-nuovo}, \ref{thm:Lemma-new-2} and \ref{thm:Lemma-new-3},
the ratio between the rates at which the two competing activation processes evolve tends quickly to infinity.
This allows the advantaged $R$-process to percolate before the competing $B$-process has managed to activate a non-negligible
fraction of nodes. In particular, for the case $q=g$ we can show that $A_B^*=O_{\mathrm{a.s.}}(g)$. This  claim is proved in two steps: first, we analyze the dynamics of an auxiliary process, the stopped process, where the $R$-activation process is stopped at a given point and only the $B$-activation process is allowed to continue; second, we infer the properties of the original process by exploiting simple coupling inequalities (see \eqref{stop-CB}).

\section{Analysis at time-scale $q$: main results} \label{Sect:q-timescalle-res}
\label{sub:Cauchy}

In this section, we report the main findings of our analysis of the
activation processes $N_S$, $S\in\{R,B\}$, in the regime
$N=\Theta_{\mathrm{a.s.}}(q)$, meaning that the total number of
activations is almost surely of the same order as the number of seeds.
As before, we assume, without loss of generality, that
$\alpha_R>\alpha_B$, and that conditions~\eqref{eq:bootcond}
and~\eqref{eq:trivial} hold.
The proofs of the results stated in this section are given in Section \ref{sec:V}.

We begin by introducing the linear interpolation $\widetilde{\bm N}(xq)=(\widetilde N_R(xq), \widetilde N_B(xq))$,  defined for $x\geq 0$,  as follows:
\begin{equation}\label{Ntildedef}
\widetilde{N}_S(xq):=N_S\big[\lfloor xq\rfloor\big]+(xq-\lfloor xq \rfloor)\left(N_S\big[\lceil xq\rceil \big]-N_S\big[\lfloor xq\rfloor\big]\right),
\end{equation}
and the sequence $\{\bm F_n(x)\}_{n\in\mathbb N}$ , where
\[
\bm F_{n}(x):=(F_{R,n}(x),F_{B,n}(x))\quad\text{with}\quad 
F_{S,n}(x):=\frac{\widetilde N_S(x q_n)}{q_n}.
\]
As usual, when no confusion arises, we will  drop the $n$ subscript from  $\bm F_n$ and $F_{S,n}$.  It turns out  that $\bm F$ converges to a vector-valued  function $\bm f$,
which is the solution of the following  Cauchy problem.
\begin{Definition}\label{def:CP}(Cauchy problem).
We denote by $	\bm f(x) = (f_R(x),f_B(x))$ the unique maximal solution of the Cauchy problem
\begin{equation}\label{eq:CP}
	\bm f'(x)=\frac{\bm{\beta}(\bm f(x))}{\beta_R(\bm f(x))+\beta_{B}(\bm f(x))},
	\quad x\in (0,\kappa_{\bm f}),\quad\text{$\bm f(0)=(0,0)$,}
\end{equation}
with $\bm{\beta}(\bm x):=
\bm{\beta}(x_R,x_B):=(\beta_R(x_R,x_B),\beta_B(x_R,x_B))$  as in \eqref{eq:beta0}. 
\end{Definition}
This is formalized by the following theorem.

\begin{Theorem}\label{prop:equadiff}
	For every $\kappa<\kappa_{\bm f}$,  we have
	\begin{equation}\label{eq:fluidlimit}
		\sup_{x\in  [0,\kappa] } \|\bm F(x)-\bm f(x)\|\to 0,\quad\text{a.s.}
	\end{equation}	
\end{Theorem}
As an immediate consequence of this theorem, we obtain the following corollary.
\begin{Corollary}\label{cor:190823}
For every $\kappa<\kappa_{\bm f}$ and $S\in\{R,B\}$,  we have
\begin{equation}\label{180823-2}
\lim \frac{\widetilde N_S(\kappa q)}{q}= f_S(\kappa), \quad\text{a.s.}
\end{equation}
\end{Corollary}
\subsection{On the solution of the Cauchy problem \eqref{eq:CP} }\label{sub:Cauchyproperties}

In this section, we summarize the key properties of the solution to the Cauchy problem \eqref{eq:CP} that are relevant to our main proofs.  A more detailed analysis of this solution,  including its connection to the solution of the simplified coupled problem \eqref{simplified:CPcoupled},  is provided in Appendix \ref{Appendix-CP}.

Recalling  Remark \ref{re:29Ott} and the fact 
that $ {\bm g}$ is the maximal solution of the Cauchy problem \eqref{simplified:CPcoupled},
we now state the following proposition.

\begin{Proposition} \label{le:coupledphi}	
The table below  reports the  values of $\kappa_{\bm f}$, 	 $\lim_{ x\uparrow \kappa_{\bm f}} f_R(x)$ and
 $\lim_{ x\uparrow \kappa_{\bm f}} f_B(x)$  for the different parameter regimes. It also provides explicit
 expressions for $f_R(x) $ and $f_B(x)$ when  
$p^{-1}\ll q \ll n$:
\begin{center}
\begin{tabular}{|c|c|c|c|c|c|c|}
\hline
    Case & Parameters  &$\kappa_{\bm f}$ & $\lim_{ x\uparrow \kappa_{\bm f}} f_R(x)$     &  $\lim_{ x\uparrow \kappa_{\bm f}} f_B(x)$   & $f_R(x)$ &  $f_B(x)$ \\
    \hline 	
$(i)$  &$q=g$ and  $\alpha_R<1$ & $z_R+z_B$ &$z_R$ &$z_B$  &-& - \\
$(ii)$  &  $q=g$ and $\alpha_R>1$ &$+\infty$  & $+\infty$  &$g_B(\kappa_{\bm g} )$  &-& -\\
$(iii)$  & $g\ll q\ll  p^{-1}$                   & $+\infty $& $+\infty$ & $g_B(\kappa_{\bm g} )$  &-&- \\
$(iv)$  & $q=p^{-1}$   & $+\infty$ & $+\infty$  &  $\overline{f}_B $ &- &- \\
$(v)$  &  $p^{-1}\ll q \ll n$  & $+\infty$ &$+\infty$ & 0 &  x& 0\\
\hline
\end{tabular}
\end{center}
Here $\kappa_{\bm g}:=\int_{0}^{\infty}\frac{\dd  y}{\beta_R(y)}<\infty$, $g_B(\kappa_{\bm g}):= \lim_{y\uparrow \kappa_{\bm g}} g_B(y)$ and  $\overline{f}_B$ is a suitable  strictly positive constant.  For case $(ii)$,  if $\alpha_B<1$,  then it follows that $g_B(\kappa_{\bm g})<z_B$. 
\end{Proposition}

{
	\begin{Remark}\label{re:70923}
	Note that if $q\ll p^{-1}$,  then $\beta_S(x_R,x_B)$ simplifies to $\beta_S (x_S)$, 
	indicating that  $\beta_S(\cdot)$ lacks dependence on the variable $x_{\overline S}$. As further clarified  in Appendix \ref{Appendix-CP}, this means that at time-scale q, the two competing activation processes largely unfold in parallel, with negligible interactions over physical time. 
	Instead,  for $q=p^{-1}$ or $q\gg p^{-1}$, $\beta_S$ depends on both  $x_R$ and $x_B$,
	indicating that $N_R$ and $N_B$  strongly interact on time-scales comparable to or asymptotically larger than
	$p^{-1}$.		
\end{Remark}
}

\subsection{Analysis of   $K^*$ and $A_S^*$} The following theorems build upon previous results by establishing both upper and lower bounds for the final number of active nodes (see \eqref{Kappa-Astar}).

\begin{Theorem}\label{prop:kminT}
$(i)$ It holds~\footnote{Of course $\lim\frac{{K^*}}{q}=\infty$  a.s.,  when 
	$\kappa_{\bm f}=\infty$. }
\begin{equation}\label{180323}
	\liminf \frac{{K^*}}{q}\ge  \kappa_{\bm f}, \quad\text{a.s.}
\end{equation}
\noindent $(ii)$ Provided that  $q=g$  and $\alpha_R>1$  or $g \ll q \ll p^{-1}$,  we have
\begin{equation} \label{180823-3}
	\liminf \frac{A^*_B}{q}\ge 	g_B(\kappa_{\bm g})+\alpha_B,\quad\text{a.s.}
\end{equation} 
where $g_B(\kappa_{\bm g})$ and $\kappa_{\bm g}$ are given in Proposition \ref{le:coupledphi}. 
\end{Theorem}

\begin{Theorem}\label{cor:BCadded2}
Let $S\in\{R,B\}$ be fixed.  If
$q=g$ and $\alpha_S<1$,  then
\[
\limsup \frac{A^*_S}{q}  \le  z_S +\alpha_S, \quad\text{a.s.}
\]
\end{Theorem}

\subsection{Analysis of the sequences $\{T_k\}_{k\in\mathbb N}$ and $\{T_k^S\}_{k\in\mathbb N}$ at time-scale $q$}

The next result describes the asymptotic behavior of $ T_{\lfloor\kappa q\rfloor}$ and  $T^S_{\lfloor \kappa_S q\rfloor}$,  
for appropriate constants $\kappa,\kappa_S>0$,  $S\in\{R,B\}$.   First, we define the scaling factor $\eta$ as follows: 
\begin{equation}\label{eta}
\eta := \left\{ \begin{array}{ll}
1  & \text{if }   q=g\\
\frac{n(qp)^r}{q} & \text{if }   g \ll q\ll p^{-1}\\
\frac{n}{q}  & \text{if either}\quad  q= p^{-1} \quad \text{or}\quad    q\gg  p^{-1}.
\end{array}\right.
\end{equation}
We then state the following theorem.
\begin{Theorem}\label{propVk}
$(i)$  For each $\kappa< \kappa_{\bm f}$,  we have
\begin{equation}\label{eqVk1}
	\eta	T_{\lfloor \kappa q \rfloor}\to \int_0^{\kappa} \frac{1}{\beta_R(\bm f(x))+\beta_B(\bm f(x)) } \dd  x, \quad \text{a.s.} 
\end{equation}\\
\noindent$(ii)$  Let $\kappa_S\in (0, \lim_{x\to \kappa_{\bm f}}f_S(x))$. Then
\begin{equation}\label{eqVk2}
\eta	T^S_{\lfloor \kappa_S q \rfloor}\to  \int_0^{f^{-1}_S(\kappa_S)} \frac{1}{\beta_R(\bm f(x))+ \beta_B(\bm f(x))} \dd x,   \quad  \text{a.s.}
\end{equation}
\end{Theorem}
Note that if $q\ll p^{-1}$,  then by \eqref{eq:CP} we have
\[
\int_0^{f_S^{-1}(\kappa_S)} \frac{1}{\beta_R(\bm f(x))+\beta_B(\bm f(x)) } \dd  x= \int_0^{\kappa_S} \frac{1}{\beta_S(y)} \dd  y. 
\]

\section{Proofs for the Time-Scale $q$ Analysis}\label{sec:V}
This section contains the proofs of Theorems \ref{prop:equadiff}, \ref{prop:kminT}, \ref{cor:BCadded2}, and \ref{propVk}, all of which build upon ancillary preliminary results. Here, we will only state these preliminary results,  deferring their (rather standard) proofs to Appendices \ref{sec:conc},  \ref{sec:prop5.3} and \ref{sec:lemmacompare}.  Finally,  we will demonstrate how Theorem \ref{thm:subcriticostoc} directly follows from Theorems \ref{prop:kminT} and \ref{cor:BCadded2}.  While the proofs of Theorems \ref{prop:kminT} and \ref{cor:BCadded2} are relatively simple,  those for Theorems \ref{prop:equadiff} and \ref{propVk} require more elaborate arguments.
Fig. \ref{log-dep1} summarizes the logical dependencies among findings.  
We suggest reading our proofs starting with the main results,  and then looking at the proofs of auxiliary results in the appendices.

\begin{figure}[t]
	\begin{tikzpicture}
	\tikzstyle{block} = [draw, fill=white, rectangle, very thick, rounded corners,	minimum height=3em, minimum width=6em]
	\node [block] at (0 ,0) (A) { \textsc{ Theorem \ref{prop:equadiff}} };
	
	\node [block] at (-2,-2) (B) { $\begin{array}{c}
			\textsc{ Proposition \ref{le:inequnifbis} }\\ 
		\end{array}$};
	\node [block] at (0 ,-3.5) (C) { $\begin{array}{c}
			\textsc{ Proposition \ref{cor:BC} }\\ 
\end{array}$};

	\node [block] at (3.5 ,0) (D) { \textsc{ Corollary \ref{cor:190823}} };
	\node [block] at (3.5 ,-1.5) (E) { \textsc{ Theorem  \ref{prop:kminT}}  };
	\node [block] at (3.5 ,-3.5) (F) { \textsc{ Theorem \ref{propVk}}  };
	\node [block] at (7 ,0) (G) { $\begin{array}{c}
			\textsc{ Proposition   \ref{le:coupledphi}}\\ 
\end{array}$ };
	\node [block] at (7 ,-3.5) (H) 
	{ $\begin{array}{c}
			\textsc{ Proposition \ref{indeptau} }\\ 
		\end{array}$ };
	\draw [->,very thick] (C) -- (A);
	\draw [->,very thick] (C) -- (B);
   \draw [->,very thick] (B) -- (A);	
	\draw [->,very thick] (A) -- (D);
	\draw [->,very thick] (C) -- (E);
		\draw [->,very thick] (A) -- (E);
		\draw [->,very thick] (A) -- (F);
	\draw [->,very thick] (D) -- (E);
	\draw [->,very thick] (C) -- (F);
	\draw [->,very thick] (D) -- (E);
	\draw [->,very thick] (E) -- (F);
	\draw [->,very thick] (G) -- (E);
	\draw [->,very thick] (H) -- (F);

	\node [block] at (10.5 ,-0.0) (L) {  $\begin{array}{c}
				\textsc{ Proposition \ref{lemma-compare} }\\ 
			\end{array}$  };
\node [block] at (10.5 ,-1.5) (M) { { \textsc{ Theorem \ref{cor:BCadded2} } } };  
\draw [->,very thick] (L) -- (M);

\node [block] at (10.5 ,-3) (N) { { \textsc{ Theorem \ref{thm:subcriticostoc} } } };  
\draw [->,very thick] (M) -- (N);
\draw [->,very thick] (E) -- (N);
\end{tikzpicture}		
\caption{Logical dependencies among scale $q$  results;  $A\to B$ means that $A$  is invoked in the  proof of  $B$. }\label{log-dep1}
\end{figure}

{
	\begin{Remark}
We emphasize that although Theorem~\ref{propVk} is not required for the derivation of Theorem~\ref{thm:subcriticostoc}, it plays a pivotal role in subsequent sections, particularly in the analysis of the system at scales larger than $q$.
	\end{Remark}	
}

\subsection{Further notation}

Letting $\bm k:=(k_R,k_B)\in  (\mathbb{N}\cup\{0\})^2$,  we define
\begin{equation}\label{170124}
	\mathbb{I}_k:=\{\bm k:\,\,k_R+k_B=k\},\quad k\in\mathbb N\cup\{0\}.
\end{equation}
In what follows, we consider $\kappa\in(0,\kappa_{\bm f})$, where $\kappa_{\bm f}$ is defined  in Definition \ref{def:CP} and computed
in Proposition \ref{le:coupledphi}.  We define the sets:
\[
\mathbb{T}(\kappa):= \left\{\begin{array}{ll}
\{\bm k:\,\,k_R+k_B\leq\kappa q\}=
\bigcup_{0\leq k\le \lfloor \kappa q\rfloor }\mathbb{I}_k  & \text{if } q \ll p^{-1} \text{ or }  q=p^{-1}\\
 \left\{ \bm k:\,\,k_R+k_B\le  \kappa q\quad \text{and} \quad \frac{k_R+\alpha_Rq}{k_B+\alpha_Bq}\ge  \frac{1}{2}+\frac{\alpha_R}{2\alpha_B}
 \right\}
   & \text{if } q\gg p^{-1},\\
\end{array} \right.\\
\]
and,  for $\bm{x}:=(x_R,x_B)\in [0,\infty)^2$,

\begin{align}
	\mathbb{T}'(\kappa):=\left\{ \begin{array}{ll}
		\{\bm x:\,\,x_R+x_B\leq\kappa \}  &  \text{if } q \ll p^{-1} \text{ or }  q=p^{-1}\\  
		\left\{\bm x:  x_R+x_B \le \kappa\quad\text{and}\quad
		\frac{x_R+\alpha_R}{x_B+\alpha_B} \ge \frac{1}{2} +\frac{\alpha_R}{2\alpha_B}\right\}  &  \text{if } q\gg p^{-1}. 
		\label{defKprimo}
	\end{array}\right.
\end{align}

Letting $z>0$ denote a constant such that $2z<\kappa$, for
$\bm{\ell}=(\ell_R,\ell_B)\in[0,\infty)^2$ such that
$\ell_R+\ell_B\leq\kappa-2z$, we define
\begin{equation}\label{170124-6}
	\mathbb{L}_{\bm{\ell}}(\kappa,z)	:=
	\left\{\bm x\in[0,\infty)^2:	x_R\geq \ell_R-z/2,\,
	x_B\geq \ell_B-z/2,\,	x_R+x_B\leq \ell_R+\ell_B+2z
	\right\}.
\end{equation}

\subsection{Auxiliary results} 

The proofs of Theorems \ref{prop:equadiff}, \ref{prop:kminT} and \ref{propVk} rely on Propositions \ref{cor:BC} and \ref{le:inequnifbis} below.  Their rather standard proofs are provided in Appendices \ref{sec:conc} and \ref{sec:prop5.3}, respectively.  The proof of Theorem \ref{cor:BCadded2} utilizes Proposition \ref{lemma-compare}, the proof of which can be found in Appendix \ref{sec:lemmacompare}.

{
\begin{Proposition}  \label{cor:BC}
	\label{prop:uniform-approximation}
	Let $\eta$ be defined as in~\eqref{eta} and let
	$\kappa\in(0,\kappa_{\bm f})$. For each $S\in\{R,B\}$, define
	\begin{equation}\label{defY}
	Y_S(\bm{k})
	:=
	 \ind_{\{\bm N[k]=\bm{k}\}}
	\left|
	U^S_{k+1}
	-
	\frac{|\beta_S(\bm{k}/q)|}
	{|\beta_R(\bm{k}/q)|+|\beta_B(\bm{k}/q)|}
	\right|,
	\end{equation}
	and
	\begin{equation}\label{defhatY}
	\widehat Y_S(\bm{k})
	:=
	 \ind_{\{\bm N[k]=\bm{k}\}}
	\left|
	Q^S_{k+1}
	-
	\eta q\,\beta_S(\bm{k}/q)
	\right|,
	\end{equation}
	where $\bm{k}=(k_R,k_B)$ and $k=k_R+k_B$.
	
	Moreover, let
	\begin{equation}\label{170124-3}
	\widehat N_S[j]
	:=
	N_S[j]-J_S[j],
	\qquad
	J_S[j]
	:=
	\sum_{h=1}^{\min\{j,n_W-1\}}U^S_h,
	\qquad
	\widehat N_S[0]:=0.
	\end{equation}
	
	Then
	\begin{equation}\label{eq:zero6nov}
	\Gamma_S(\kappa)
	:=
	\max\left\{
	\sup_{\bm{k}\in\mathbb T(\kappa)}Y_S(\bm{k}),
	\frac{1}{\eta q}
	\sup_{\bm{k}\in\mathbb T(\kappa)}
	\widehat Y_S(\bm{k}),
	\frac{1}{q}
	\sup_{j\le\kappa q}
	|\widehat N_S[j]|
	\right\}
	\longrightarrow0,
	\qquad\text{a.s.}
	\end{equation}
\end{Proposition}
}

Hereafter,  for 
$\kappa\in (0,\kappa_{\bm f})$,  we set
\begin{equation}\label{170124-4}
{\Omega_\kappa}:=\{\omega\in\Omega:  \max\left\{\Gamma_R(\kappa),\Gamma_B(\kappa)\right\} \to 0\}.
\end{equation}
Note that as an immediate consequence of  Proposition \ref{cor:BC}  it turns out that  $\mathbb{P}(\Omega_\kappa)=1$.

\begin{Proposition}\label{le:inequnifbis}
	{\it
For every $y,z>0$ such that $y+2z\le \kappa<\kappa_{\bm f}$,  $S \in \{R,B\}$ and $\omega \in \Omega_{\kappa}$,
we have:
\begin{align}
	&z\liminf\sum_{\bm k\in\mathbb{I}_{\lfloor yq\rfloor}} 
	\underline{\beta}_{S,\mathbb{L}_{{\bm k/q}}(\kappa,z)}
	\ind_{\{\bm N\big[\lfloor yq\rfloor\big]=\bm k\}}\le 
	\liminf\frac{\widetilde  N_S( yq+ zq)-\widetilde N_S( yq)}{ q} \nonumber \\
		\le & \limsup\frac{\widetilde N_S(yq+\ zq )-\widetilde N_S( yq )}{q}\leq z\limsup\sum_ {\bm k\in\mathbb{I}_{\lfloor yq\rfloor}}
	\overline{\beta}_{S,\mathbb{L}_{\bm k/q}(\kappa,z)}\ind_{\{\bm N\big[\lfloor yq\rfloor\big]=\bm k\}},\quad a.s.
	\label{eq:liminf}
\end{align}
Here,  for  $q\ll p^{-1}$ or $q=p^{-1}$:
\begin{align}\label{170124-7}
	\overline{\beta}_{S,\mathbb{L}_{\bm{\ell}}(\kappa,z)}:=
	\max_{{\bm x} \in\mathbb{L}_{\bm{\ell}} (\kappa,z)} \frac{|\beta_S({\bm x}) |}{|\beta_R({\bm x} )|+|\beta_{B}( 
		{\bm x} )|}, \qquad \underline{\beta}_{S,\mathbb{L}_{\bm{\ell}}(\kappa,z)}:=\min_{{\bm x}\in\mathbb{L}_{\bm{\ell}}(\kappa,z)}\frac{|\beta_S({\bm x})|}{|\beta_{R}({\bm x})|+|\beta_{B}({\bm x})|},
\end{align}
and,  for $q\gg p^{-1}$:
\begin{equation}  \label{170124-8}
	\overline{\beta}_{S,\mathbb{L}_{\bm{\ell}}(\kappa,z)}:=\max_{{\bm x}\in\mathbb{L}_{\bm{\ell}}(\kappa,z)}\frac{|\beta_S({\bm x})|}{|\beta_R({\bm x})|+|\beta_{B}({\bm x})|}\ind_{\{ \mathbb{L}_{\bm{\ell}}(\kappa,z)\subseteq \mathbb{T}'(\kappa)\}}+ \ind_{\{ \mathbb{L}_{\bm{\ell}}(\kappa,z)\not \subseteq \mathbb{T}'(\kappa)\}}
\end{equation}
and
\begin{equation}  \label{170124-9}
	\underline{\beta}_{S,\mathbb{L}_{\bm{\ell}}(\kappa,z)}:=\min_{{\bm x}\in\mathbb{L}_{\bm{\ell}}(\kappa,z)}\frac{|\beta_S({\bm x})|}{|\beta_{R}({\bm x})|+|\beta_{B}({\bm x})|}\ind_{\{ \mathbb{L}_{\bm{\ell}}(\kappa,z)\subseteq \mathbb{T}'(\kappa)\}}.
\end{equation}
}
\end{Proposition}

Exploiting standard coupling arguments,  one can compare the final number of $S$-active nodes,  $A^*_{S,h}$,  $h \in \{1,2\}$,  resulting from two activation processes 
with different numbers of $R$ and $B$ seeds. 
More precisely,  let $a_{S,h}$ denote the initial number of $S$-seeds for the $h$-th $S$-activation process. 
The following proposition holds.

\begin{Proposition}\label{lemma-compare}
	If $a_{R,1}\leq a_{R,2}$ and $a_{B,1}\geq a_{B,2}$, then 
	\[
	A_{R,1}^* \leq_{\mathrm{st}}A_{R,2}^*\qquad\text{and}\qquad  A_{B,2}^* \leq_{\mathrm{st}}A_{B,1}^*.
	\]
\end{Proposition}

\subsection{Proof of Theorem \ref{prop:equadiff}}

\subsubsection{Highlighting main conceptual steps}
To prove the uniform convergence of $\bm F(\cdot,\omega)$ to $\bm f(\cdot)$,  for almost all $\omega$,  we distinguish two cases: the case in which either $q  \ll p^{-1}$ or $q=p^{-1}$,  and the case in which $p^{-1}\ll q \ll n$.  In the first case the proof consists of four steps:
	\begin{enumerate}
\item[{\bf Step 1}.] We show that  functions $F_{S}(\cdot,\omega)$ are  a.s.  Lipschitz continuous and uniformly bounded  over compact domains.
\item[{\bf Step 2}.] By applying the Ascoli-Arzelà theorem, we prove that a subsequence of $\bm F(\cdot,\omega)$ converges pointwise to a limiting function,  for almost all $\omega$. 
\item[{\bf Step 3}.] We provide sufficiently tight upper and lower bounds for the incremental ratio of $\bm F(\cdot,\omega)$ near a fixed point,  for almost all $\omega$.  This allows us to show that the limiting trajectory is differentiable and that it is indeed the solution to the Cauchy problem in Definition \ref{def:CP}. 
\item[{\bf Step 4}.] The uniqueness of the solution of the Cauchy problem allows us to conclude that the whole sequence $\bm F(\cdot,\omega)$ converges pointwise to the limiting function,  almost surely.  Finally,  thanks to the regularity of both $\bm F(\cdot,\omega)$  and $\bm f(\cdot)$,  we lift the pointwise convergence to a uniform convergence over compacts.
\end{enumerate}

Apart from a few minor technical adjustments, the proof of the second
case proceeds as in the first one.

\subsubsection{Detailed proof}
We analyze separately the previously mentioned cases.
\\
\\
\noindent{\it Case $q  \ll p^{-1}$ or $q=p^{-1}$.}\\
\noindent{\bf Step 1.}
Since by Proposition \ref{cor:BC} we have $\mathbb{P}(\Omega_\kappa)=1$,  it suffices to prove \eqref{eq:fluidlimit} for all $\omega\in\Omega_\kappa$. 
{ By the definition of the linear interpolation, $F_{S,n}$ is
	piecewise affine. More precisely, on every interval
	$[k/q,(k+1)/q]$ its slope is
	\[
	N_S[k+1]-N_S[k]\in\{0,1\}.
	\]
	Consequently,
	\[
	|F_{S,n}(x_1)-F_{S,n}(x_2)|
	\le |x_1-x_2|,
	\qquad x_1,x_2\in[0,\kappa],
	\]
	and hence the family $\{F_{S,n}\}_n$ is uniformly
	$1$-Lipschitz. Moreover, $0\le F_{S,n}(x)\le x\le\kappa$.}

Thus, for any $\omega\in\Omega_\kappa$,  the functions $F_{S}(\cdot,\omega)$ are 1-Lipschitz  (i.e., Lipschitz   continuous with Lipschitz constant  equal to $1$) 
and uniformly bounded. 
From this point onward,  when it is necessary to avoid ambiguity,  we explicitly indicate  the dependence on $n$ of the various quantities.\\
\noindent{\bf Step 2.} Step 1 allows us to invoke the Ascoli-Arzel\'a theorem,  which guarantees the existence of a subsequence $\{F_{S,n'}(\cdot,\omega)\}_{n'}$ converging to some
function $f_S(\cdot,\omega)$,  uniformly on $[0,\kappa]$ ($f_S(\cdot,\omega)$ is Lipschitz continuous with Lipschitz constant equal to $1$ and is bounded above by $\kappa$).\\
\noindent{\bf Step 3.} 
For an arbitrarily fixed $x\in (0,\kappa)$ and  $z\in \left(x,\frac{\kappa+x}{2}\right)$ 
we have

\begin{align}
	f_S(z,\omega)-f_S(x,\omega)&=\lim_{n'\to\infty}[F_{S,n'}(z,\omega)-F_{S,n'}(x,\omega)]\nonumber\\
&=\lim_{n'\to\infty}q_{n'}^{-1}[\widetilde N_S( xq_{n'}+ (z-x)q_{n'} )(\omega)-\widetilde N_S( x q_{n'})(\omega)]\nonumber\\
&	\leq(z-x)\lim_{n'\to\infty}\sum_ {\bm k\in\mathbb{I}_{\lfloor xq_{n'}\rfloor}}\overline{\beta}_{S,\mathbb{L}_{\bm k/q_{n'}}(\kappa,z-x)}
	\ind_{\{\bm N\big[\lfloor x q_{n'}\rfloor\big](\omega)=\bm k\}},\label{added-EL1}	
\end{align}
where the inequality follows from Proposition \ref{le:inequnifbis}
(we refer the reader to \eqref{170124-6} for the definition of the set $\mathbb{L}_{\cdot}(\cdot,\cdot)$).
Let
$
x_{n'}:=\frac{\lfloor xq_{n'}\rfloor}{q_{n'}},
$
by construction we have
\begin{equation} \label{100825}
N_{S}\big[\lfloor x q_{n'}\rfloor\big] {(\omega)} =\widetilde N_S(\lfloor x q_{n'}\rfloor)(\omega)=\widetilde N_S(x_{n'}q_{n'})(\omega)=F_{S,n'}(x_{n'},\omega)q_{n'},
\end{equation}
and recalling the monotonicity and the Lipschitzianity of $F_{S,n}(\cdot,\omega)$ we obtain
\[
F_{S,n'}(x,\omega)-\frac{1}{q_{n'}}\leq F_{S,n'}(x,\omega)-(x-x_{n'})\leq F_{S,n'}(x_{n'},\omega)\leq F_{S,n'}(x,\omega).
\]
This implies
\begin{equation}\label{eq:lim7Nov}
\lim_{n'\to\infty}\bm F_{n'}(x_{n'},\omega)=\lim_{n'\to\infty}\bm F_{n'}(x,\omega)=\bm f(x,\omega),
\end{equation}
and therefore,  for any $\omega\in\Omega_\kappa$,  we have
\begin{align} \label{201023}
f_S(z,\omega)-f_S(x,\omega)&\leq
(z-x)\limsup_{n'\to\infty}\sum_ {\bm k\in\mathbb{I}_{\lfloor x q_{n'}\rfloor}}\overline{\beta}_{S,\mathbb{L}_{\bm k/q_{n'}}(\kappa,z-x)}
\ind_{\{\bm N [ x_{n'} q_{n'} ](\omega)=\bm k\}}\nonumber\\
&=(z-x)\limsup_{n'\to\infty}\overline{\beta}_{S,\mathbb{L}_{\bm F_{n'}(x_{n'},\omega)}(\kappa,z-x)}
=(z-x)\overline{\beta}_{S,\mathbb{L}_{\bm f(x,\omega)}(\kappa,z-x)},
\end{align}
where  the first equality follows from \eqref{100825}, and the identity, $\widetilde N_S(x_{n'}q_{n'})(\omega)= N_S[x_{n'}q_{n'}](\omega)$, while the second is a consequence of \eqref{eq:lim7Nov} and the continuity of the function $u\mapsto\overline{\beta}_{S,\mathbb{L}_{u}(\kappa,z-x)}$ (indeed the set ${\mathbb  L}_u$ varies continuously in the Hausdorff sense with $u$). 
Similarly,  for any $\omega\in\Omega_\kappa$, we have
\begin{equation*}
f_S(z,\omega)-f_S(x,\omega)\geq (z-x)\underline{\beta}_{S,\mathbb{L}_{\bm f{(x,\omega)}}\left(\kappa,z-x\right)},
\quad\forall\,z\in\left(x,\frac{\kappa+x}{2}\right).
\end{equation*}
Thus, for any $\omega\in\Omega_\kappa$,  any $x\in (0,\kappa)$ and any $z\in \left(x,\frac{\kappa+x}{2}\right)$,  we have
\begin{equation}\label{eq:UB}
\frac{f_S(z,\omega)-f_S(x,\omega)}{z-x}\leq\overline{\beta}_{S,\mathbb{L}_{\bm f{(x,\omega)}}(\kappa,z-x)},
\quad \frac{f_S(z,\omega)-f_S(x,\omega)}{z-x}\geq\underline{\beta}_{S,\mathbb{L}_{\bm f{(x,\omega)}}(\kappa,z-x)}.
\end{equation}
Since the set $\mathbb{L}_{\bm f{(x,\omega)}}(\kappa,z-x)$ is compact, it holds
\[
\overline{\beta}_{S,\mathbb{L}_{\bm f{(x,\omega)}}(\kappa,z-x)}=\frac{|\beta_S(\bm v)|}{|\beta_R(\bm v)|+
|\beta_{B}(\bm v)|}\quad\text{and}
\quad\underline{\beta}_{S,\mathbb{L}_{\bm f(x,\omega)}(\kappa,z-x)}=\frac{|\beta_S(\bm w)|}{|\beta_R(\bm w)|+|\beta_{B}(\bm w)|},
\]
for some 
\[
\bm v=(v_R,v_B),\bm w=(w_R,w_B) \in\mathbb{L}_{\bm f{(x,\omega)}(\kappa,z-x)}.
\]
By the definition of the set $\mathbb{L}_{\bm f(x,\omega)}(\kappa,z-x)$ it follows
\begin{equation}\label{eq:convmaxmin}
v_R,w_R\to f_R(x,\omega)\quad\text{and}\quad v_B,w_B\to f_{B}(x,\omega),\quad\text{as $z\downarrow x$.}
\end{equation}
{ Therefore, taking the $\limsup$ as $z\downarrow x$ in the first
	inequality in \eqref{eq:UB} and the $\liminf$ as $z\downarrow x$
	in the second one, and using \eqref{eq:convmaxmin} together with
	the continuity of $\beta_S$, we obtain
	\begin{equation}
		f_S^{'+}(x,\omega)
		=
		\Phi_S(\bm f(x,\omega)),
		\qquad x\in(0,\kappa),
		\label{eq:fS'+}
	\end{equation}
	where
	\[
	\Phi_S(\bm x)
	:=
	\frac{|\beta_S(\bm x)|}
	{|\beta_R(\bm x)|+|\beta_B(\bm x)|},
	\qquad S\in\{R,B\}.
	\]
	
	Observe that, since
$N_R[k]+N_B[k]=k,$
	the linear interpolation satisfies
$	F_R(x,\omega)+F_B(x,\omega)=x, $ 	and therefore every subsequential limit satisfies
$	f_R(x,\omega)+f_B(x,\omega)=x.$
	Hence, for $x\in[0,\kappa]$, the trajectory
	$\bm f(x,\omega)$ remains in $\mathbb T'(\kappa)$.
	
	Moreover, since $\kappa<\kappa_{\bm f}$,
	\[
	\inf_{\bm x\in\mathbb T'(\kappa)}
	\left(
	|\beta_R(\bm x)|+|\beta_B(\bm x)|
	\right)>0.
	\]
	Consequently, $\bm \Phi=(\Phi_R,\Phi_B)$ is locally
	Lipschitz on a neighborhood of $\mathbb T'(\kappa)$, and the
	Cauchy problem
	\begin{equation}
		\bm h'(x)
		=
		\bm\Phi(\bm h(x)),
		\qquad
		\bm  h(0)=(0,0),
		\label{eq:auxiliary-CP}
	\end{equation}
	admits a unique solution on $[0,\kappa]$.
	{ Moreover, since $\bm f(\cdot,\omega)$ is continuous and
		$\bm\Phi$ is locally Lipschitz, the mapping
		$
		x\mapsto
		\bm\Phi(\bm f(x,\omega))
		$
		is continuous. Since, by \eqref{eq:fS'+}, it coincides with the
		right-hand derivative of $\bm f(\cdot,\omega)$ on
		$(0,\kappa)$, it follows that $\bm f(\cdot,\omega)$ is
		differentiable on $(0,\kappa)$ and satisfies
		$
		\bm f'(x,\omega)
		=
		\bm\Phi(\bm f(x,\omega))
		$  (see e.g. Theorem A22 p. 541 of \cite{Bil}).
		Hence $\bm f(\cdot,\omega)$ is a solution of
		\eqref{eq:auxiliary-CP}.
		 }
	On the other hand, by Remark~\ref{re:29Ott} and
	Proposition~\ref{le:coupledphi}, along the deterministic
	trajectory $\bm f$ we have
	\[
	\beta_S(\bm f(x))\ge0,
	\qquad S\in\{R,B\},\quad x\in[0,\kappa].
	\]
	Therefore, along the deterministic trajectory,
	\[
	\Phi_S(\bm f(x))
	=
	\frac{\beta_S(\bm f(x))}
	{\beta_R(\bm f(x))+\beta_B(\bm f(x))},
	\]
	and thus $\bm  f$ is also a solution of
	\eqref{eq:auxiliary-CP}.
	
	Since the subsequential limit $\bm  f(\cdot,\omega)$
	satisfies \eqref{eq:auxiliary-CP} with the same initial
	condition, uniqueness implies
	\[
	\bm  f(x,\omega)=\bm  f(x),
	\qquad x\in[0,\kappa].
	\]
	Thus every convergent subsequence has the same deterministic
	limit $\bm  f$.}

\noindent{\bf Step 4.}{
	Suppose that $\bm F_n(\cdot,\omega)$ does not converge uniformly to
	$\bm f$ on $[0,\kappa]$. Then, for some $\varepsilon_0>0$, there
	exists a subsequence $\{n'\}$ such that
	\begin{equation}\label{hyp-step4-main}
	\sup_{x\in[0,\kappa]}
	\|\bm F_{n'}(x,\omega)-\bm f(x)\|
	\ge\varepsilon_0.
	\end{equation}
	By Step~1 and the Ascoli--Arzelà theorem, a further subsequence
	converges uniformly to some function $\bm h$. Repeating Step~3 shows
	that $\bm h$ solves \eqref{eq:auxiliary-CP}; hence, by uniqueness, we have necessarily
	$\bm h=\bm f$, in contradiction with \eqref{hyp-step4-main}. Therefore,
	$\bm F_n(\cdot,\omega)\to\bm f$ uniformly on $[0,\kappa]$.} 
\\
\\
{
\noindent{\it Case $p^{-1}\ll q\ll n$.}\\
Let $\mathring{\mathbb T}'(\kappa)$ denote the interior of
$\mathbb T'(\kappa)$ relative to $[0,\infty)^2$, where
$\mathbb T'(\kappa)$ is defined in \eqref{defKprimo}.
In this regime, $\bm{\beta}$ is discontinuous for $ \frac{x_R+\alpha_R}{  x_B+\alpha_{B}} =1$. We therefore argue locally inside
$\mathring{\mathbb T}'(\kappa)$.

Let $\widehat{\bm f}(\cdot)$ be any subsequential limit
obtained as in Steps 1--2. Since
\[
\widehat{\bm f}(0)=(0,0)
\]
and $\alpha_R>\alpha_B$, we have
$(0,0)\in\mathring{\mathbb T}'(\kappa)$.
Define
$
\tau_*:=
\inf\left\{
x\in[0,\kappa]:
\widehat{\bm f}(x)\notin
\mathring{\mathbb T}'(\kappa)
\right\},
$
with the convention $\inf\emptyset=\kappa$.
For every $x<\tau_*$, by continuity of
$\widehat{\bm f}$ we can choose $z>x$ sufficiently close
to $x$ so that
\[
\mathbb L_{\widehat{\bm f}(x)}(\kappa,z-x)
\subseteq
\mathbb T'(\kappa).
\]
Hence Proposition~\ref{le:inequnifbis} can be applied exactly as
in the previous case. Since, in the relative interior of
$\mathbb T'(\kappa)$, we have
$\beta_R(\bm x)=1$ and 
$ \beta_B(\bm x)=0, $
we obtain
\[
\widehat f_R'(x)=1,
\qquad
\widehat f_B'(x)=0,
\qquad x<\tau_*.
\]
Together with
$\widehat{\bm f}(0)=(0,0)$ this yields
\[
\widehat{\bm f}(x)=(x,0),
\qquad 0\le x<\tau_*.
\]	
	However, by Proposition~\ref{le:coupledphi}, $(x,0)$ belongs to the relative interior of 	$\mathbb T'(\kappa)$ for every $x\in[0,\kappa)$.
	If $\tau_*<\kappa$, continuity of
	$\widehat{\bm f}$ and the identity
	$\widehat{\bm f}(x)=(x,0)$ for $x<\tau_*$
	yield
	\[
	\widehat{\bm f}(\tau_*)=(\tau_*,0) 	\in\mathring{\mathbb T}'(\kappa),
	\]
	contradicting the definition of $\tau_*$. Hence
	$\tau_*=\kappa$. By continuity,
$	\widehat{\bm f}(\kappa)=(\kappa,0),$
	and consequently
	$
	\widehat{\bm f}(x)	=(x,0)=\bm f(x),	\text{ for } x\in[0,\kappa].
	$
The remainder of the proof proceeds as in the previous case.}

\subsection{Proof of Theorem \ref{prop:kminT} }
First we prove \eqref{180323} and then \eqref{180823-3}.

\subsubsection{ Proof of \eqref{180323}: Highlighting main conceptual steps}

The  proof of relation \eqref{180323}  is  divided  in two steps.
\begin{enumerate}
 \item[{\bf Step 1}.] Exploiting the properties of 
	$\bm f(x)$ (see Proposition \ref{le:coupledphi}) and the convergence results in Proposition \ref{cor:BC} and Theorem \ref{prop:equadiff},	
we show that,  for sufficiently large $n$,
	\[
	\min_{k\in [0,\kappa q]} \frac{{\max}\{Q_{k+1}^{R},Q_{k+1}^{B}\}}{\eta q}> 0,\qquad     \text{a.s.}
	\]
\item[{\bf Step 2}.]  To conclude the proof of relation \eqref{180323}, we observe that, since  
$Q_{K^*+1}^R=Q_{K^*+1}^{B}=0$, it necessarily follows that $K^*\ge  \lfloor\kappa q\rfloor$, 
$ \text{a.s.}$,
 for all sufficiently large $n$  and for any  $\kappa<\kappa_{\bm f}$. 

 \end{enumerate}
We emphasize that the uniform convergence of Theorem \ref{prop:equadiff} plays a key role in the proof of \eqref{180323}.
The proof of \eqref{180823-3} follows directly by \eqref{180323},  Corollary \ref{cor:190823} and Proposition \ref{le:coupledphi}.

\subsubsection{Detailed proof of \eqref{180323}}

We show the previously mentioned steps.
\\
\\
\noindent{\bf Step 1.}
Fix arbitrarily
\[
0<\kappa<\kappa_1<\kappa_{\bm f}.
\]
Let $\bm f$ be  as in \eqref{eq:CP} , we define the function
\[
b(\kappa):=\min_{x\in [0,\kappa]}{\max}\{\beta_R(\bm f(x)),\beta_{B}(\bm f(x))\} >0.
\]
The strict positivity  of $b$ follows  immediately from Remark \ref{re:29Ott} and Proposition \ref{le:coupledphi}.
For an arbitrarily fixed $\delta>0$,  we define the set 
\begin{equation}\label{eq:Bdelta}
	\mathbb{B}'_{{\bm f}}(\kappa,\delta):=\{\bm y\in [0, \infty)^2:
	\quad\min_{x\in[0, \kappa]}\|\bm y-\bm f(x)\|\leq\delta\}.
\end{equation}
For $q\ll p^{-1}$  or $q=p^{-1}$, $\bm{\beta}$ is uniformly continuous on
the compact set $\mathbb T'(\kappa_1)$. On the other hand, for
$p^{-1}\ll q\ll n$, we have
\[
\beta_R(\bm x)=1,\qquad \beta_B(\bm x)=0,
\qquad \bm x\in\mathbb T'(\kappa_1).
\]
Therefore, in all cases, $\delta_0>0$ can be chosen sufficiently
small so that
\[
\mathbb B'_{\bm f}(\kappa,\delta_0)
\subseteq\mathbb T'(\kappa_1)
\]
and
\[
\max_{\bm y\in\mathbb B'_{\bm f}(\kappa,\delta_0)}
\min_{x\in[0,\kappa]}
\|\bm\beta(\bm y)-\bm\beta(\bm f(x))\|
<\frac{b(\kappa)}4.
\]

This choice of $\delta_0$
leads to 
\begin{equation}\label{eq:11Novmat1}
\min_{\bm x\in\mathbb{B}'_{\bm f}(\kappa,\delta_0)}{\max}\{\beta_R(\bm{x}),\beta_{B}(\bm{x})\}\geq 3b(\kappa)/4.
\end{equation}
Based on Proposition \ref{cor:BC} and Theorem \ref{prop:equadiff}, 
we know that
\[
\sup_{\bm k\in\mathbb{T}(\kappa)}\max\{ \widehat{Y}_{R}(\bm k),\widehat{Y}_{B}(\bm k)\}/(\eta q)\to 0
\quad\text{and}\quad 
\sup_{x\in [0,\kappa]}\|\bm F(x)-\bm f(x)\|\to 0,\quad\text{a.s.}.
\]
This implies that for almost
every $\omega \in \Omega_\kappa$ there exists $n_0(\omega)$ such that for all $n>n_0(\omega)$:
\begin{equation}\label{eq:12Novsera1}
\text{$\bm F(x,\omega)\in\mathbb{B}'_{\bm f}(\kappa,{\delta_0})$\quad $\forall$ $x\in [0,\kappa]$}
\quad \text{and } \quad  	
  \sup_{\bm k\in\mathbb{T}(\kappa)}\max\{ \widehat{Y}_{R}(\bm k),\widehat{Y}_{B}(\bm k)\}/(\eta q)<b(\kappa)/4.
\end{equation}
{
Since $\bm N[k]/q=\bm F(k/q)$ and $k/q\in[0,\kappa]$ for
$0\le k\le\lfloor\kappa q\rfloor$, relations
\eqref{eq:11Novmat1}--\eqref{eq:12Novsera1} yield
\begin{equation}\label{eq:12Novpom1}
\min_{0\le k\le\lfloor\kappa q\rfloor}
\max\left\{
\beta_R\left(\frac{\bm N[k]}q\right),
\beta_B\left(\frac{\bm N[k]}q\right)
\right\}
\ge\frac{3b(\kappa)}4
\end{equation}
for almost every $\omega \in \Omega_\kappa$ and  $n>n_0(\omega)$.
}	
{
Moreover, for all sufficiently large $n$,
\[
\frac{\bm N[k]}{q}\in\mathbb T'(\kappa),
\qquad
0\leq k\leq\lfloor\kappa q\rfloor.
\]
Indeed, the constraint on the total number of activations follows
from
\[
N_R[k]+N_B[k]=k\leq\kappa q,
\]
whereas, for $p^{-1}\ll q\ll n$, the additional constraint defining
$\mathbb T'(\kappa)$ follows from
\[
\bm F(k/q)\in
\mathbb B'_{\bm f}(\kappa,\delta_0)
\subseteq\mathbb T'(\kappa_1).
\]
Hence, Proposition~\ref{cor:BC} yields, almost surely for all
sufficiently large $n$,
\[
\sup_{0\leq k\leq\lfloor\kappa q\rfloor}
\max_{S\in\{R,B\}}
\left|
\frac{Q^S_{k+1}}{\eta q}
-
\beta_S\left(\frac{\bm N[k]}{q}\right)
\right|
<
\frac{b(\kappa)}{4}.
\]
Combining this relation with \eqref{eq:12Novpom1}, we obtain
\begin{equation}\label{eq:minkappaq}
	\min_{0\leq k\leq\lfloor\kappa q\rfloor}
	\frac{\max\{Q^R_{k+1},Q^B_{k+1}\}}{\eta q}
	>
	\frac{b(\kappa)}{2}>0,
	\qquad\text{a.s.}
\end{equation}
}

\noindent{\bf Step 2.}
From the definition of $K^*$ and \eqref{eq:QSnew}, we have
$Q^R_{K^*+1}=Q^B_{K^*+1}=0.$ Hence, by \eqref{eq:minkappaq}, almost surely, for all sufficiently large $n$,
$
K^*>\lfloor\kappa q\rfloor.
$
Since $\kappa<\kappa_{\bm f}$ is arbitrary, \eqref{180323}
follows.
\\
\\
\noindent{\it Proof of  \eqref{180823-3}}\\
The proof of \eqref{180823-3} is  straightforward. We begin by observing  that:
\begin{align*} 
\liminf \frac{A^*_B}{q} 
= \liminf  \frac{N_B[K^*]}{q}+\alpha_B
\ge  \liminf  \frac{\widetilde N_B(\kappa q)}{q}+\alpha_B,   \quad  \forall \kappa>0 \quad\text{a.s.,}
\end{align*} 
where the inequality follows from \eqref{180323} and the monotonicity of $N_B(\cdot)$.  Therefore,  by Corollary \ref{cor:190823} and Proposition \ref{le:coupledphi},  we have
\[
\liminf  \frac{N_B[K^*]}{q}\ge   \lim_{\kappa\to \infty}  \liminf  \frac{\widetilde N_B(\kappa q)}{q}=\lim_{\kappa \to \infty} f_B(\kappa)= g_B(\kappa_{\bm g}),\quad a.s.
\]
and the proof is completed.

\subsection{ Proof of Theorem  \ref{cor:BCadded2} }
We will adopt  the notation of Proposition \ref{lemma-compare}.
The proof of  Theorem \ref{cor:BCadded2} relies on  comparing  the dynamics of two systems:
$(i)$ the original system (say system 1); 
$(ii)$ a companion system (say system 2)  where   $a_{\overline{S},2}=0$,   while  $a_{S,2}=a_{S,1}$.  As already noted in Remark \ref{remark-bootstrap},
the final size of $S$-active nodes in the companion system, say $A_{S,2}^*$,  equals the final size of active nodes
in a classical bootstrap percolation  process. 
 Using Proposition~\ref{lemma-compare}, after interchanging the roles
of $R$ and $B$ when $S=B$, and Theorem~3.2 in \cite{Bernoulli}, we
obtain
for any $\delta>0$ {there exist $c(\delta)>0$ and $n_\delta$ such that, for any $n\geq n_{\delta}$,}
\[
\mathbb{P}\left(\frac{A^*_{S,1}}{q}>z_S+\alpha_S+\delta\right)
\le \mathbb{P}\left(\frac{A^*_{S,2}}{q}>z_S+\alpha_S+\delta\right)=O(\exp(- c(\delta) q)).
\] 
Since $q\ge g$ and $g/\log n\to\infty$, the preceding bound is
summable; the claim follows from Borel--Cantelli.

\subsection{Proof of Theorem \ref{thm:subcriticostoc}}
The claim  is an immediate consequence of  Theorem  \ref{prop:kminT}$(i)$,  Theorem \ref{cor:BCadded2}  and \eqref{Kappa-Astar}. 
Indeed  recalling that $\kappa_{\bm f}=z_R+z_B$,  we have
\begin{align}
	z_S+\alpha_S\geq\limsup\frac{A_S^*}{q}\geq&\liminf\frac{A_S^*}{q}\geq\liminf\left(\frac{A^*}{q}-\frac{A_{\overline{S}}^*}{q}\right)\geq
	\liminf\frac{A^*}{q}+\liminf\left(-\frac{A_{\overline{S}}^*}{q}\right)\nonumber\\
	&\geq z_R+z_{B}+\alpha_R+\alpha_{B}-\limsup\frac{A_{\overline{S}}^*}{q}
	\geq z_S+\alpha_S,\quad a.s.\nonumber   
\end{align}

\subsection{Proof of Theorem \ref{propVk}}

Let $\bm Z$ be the Markov chain in Proposition
\ref{MarkovChain}. We recall that
\[
R(\bm z)=R^R(\bm z)+R^B(\bm z)\geq0,
\qquad \bm z\in\mathbb S,
\]
where
\begin{equation}\label{200524}
	R^S(\bm z)
	:=
	Q^S(\bm z)\ind_{\{\bm z^{(E)}=0\}}
	+
	(n_W-N(\bm z))U^S(\bm z)
	\ind_{\{\bm z^{(E)}=1\}}.
\end{equation}
Here $R^S(\bm z)$ represents the global rate at which the next
node to activate gets color $S$. For ease of notation, we set
\[
R_{k+1}:=R(\bm Z_k),
\qquad
R^S_{k+1}:=R^S(\bm Z_k).
\]

We first prove Part $(i)$.

\subsubsection{Highlighting main conceptual steps}

The proof of Theorem \ref{propVk}$(i)$ proceeds in six steps.
\begin{enumerate}
	\item[{\bf Step 1}.]
	We use Theorem \ref{prop:kminT}, Proposition \ref{cor:BC},
	and Theorem \ref{prop:equadiff} to establish deterministic
	upper and lower bounds for the total rate $R_{k+1}$, for large
	values of $n$.
	
	\item[{\bf Step 2}.]
	We note that, thanks to Proposition \ref{indeptau}, the
	sojourn times $\{W_k\}_{0\leq k<\lfloor\kappa q\rfloor}$ are
	conditionally independent given the corresponding sequence of
	rates $\{R_{k+1}\}_{0\leq k<\lfloor\kappa q\rfloor}$, and
	\[
	W_k\mid R_{k+1}
	\overset{\rm L}{=}
	{\rm Exp}(R_{k+1}).
	\]
	
	\item[{\bf Step 3}.]
	We prove that, for all $n$ sufficiently large and every
	sufficiently small $\varepsilon>0$, the random variables
	$\eta W_k$ can be upper and lower bounded by auxiliary random
	variables $\overline W_k^{(\varepsilon)}$ and
	$\underline W_k^{(\varepsilon)}$, respectively.
	
	\item[{\bf Step 4}.]
	As a consequence of Step 3, the sums of the auxiliary random
	variables provide upper and lower bounds for
	$\eta T_{\lfloor\kappa q\rfloor}$.
	
	\item[{\bf Step 5}.]
	We show that the sums of the auxiliary random variables are
	asymptotically concentrated around their expectations,
	denoted by $\overline\mu^{(\varepsilon)}(\kappa)$ and
	$\underline\mu^{(\varepsilon)}(\kappa)$.
	
	\item[{\bf Step 6}.]
	We conclude by showing that, letting first $n\to\infty$ and
	then $\varepsilon\downarrow0$, both
	$\overline\mu^{(\varepsilon)}(\kappa)$ and
	$\underline\mu^{(\varepsilon)}(\kappa)$ converge to
	\[
	\int_0^\kappa
	\frac{1}{
		\beta_R(\bm f(x))+\beta_B(\bm f(x))
	}\,\dd x.
	\]
\end{enumerate}

\subsubsection{Detailed proof}

\noindent{\bf Step 1.}

Fix
$0<\kappa<\kappa_1<\kappa_{\bm f},$
and define
\[
\beta_\Sigma(\bm x)
:=
\beta_R(\bm x)+\beta_B(\bm x),
\qquad
b(x):=\beta_\Sigma(\bm f(x)).
\]

By the properties of the solution of the Cauchy problem
\eqref{eq:CP}, $b(x)>0$ for every $x\in[0,\kappa_1]$.
Moreover, for $q\ll p^{-1}$  or $q=p^{-1}$, $b(\cdot)$ is continuous on
$[0,\kappa_1]$, whereas for $p^{-1}\ll q\ll n$,
Proposition~\ref{le:coupledphi} gives
\[
\bm f(x)=(x,0),
\qquad
b(x)=1.
\]
Hence
\begin{equation}\label{eq:propVk-bstar}
	b_*:=
	\inf_{x\in[0,\kappa_1]}b(x)>0.
\end{equation}

We now prove that
\begin{equation}\label{eq:propVk-total-rate}
	\sup_{0\leq k<\lfloor\kappa q\rfloor}
	\left|
	\frac{R_{k+1}}{\eta q}
	-
	b(k/q)
	\right|
	\longrightarrow0,
	\qquad\text{a.s.}
\end{equation}

Indeed, by Theorem~\ref{prop:kminT}, almost surely, for all
$n$ sufficiently large,
$
K^*>\lfloor\kappa_1q\rfloor.
$
Therefore,
\[
R_{k+1}
=
Q^R_{k+1}+Q^B_{k+1},
\qquad
0\leq k<\lfloor\kappa q\rfloor.
\]

Moreover, Theorem~\ref{prop:equadiff} yields
\begin{equation}\label{110223}
	\sup_{0\leq k\leq\lfloor\kappa q\rfloor}
	\left\|
	\frac{\bm N[k]}{q}
	-
	\bm f(k/q)
	\right\|
	\longrightarrow0,
	\qquad\text{a.s.}
\end{equation}
For $q\ll p^{-1}$  or $q=p^{-1}$, Proposition~\ref{cor:BC} gives, for
$S\in\{R,B\}$,
\begin{equation}\label{260826}
\sup_{0\leq k<\lfloor\kappa q\rfloor}
\left|
\frac{Q^S_{k+1}}{\eta q}
-
\beta_S\left(\frac{\bm N[k]}{q}\right)
\right|
\longrightarrow0,
\qquad\text{a.s.}
\end{equation}
Hence, by the uniform continuity of $\beta_\Sigma$ on a compact
neighborhood of $\{\bm f(x):0\leq x\leq\kappa\}$ and
\eqref{110223}, relation \eqref{eq:propVk-total-rate} follows.

For $p^{-1}\ll q\ll n$, the trajectory
	$\{\bm f(x):0\leq x\leq\kappa\}$ lies in the relative
	interior of $\mathbb T'(\kappa_1)$. Hence, by
	\eqref{110223}, for all sufficiently large $n$,
	\[
	\bm N[k]/q\in\mathbb T'(\kappa_1),
	\qquad
	0\leq k\leq\lfloor\kappa q\rfloor.
	\]
	Equivalently,
	$\bm N[k]\in\mathbb T(\kappa_1)$.
	We can therefore apply Proposition~\ref{cor:BC} with
	$\kappa_1$, obtaining that also in this case    \eqref{260826} holds.
	Since $\beta_R=1$ and $\beta_B=0$ on
	$\mathbb T'(\kappa_1)$, relation
	\eqref{eq:propVk-total-rate} follows.

Consequently, there exists a constant $c'>0$ such that, for
every sufficiently small $\varepsilon>0$, almost surely, for
all $n$ sufficiently large,
\begin{equation}\label{130825}
	\underline R_{k+1}(\varepsilon)
	\leq
	R_{k+1}
	\leq
	\overline R_{k+1}(\varepsilon),
	\qquad
	0\leq k<\lfloor\kappa q\rfloor,
\end{equation}
where
\[
\underline R_{k+1}(\varepsilon)
:=
\eta q
\left[
b(k/q)-c'\varepsilon
\right]
\quad \text{and}\quad 
\overline R_{k+1}(\varepsilon)
:=
\eta q
\left[
b(k/q)+c'\varepsilon
\right].
\]
We choose $\varepsilon>0$ sufficiently small so that
$
c'\varepsilon<\frac{b_*}{2}.
$
Thus
$
\underline R_{k+1}(\varepsilon)
\geq
\frac{b_*}{2}\eta q>0.
$

\noindent{\bf Step 2.}
Recall that
\[
W_k:=T_{k+1}-T_k.
\]
By Proposition~\ref{indeptau}, conditionally on the sequence
\[
\{R_{k+1}=r_{k+1}\}_{0\leq k<\lfloor\kappa q\rfloor},
\]
the random variables
$\{W_k\}_{0\leq k<\lfloor\kappa q\rfloor}$ are independent and
$
W_k
\overset{\rm L}{=}
{\rm Exp}(r_{k+1}).
$

\noindent{\bf Step 3.}
For
$
0\leq k<\lfloor\kappa q\rfloor,
$
define
\begin{equation}\label{eq:02072025uno}
	\underline W_k^{(\varepsilon)}
	:=
	\eta
	\frac{R_{k+1}}
	{\overline R_{k+1}(\varepsilon)}
	W_k,
	\qquad
	\overline W_k^{(\varepsilon)}
	:=
	\eta
	\frac{R_{k+1}}
	{\underline R_{k+1}(\varepsilon)}
	W_k.
\end{equation}
By \eqref{130825},  for all sufficiently large $n$
\begin{equation}\label{eq:10022023uno}
	\underline W_k^{(\varepsilon)}
	\leq
	\eta W_k
	\leq
	\overline W_k^{(\varepsilon)} \qquad \text{a.s.}
\end{equation}

Furthermore,
\begin{align}
	\overline W_k^{(\varepsilon)}
	\mid R_{k+1}
	&\overset{\rm L}{=}
	{\rm Exp}\left(
	\frac{\underline R_{k+1}(\varepsilon)}{\eta}
	\right),
	\label{eq:10022023due}
	\\
	\underline W_k^{(\varepsilon)}
	\mid R_{k+1}
	&\overset{\rm L}{=}
	{\rm Exp}\left(
	\frac{\overline R_{k+1}(\varepsilon)}{\eta}
	\right).
	\label{eq:10022023tre}
\end{align}

Conditionally on the whole sequence
$\{R_{k+1}\}_{0\leq k<\lfloor\kappa q\rfloor}$, the variables
$\{\underline W_k^{(\varepsilon)}\}$ are independent.
Moreover, their conditional joint distribution is
\[
\bigotimes_{0\leq k<\lfloor\kappa q\rfloor}
{\rm Exp}\left(
\frac{\overline R_{k+1}(\varepsilon)}{\eta}
\right),
\]
which does not depend on the realized values of
$\{R_{k+1}\}$. Therefore, after unconditioning,
$\{\underline W_k^{(\varepsilon)}\}$ are jointly independent.

The same argument shows that
$\{\overline W_k^{(\varepsilon)}\}$ are jointly independent.
In particular,
\[
\underline W_k^{(\varepsilon)}
\overset{\rm L}{=}
{\rm Exp}\left(
q[b(k/q)+c'\varepsilon]
\right)
\qquad \text{and}\qquad 
\overline W_k^{(\varepsilon)}
\overset{\rm L}{=}
{\rm Exp}\left(
q[b(k/q)-c'\varepsilon]
\right).
\]

\noindent{\bf Step 4.}

Since
\[
T_{\lfloor\kappa q\rfloor}
=
\sum_{k=0}^{\lfloor\kappa q\rfloor-1}W_k,
\]
relation \eqref{eq:10022023uno} implies that, almost surely,
for all sufficiently large $n$,
\[
\sum_{k=0}^{\lfloor\kappa q\rfloor-1}
\underline W_k^{(\varepsilon)}
\leq
\eta T_{\lfloor\kappa q\rfloor}
\leq
\sum_{k=0}^{\lfloor\kappa q\rfloor-1}
\overline W_k^{(\varepsilon)}.
\]
Therefore,
\[
	\liminf \sum_{k=0}^{\lfloor\kappa q\rfloor-1}
	\underline W_k^{(\varepsilon)}
	\leq
	\liminf	\eta T_{\lfloor\kappa q\rfloor}
	\leq
	\limsup	\eta T_{\lfloor\kappa q\rfloor}
	\leq
	\limsup \sum_{k=0}^{\lfloor\kappa q\rfloor-1}
	\overline W_k^{(\varepsilon)},
	\qquad\text{a.s.}
\]

\noindent{\bf Step 5.}

Define
\[
\mu_*(\kappa)
:=
\int_0^\kappa
\frac{1}{b(x)}
\,\dd x
=
\int_0^\kappa
\frac{1}{
	\beta_R(\bm f(x))+\beta_B(\bm f(x))
}
\,\dd x.
\]

The expectations of the two bounding sums are
\begin{equation}\label{eq:propVk:ave}
	\underline\mu^{(\varepsilon)}(\kappa)
	:=
	\sum_{k=0}^{\lfloor\kappa q\rfloor-1}
	\frac{\eta}
	{\overline R_{k+1}(\varepsilon)}
	=
	\frac{1}{q}
	\sum_{k=0}^{\lfloor\kappa q\rfloor-1}
	\frac{1}{b(k/q)+c'\varepsilon},
\end{equation}
and
\[
\overline\mu^{(\varepsilon)}(\kappa)
:=
\sum_{k=0}^{\lfloor\kappa q\rfloor-1}
\frac{\eta}
{\underline R_{k+1}(\varepsilon)}
=
\frac{1}{q}
\sum_{k=0}^{\lfloor\kappa q\rfloor-1}
\frac{1}{b(k/q)-c'\varepsilon}.
\]

Since
\[
b(k/q)-c'\varepsilon
\geq
\frac{b_*}{2},
\]
Thus, the rates of all the exponential random variables involved are at least   $b_*q/2$.

Hence, by the exponential tail bounds in
\cite{Janson-exp}, for every fixed sufficiently small
$\varepsilon>0$ and every $\delta>0$, there exists
$c_{\varepsilon,\delta}>0$ such that, for all sufficiently
large $n$,
\[
\mathbb P\left(
\left|
\sum_{k=0}^{\lfloor\kappa q\rfloor-1}
\underline W_k^{(\varepsilon)}
-
\underline\mu^{(\varepsilon)}(\kappa)
\right|
>
\delta
\right)
\leq
2{\rm e}^{-c_{\varepsilon,\delta}q},
\]
and an analogous bound holds for
$\sum_k\overline W_k^{(\varepsilon)}$.

Since $q/\log n\to\infty$, these probabilities are summable
with respect to $n$. 
{ Therefore, by the Borel--Cantelli lemma,
	\begin{equation}\label{eq:02072025tre}
		\sum_{k=0}^{\lfloor\kappa q\rfloor-1}
		\underline W_k^{(\varepsilon)}
		-
		\underline\mu^{(\varepsilon)}(\kappa)
		\longrightarrow0,
		\qquad
		\sum_{k=0}^{\lfloor\kappa q\rfloor-1}
		\overline W_k^{(\varepsilon)}
		-
		\overline\mu^{(\varepsilon)}(\kappa)
		\longrightarrow0,
		\quad\text{a.s.}
\end{equation}}


\noindent{\bf Step 6.}

For every fixed sufficiently small $\varepsilon>0$,
the quantities in \eqref{eq:propVk:ave} are Riemann sums.
Consequently,
\begin{equation}
	\underline\mu^{(\varepsilon)}(\kappa)
	\longrightarrow
	\int_0^\kappa
	\frac{1}{b(x)+c'\varepsilon}
	\,\dd x,
\quad \text{and}\quad
	\overline\mu^{(\varepsilon)}(\kappa)
	\longrightarrow
	\int_0^\kappa
	\frac{1}{b(x)-c'\varepsilon}
	\,\dd x,
	\qquad n\to\infty.
	\label{eq:02072025quattro}
\end{equation}

Combining Steps 4--6 gives
\[
\int_0^\kappa
\frac{1}{b(x)+c'\varepsilon}
\,\dd x
\leq
\liminf_{n\to\infty}
\eta T_{\lfloor\kappa q\rfloor}
\quad
\text{and}
\quad
\limsup_{n\to\infty}
\eta T_{\lfloor\kappa q\rfloor}
\leq
\int_0^\kappa
\frac{1}{b(x)-c'\varepsilon}
\,\dd x,
\qquad\text{a.s.}
\]

Finally, letting $\varepsilon\downarrow0$ and recalling that
$b(x)\geq b_*>0$, both bounding integrals converge to
$\mu_*(\kappa)$. Hence
\[
\eta T_{\lfloor\kappa q\rfloor}
\longrightarrow
\int_0^\kappa
\frac{1}{
	\beta_R(\bm f(x))+\beta_B(\bm f(x))
}
\,\dd x,
\qquad\text{a.s.},
\]
which proves \eqref{eqVk1}.

\medskip
\noindent{\it Proof of Part $(ii)$.}
Fix $S\in\{R,B\}$ and let
$
\kappa_S
\in
\left(
0,
\lim_{x\uparrow\kappa_{\bm f}}f_S(x)
\right).
$
Set
\[
x_S:=f_S^{-1}(\kappa_S).
\]
By the properties of $f_S$, $x_S\in(0,\kappa_{\bm f})$.
For an arbitrarily small $\delta>0$ such that
\[
0<x_S-\delta<x_S+\delta<\kappa_{\bm f},
\]
we have
\[
f_S(x_S-\delta)
<
\kappa_S
<
f_S(x_S+\delta).
\]

By Theorem~\ref{prop:equadiff}, almost surely, for all
sufficiently large $n$,
\[
N_S[\lfloor(x_S-\delta)q\rfloor]
<
\lfloor\kappa_Sq\rfloor
<
N_S[\lfloor(x_S+\delta)q\rfloor].
\]
Therefore,
\[
T_{\lfloor(x_S-\delta)q\rfloor}
\leq
T^S_{\lfloor\kappa_Sq\rfloor}
\leq
T_{\lfloor(x_S+\delta)q\rfloor}.
\]

Applying Part $(i)$ to the two terms on the right and on the
left gives
\[
\begin{split}
	\int_0^{x_S-\delta}
	\frac{1}{
		\beta_R(\bm f(x))+\beta_B(\bm f(x))
	}
	\,\dd x
	&\leq
	\liminf
	\eta T^S_{\lfloor\kappa_Sq\rfloor}
 \leq
	\limsup
	\eta T^S_{\lfloor\kappa_Sq\rfloor}
	\\
	&\leq
	\int_0^{x_S+\delta}
	\frac{1}{
		\beta_R(\bm f(x))+\beta_B(\bm f(x))
	}
	\,\dd x,
	\qquad\text{a.s.}
\end{split}
\]
Letting $\delta\downarrow0$ and recalling that
$x_S=f_S^{-1}(\kappa_S)$, we obtain
\[
\eta T^S_{\lfloor\kappa_Sq\rfloor}
\longrightarrow
\int_0^{f_S^{-1}(\kappa_S)}
\frac{1}{
	\beta_R(\bm f(x))+\beta_B(\bm f(x))
}
\,\dd x,
\qquad\text{a.s.},
\]
which proves \eqref{eqVk2}.

\section{Analysis at time-scales greater than $q$:  main results}\label{sub:scaling}

In this section, we analyze the joint dynamics of $\bm N[\cdot]$ and the pair $(|\mathcal{S}_R[\cdot]|, |\mathcal{S}_B[\cdot]|)$ over time scales that are  asymptotically larger than the number of seeds.
Recall that the function $g_B$ and the constant $\kappa_{\bm g}$ are defined in Proposition \ref{le:coupledphi}.
The following theorems hold (their proofs are provided in Section \ref{sect-theo->q}).

\begin{Theorem}\label{lemma-nuovo}
	If  either $(i)$ $q=g$  and $\alpha_R>1$ or $(ii)$  $g\ll q\ll p^{-1}$,  then 
	\begin{equation}\label{new-strong-timing-ineq-new}
		\text{$\forall$ $\varepsilon>0$, \quad $\mathbb{P}\left(\liminf\Big\{ \{ N_B\big[\lfloor f(n) p^{-1}\rfloor \big]  \le  \lfloor (g_B(\kappa_{\bm g})+\varepsilon)q\rfloor \}\cap\{
			{K^*}\ge 	\lfloor f(n) p^{-1}\rfloor   \}\Big\} \right)=1$,}
	\end{equation}
	where: under the assumption $(i)$, 
	$f$ is a (generic) function such that $f(n)\to \infty$ and $f(n)p^{-1}=o(n)$ and,  
	under the assumption $(ii)$,  $f(n):=c_0/(qp)^{r-1}\to \infty$,  for  a sufficiently small positive constant $c_0>0$.
\end{Theorem}
Informally, Theorem \ref{lemma-nuovo} states that for $q\ll p^{-1}$  the percolation process does not terminate before time-index
$\lfloor f(n)p^{-1}\rfloor$. Meanwhile, the number of $B$-activated nodes  remains $O_{\mathrm{a.s.}}(q)$.

\begin{Theorem}\label{thm:Lemma-new-2}
	Assume $q=g$ and  $\alpha_R>1$.
	Then
	\begin{equation}\label{new-third-part}
		\text{$\forall$ $\varepsilon>0$ and $c\in(0,1)$,  $\mathbb{P}\left(\liminf\Big\{\{ N_B[ K^*]\le \lfloor (g_B(\kappa_{\bm g})+\varepsilon)g \rfloor \} \cap  \{ K^*\ge \lfloor c n\rfloor \}\Big\} \right)=1.$}
	\end{equation}
\end{Theorem}

In the supercritical regime where $q=g$,  Theorem \ref{thm:Lemma-new-2} strengthens Theorem \ref{lemma-nuovo} by showing that the percolation process reaches time-index 
$\lfloor cn\rfloor$, before terminating. At the same time, the number of $B$-activated nodes  remains 
 $O_{\mathrm{a.s.}}(q)$.

\begin{Theorem}\label{thm:Lemma-new-3}
	Assume $g \ll q\ll n$.  
	Then
{ 
		\begin{equation}\label{eq:11082022-2-mod}
		\text{$\forall \, c\in (0,1)$,  $\mathbb{P}\left(
			\left\{\lim \frac{N_B\big[ K^* \big]}  {{ n } }=0\right\}\cap \liminf\{{K^*\ge  \lfloor cn\rfloor}\}   \right)=1.$}
	\end{equation}
}	
\end{Theorem}
Theorem \ref{thm:Lemma-new-3} applies  to the case $q\gg g$, demonstrating that the percolation process reaches time-index 
$\lfloor cn\rfloor$,  before terminating,  while the number of $B$-activated nodes stays within $o_{\mathrm{a.s.}}(n)$. 

{\begin{Remark}
		If $q\ll p^{-1}$,  our analysis is split into two stages. First, we examine the dynamics over time-scales $q'$ up to  an intermediate time scale denoted by $f(n)p^{-1}$ 
		(see Theorem \ref{lemma-nuovo}).  Subsequently,  we analyze the dynamics over time-scales greater than or equal to $f(n)p^{-1}$,  using Theorem \ref{thm:Lemma-new-2} when $q=g$,  and Theorem \ref{thm:Lemma-new-3} when  $g\ll q\ll p^{-1}$. For cases where  $q=p^{-1}$ or $ p^{-1}\ll q \ll n$,  we perform a direct analysis,  across all time-scales,  by applying Theorem \ref{thm:Lemma-new-3}. 
\end{Remark} }

\section{Proofs for the Analysis at Time-Scales Larger than $q$}\label{sect-theo->q}
This section contains the proofs of Theorems \ref{lemma-nuovo}, \ref{thm:Lemma-new-2},  \ref{thm:Lemma-new-3} and \ref{prop:supercritical}. 
The proofs of Theorems \ref{lemma-nuovo}, \ref{thm:Lemma-new-2} and \ref{thm:Lemma-new-3} rely on some ancillary results which are stated in  
Section \ref{preliminary-theo>q}.  Regarding the proof of Theorem \ref{prop:supercritical}: when $q=g$, it directly follows from Theorems \ref{thm:Lemma-new-2} and \ref{prop:kminT}$(ii)$. When $g\ll q$, it is an immediate consequence of Theorem \ref{thm:Lemma-new-3}.  Fig. \ref{log-dep2} summarizes the logical dependencies among our main findings.
\begin{figure}
	\begin{tikzpicture}
	\tikzstyle{block} = [draw, fill=white, rectangle, very thick, rounded corners,	minimum height=3em, minimum width=9em]
	\tikzstyle{block2} = [draw, blue, fill=white, rectangle, very thick, rounded corners,	minimum height=3em, minimum width=9em]
	\node[block2]   at (0,1.5)           {\bf Case $q=g$};
	\node [block] at (0 ,0.0) (F) { $\begin{array}{c}
				\textsc{ Theorems \ref{prop:kminT} }(i) \\ 
				\textsc{and \ref{propVk} }
			\end{array}$} ;          
	\node [block] at (3.5 ,0) (A) { \textsc{ Theorem \ref{lemma-nuovo}} $(i)$ };
	\node [block] at (7 ,0) (B) { \textsc{ Theorem \ref{thm:Lemma-new-2} } };
	\node [block] at (3.5 ,1.5) (U) { $\begin{array}{c}
			\textsc{ Lemma \ref{lemma:stop-CB}}\\ 
		\end{array}$};
	\node [block] at (7 ,1.5) (V) { \textsc{ Lemma \ref{thm:12072022}} };
	\node [block] at (3.5 ,-1.5) (W) { \textsc{ Lemma \ref{le:200523}} };
	\node [block] at (7, -1.5) (Z) { $\begin{array}{c}
			\textsc{ Lemma 	\ref{lemma-140223} }\\ 
		\end{array}$};
	
	\node [block] at (10.5 ,0) (Q) {
			\textsc{ Theorem \ref{prop:kminT}} $(ii)$};          	
\node [block] at (10.5 ,-1.5) (L) { \textsc{ Theorem  \ref{prop:supercritical} } };
	
	\draw [->,very thick] (A) -- (B);
	\draw [->,very thick] (F) -- (A);  
	\draw [->,very thick] (U) -- (A);  
	\draw [->,very thick] (V) -- (A);      
	\draw [->,very thick] (W) -- (A);        
	\draw [->,very thick] (Z) -- (B);       
	\draw [->,very thick] (B) -- (L);       
		\draw [->,very thick] (Q) -- (L);       	
	
	\node[block2]   at (0,-3.25)           {\bf Case $g\ll q \ll p^{-1}$};
	\node [block] at (0 ,-4.75) (F1) { \textsc{$\begin{array}{c}
				\textsc{ Theorems \ref{prop:kminT} }(i) \\ 
				\textsc{and \ref{propVk} }
			\end{array}$} };          
	\node [block] at (3.5 ,-4.75) (A1) { \textsc{ Theorem \ref{lemma-nuovo}} $(ii)$ };
	\node [block] at (7 ,-4.75) (B1) { \textsc{ Theorem \ref{thm:Lemma-new-3} } };
	
		\node [block] at (10.5 ,-4.75) (Q1) { \textsc{ Theorem \ref{prop:supercritical} } };
	\node [block] at (3.5 ,-3.25) (U1) { $\begin{array}{c}
			\textsc{ Lemma \ref{lemma:stop-CB} }\\ 
\end{array}$};
	\node [block] at (7 ,-3.25) (V1) { \textsc{ Lemma  \ref{thm:12072022}} };
	\node [block] at (3.5 ,-6.25) (W1) { \textsc{ Lemma \ref{le:200523}} };
	
	\draw [->,very thick] (A1) -- (B1);
	\draw [->,very thick] (F1) -- (A1);  
	\draw [->,very thick] (U1) -- (A1);  
	\draw [->,very thick] (V1) -- (A1);      
	\draw [->,very thick] (W1) -- (A1);        
	\draw [->,very thick] (B1) -- (Q1);

	\node[block2]   at (0,-8)           {\bf Cases $q= p^{-1}$  and  $p^{-1}\ll q\ll n$ };
	\node [block] at (1 ,-9.5) (F2) { \textsc{ $\begin{array}{c}
				\textsc{ Theorem \ref{prop:kminT} }(i) \\ 
				\textsc{and Corollary \ref{cor:190823} }
			\end{array}$ } };          
	\node [block] at (4.5 ,-9.5) (B2) { \textsc{ Theorem \ref{thm:Lemma-new-3} } };
	\node [block] at (4.5 ,-8) (W2) { \textsc{ Lemma \ref{le:200523}} };
		\node [block] at (8 ,-9.5) (Q2) { \textsc{ Theorem \ref{prop:supercritical} } };
	
	\draw [->,very thick] (F2) -- (B2);  
	\draw [->,very thick] (W2) -- (B2);  
	\draw [->,very thick] (B2) -- (Q2);	                            
\end{tikzpicture}
\caption{Logical dependencies among scale $q$  results;  $A\to B$ means that $A$  is invoked in the  proof of  $B$. }\label{log-dep2}		
\end{figure}

\subsection{Auxiliary results}  \label{preliminary-theo>q}  In this section, we present the auxiliary results that are invoked in the proofs of Theorems \ref{lemma-nuovo}, \ref{thm:Lemma-new-2}, and \ref{thm:Lemma-new-3}.

\subsubsection{A stochastic bound on \Ssupra{} nodes}

{
\begin{Lemma}\label{le:200523}
	\label{lemma:suprathreshold-bound}
	Fix $k\le n_W$ and two deterministic integers
	$\ell_R,\ell_B\ge0$. Define
	\begin{equation}
		\label{general-N-event}
		\mathcal N_{\ell_R,\ell_B}^{(k)}
		:=
		\left\{
		N_R[k]\ge\ell_R,\;
		N_B[k]\le\ell_B
		\right\}.
	\end{equation}
	Then, on the event
	$\mathcal N_{\ell_R,\ell_B}^{(k)}$, the following pathwise
	inclusions hold:
	\begin{equation}
		\label{pathwise-RB-bound}
		\mathcal S_R[k]
		\supseteq
		\widehat{\mathcal S}_R(\ell_R,\ell_B),
\qquad\text{and}\qquad
		\mathcal S_B[k]
		\subseteq
		\widehat{\mathcal S}_B(\ell_R,\ell_B).
	\end{equation}
	Consequently, for every $x>0$,
	\begin{equation}
		\label{prob-R-bound}
		\mathbb P\left(
		|\mathcal S_R[k]|\le x,\,
		\mathcal N_{\ell_R,\ell_B}^{(k)}
		\right)
		\le
		\mathbb P\left(
		\text{Bin}
		(n_W,\pi_R(\ell_R,\ell_B))
		\le x
		\right),
	\end{equation}
	and
	\begin{equation}
		\label{prob-B-bound}
		\mathbb P\left(
		|\mathcal S_B[k]|\ge x,\,
		\mathcal N_{\ell_R,\ell_B}^{(k)}
		\right)
		\le
		\mathbb P\left(
		\text{Bin}
		(n_W,\pi_B(\ell_R,\ell_B))
		\ge x
		\right).
	\end{equation}
\end{Lemma}
The detailed proof can be found in Appendix \ref{appendix-stop}.
Note that \eqref{prob-R-bound}  and  \eqref{prob-B-bound}  can be rewritten as:
\begin{equation}
	\label{prob-R-bound2}
	|\mathcal S_R[k]|\ind_{\mathcal N_{\ell_R,\ell_B}^{(k)}}+ n_W\left(1-  \ind_{\mathcal N_{\ell_R,\ell_B}^{(k)}}\right) 
	 \geq_{\mathrm{st}}
	\text{Bin}	(n_W,\pi_R(\ell_R,\ell_B))
\end{equation}

\begin{equation}
	\label{prob-B-bound2}
	|\mathcal S_B[k]|\ind_{\mathcal N_{\ell_R,\ell_B}^{(k)}}\,
	 \leq_{\mathrm{st}}
	\text{Bin}	(n_W,\pi_B(\ell_R,\ell_B)).
\end{equation}
}

{
\subsubsection{The stopped activation process}
We now introduce an auxiliary process,  hereafter called stopped process,  which is easier to analyze.  
In essence, the  stopped activation process  $N^{\text{stop}}=N^{\text{stop}}_R+N^{\text{stop}}_B$
proceeds as follows: up to a stopping time $Z_{\text{stop}}$ (either fixed or a point in the original process $N$)
$N^{\text{stop}}$ mirrors the original process $N^\circ$. If $Z_{\text{stop}}$  occurs before time $T_{K^*}$  the $R$-activation process halts at 
$Z_{\text{stop}}$  (no new $R$-active nodes).  Meanwhile, $B$-activation continues normally (i.e. following usual rules):  any jointly $W$ and \Bsupra{} node becomes $B$-active upon wake-up, until no jointly $W$ and  \Bsupra{}  nodes remain.
Formally,  on $\{t\le Z_{\text{stop}}\}$,  points in $N^{\text{stop}}_S$,  $S\in\{B,R\}$,  are obtained by thinning $\{(T'_k,V'_k)\}_{k\in\mathbb N}$,  retaining only those couples
$(T'_k,V'_k)$, $k\in\mathbb N$, for which,  at time $(T'_k)^{-}$,  the $W$ node $V'_k$
satisfies the \lq\lq threshold condition with respect to $S$\rq\rq.
On $\{t> Z_{\text{stop}}\}$,  $N_B^{\text{stop}}$ retains points satisfying the $B$-threshold condition, while 
$N_R^{\text{stop}}$ adds no new points, i.e., $N^{\text{stop}}_R(t)= N^\circ_R(\min\{t,  Z_{\text{stop}}\})$.
{
The stopped process is naturally defined up to its stopping index
$K^{*,\mathrm{stop}}$. For later use, we artificially prolong it
beyond its natural stopping time by allowing every remaining white
node to activate upon wake-up and assigning it color $B$. This
prolongation only serves to define discrete-time variables at indices
larger than $K^{*,\mathrm{stop}}$. The quantities
$K^{*,\mathrm{stop}}$, $A_S^{*,\mathrm{stop}}$, and
$T_k^{S,\mathrm{stop}}$ always refer to the natural, non-prolonged
stopped process, with $\inf\emptyset:=+\infty$ in the definition of
$T_k^{S,\mathrm{stop}}$.
}
 Variables associated with the stopped process will be denoted with a \lq\lq stop\rq\rq\, superscript or subscript to distinguish them from those of the original process. The properties defined in \eqref{eq:QSnew}, \eqref{180523-1}, and Lemma \ref{le:200523} all apply to this new stopped process.}
Additionally, a standard coupling argument, detailed in Appendix \ref{appendix-stop}, leads to the following lemma.
\begin{Lemma}\label{lemma:stop-CB}
	{
		\begin{equation} \label{stop-CB}
		A_B^{*,\mathrm{stop}}\geq A_B^*
		\qquad\text{and}\qquad
		T_k^{B,\mathrm{stop}}\leq T_k^{B,\circ},
		\qquad k\in\mathbb N\cup\{0\},
		\quad\text{a.s.}
		\end{equation}}
\end{Lemma}
{
The following lemma identifies the asymptotic black
activation times when the number of red activations
is kept bounded by a fixed multiple of $q$.

\begin{Lemma}\label{thm:12072022}
	Assume $q=g$ or $g\ll q\ll p^{-1}$.
	Let $Z_{\mathrm{stop}}$ be the stopping time used
	to define the stopped process, and suppose that,
	for some fixed $\kappa\geq0$,
	\[
	Z_{\mathrm{stop}}
	\leq T^{R,\circ}_{\lfloor\kappa q\rfloor},
	\qquad\text{a.s.}
	\]
	Fix $\kappa_B>0$ such that $\beta_B$ is strictly
	positive on $[0,\kappa_B]$. Then
	\begin{equation}\label{eqVk2stop}
		\eta T^{B,\mathrm{stop}}_{\lfloor\kappa_Bq\rfloor}
		\longrightarrow
		\int_0^{\kappa_B}\frac{1}{\beta_B(y)}\,\dd y,
		\qquad\text{a.s.}
	\end{equation}
	Equivalently, if $q=g$ and $\alpha_B\leq1$,
	we require $\kappa_B\in(0,z_B)$; in the other
	cases, $\kappa_B$ can be any positive constant.
\end{Lemma}
The proof follows the same lines as the proof
of Theorem \ref{propVk}, and therefore it is omitted.
}
%

\subsubsection{Asymptotic behavior of ordered non-negative random variables}

\begin{Lemma}	\label{lemma-140223}
	
Let $\{X_n\}_{n\geq 1}$  and $\{Y_n\}_{n\geq 1}$ be two sequences of non-negative random variables such that
	 $X_n  \leq_{\mathrm{st}}Y_n$ for any $n$.   If the random variables $\{Y_n\}_{n\geq 1}$ are independent and
	   $Y_n\to 0$ a.s. ,  then	$X_n \to  0$ a.s.

\end{Lemma}
The proof of this lemma is given in Appendix \ref{appendix-stop}.

	\subsection{Proof of Theorem \ref{lemma-nuovo} }

\subsubsection{Highlighting the main conceptual steps}

The proof of Theorem \ref{lemma-nuovo} can be divided in five steps.
{
	\begin{enumerate}
		\item[\bf{Step 1}.]		Our analysis at time-scale $q$ reveals that 
		$\mathbb{P}(\limsup\mathcal{A}_0^c)=0$,  where 
		{	\begin{equation} \label{030825-1}
			\mathcal{A}_0:=\left\{
			T^{B,\circ}_{{\primotheokappab}}>\tau_2, \;  T^\circ_{\lfloor {\kappa} q \rfloor }\le \tau_1,\; {K^*}> { \lfloor {\kappa}q \rfloor} \right\}. 
		\end{equation}  }
		Here ${\primotheokappab}:=\lfloor {(g_B(\kappa_{\bm g})}+\varepsilon)q \rfloor $  and  $0<\tau_1<\tau_2$ are suitable deterministic time levels.
		
		\item[\bf{Step 2}.] 
		{
		We define the sequence of random times  $Z_i$: 
\begin{equation}\label{030825-2}
	Z_{i+1}	:=	\min\left\{	T^\circ_{2N^\circ(Z_i)},	Z_i+\delta_i \right\},
	\quad
	Z_0:=T^\circ_{\lfloor\kappa q\rfloor} \quad  0\leq i<i_1
	:=
	\left\lceil
	\log_2
	\frac{\lfloor f(n)p^{-1}\rfloor}
	{\lfloor\kappa q\rfloor}
	\right\rceil .
\end{equation}

		}
		with  constants $\delta_i$  specified in \eqref{03022023-1} and satisfying $\sum_{i=0}^{i_1-1}\delta_i<\tau_2-\tau_1$. 
		
		\item[\bf{Step 3}.]  We prove  that 
{$\mathcal{A}_0 \subseteq \{N^\circ_B(Z_{i_1}) \le \primotheokappab \}$}. This inclusion implies  that the average number of \Rsupra{} nodes  stays high across any interval $[Z_i, Z_{i+1} )$,  $0\le i<i_1$. This in turn ensures a sufficiently high $R$-activation rate to guarantee {$ T^\circ_{2N^\circ(Z_i)}<Z_i+\delta_i$} a.s., while also preventing the percolation process from halting. We will  formalize this in the next two steps.
		\item[\bf{Step 4}.]  Defined  events
		$\mathcal{K}_{i}$ and  $\mathcal{Z}_i $  respectively  as in \eqref{030825-3} and \eqref{030825-4}, first we show that: 
	{	\begin{equation*}
			\mathcal{A}_0 \cap [\cap_{i=0}^{{i_1}-1}(\mathcal{K}_i\cap  \mathcal{Z}_i)]\subseteq  
			\{ N^\circ_B( \lfloor f(n)p^{-1} \rfloor)< {{\primotheokappab}}\}\cap\{  {K^*}\ge  \lfloor f(n)p^{-1} \rfloor \}.
		\end{equation*}
}		
		\item[\bf{Step 5}.]  Then we prove that $\mathbb{P}(	\liminf \mathcal{A}_0 \cap [\cap_{i=0}^{{i_1}-1}(\mathcal{K}_i\cap  \mathcal{Z}_i)] )=1$.
	\end{enumerate}
}	

\subsubsection{Detailed proof}

We prove Steps 1-5 previously described.
\\
\\
\noindent{\bf Step 1.}
By monotonicity in $\varepsilon$, it suffices to prove the claim for all sufficiently small  $\varepsilon>0$.
Hereafter, $\varepsilon>0$ is therefore fixed and sufficiently
small. 
In particular, when $q=g$ and $\alpha_B\le 1$, we choose
$\varepsilon$ such that
\[
g_B(\kappa_{\bm g})+\varepsilon<z_B.
\]
For $\alpha_B<1$, the inequality $g_B(\kappa_{\bm g})<z_B$ follows from
Proposition~\ref{le:coupledphi}.   The case $\alpha_B=1$ is discussed in Appendix \ref{proof-psi-tau}.

{
By Theorem~\ref{prop:kminT}$(i)$, and since
$\kappa_{\bm f}=+\infty$ in the regimes considered here, for every
fixed $\kappa>0$ we have, almost surely, eventually
$ K^*>\lfloor\kappa q\rfloor.$
Consequently, the original and extended processes coincide up to
the activation index $\lfloor\kappa q\rfloor$, and hence
\[
T^\circ_{\lfloor\kappa q\rfloor}
=
T_{\lfloor\kappa q\rfloor}
\qquad\text{eventually, a.s.}
\]
Therefore, Theorem~\ref{propVk} yields
\begin{equation}
	\eta T^\circ_{\lfloor\kappa q\rfloor}
	\longrightarrow
	\tau
	:=
	\int_0^\kappa
	\frac{1}{
		\displaystyle\sum_{S\in\{R,B\}}
		\beta_S(\bm f(y))
	}
	\,\dd y
	<\infty,
	\qquad\text{a.s.}
	\label{170424}
	\end{equation}}
For the application of Lemma~\ref{thm:12072022}, consider the
stopped activation process introduced above with
{
$
Z_{\mathrm{stop}}
:=
T^{R,\circ}_{\lfloor\kappa q\rfloor}.
$
Thus, the hypothesis
$Z_{\mathrm{stop}}\leq
T^{R,\circ}_{\lfloor\kappa q\rfloor}$  of Lemma~\ref{thm:12072022}
is satisfied with equality.
}
Recalling that
$\primotheokappab =\primotheokappab(\varepsilon) :=\lfloor(g_B(\kappa_{\bm g})+\varepsilon)q\rfloor, $
relation~\eqref{stop-CB} and Lemma~\ref{thm:12072022} yield,
almost surely,

\begin{equation}
	T^{B,\mathrm{stop}}_{\primotheokappab}
	\le
	T^{B,\circ}_{\primotheokappab},
	\qquad
	\eta T^{B,\mathrm{stop}}_{\primotheokappab}
	\longrightarrow
	\psi
	:=
	\int_0^{g_B(\kappa_{\bm g})+\varepsilon}
	\frac{1}{\beta_B(y)}\,\dd y.
	\label{170424-2}
	\end{equation}
As shown in Appendix~\ref{proof-psi-tau},
$ \psi>\tau. $
Therefore, recalling the definition of the event
$\mathcal A_0$ in \eqref{030825-1}, and setting
\[
\tau_1=\frac{\psi+2\tau}{3\eta},
\qquad
\tau_2=\frac{2\psi+\tau}{3\eta},
\]
relations \eqref{170424}, \eqref{170424-2}, and
Theorem~\ref{prop:kminT}$(i)$ imply
\begin{equation}
	\mathbb P(\limsup\mathcal A_0^c)=0.
	\label{090223}
\end{equation}

\noindent{\bf Step 2.}
{
Throughout Steps 2--5, the auxiliary quantities
$Z_i$, $K_i$, $\mathcal G^{(k)}$, $\mathcal K_i$,
$\mathcal Z_i$, and $\mathcal D_i$ are constructed from the
original physical process $N^\circ$. To avoid overburdening the
notation, we use unadorned symbols for these auxiliary quantities,
while retaining the superscript $\circ$ on the underlying physical
processes $N^\circ$ and $T^\circ$.}

Let $[Z_i,Z_{i+1})$ be the  intervals defined by  \eqref{030825-2}  with
	
	\begin{equation}\label{03022023-1}
		{\delta_i}:= \frac{2^{i+1}{\kappa}q}{\lambda_i}, \quad
		\lambda_i :=\left\{  \begin{array}{ll}
			\frac{\text{e}^{-1}}{2} n \frac{[(2^i\lfloor {\kappa}q \rfloor - \primotheokappab)p]^r}{r!} - 2^{i+1}\kappa q	 &\qquad\text{if } \quad  0\le i < {i_0}  \\
			c_1 n/3&\qquad \text{if } {i_0}  \le  i < i_1,
		\end{array}\right.
	\end{equation}	
	$ {i_0}:=
	\lfloor  \log_2 \frac{ \lfloor p^{-1}\rfloor }{\lfloor 2{\kappa}q\rfloor} \rfloor$   and $c_1$  is an appropriate strictly positive constant (specified in \eqref{2404023-2}).  Hereafter,  for the case $q=g$ we assume that 
	${\kappa}$ is chosen sufficiently large to guarantee $\lambda_i>0$ for any $i<i_0$.\footnote{For $q\gg g$, $\lambda_i>0$ is guaranteed for $n$ large enough since the second (negative) term is negligible with respect to the first (positive) term.}
	\\	
	\noindent{\bf Step 3.}
	{For sufficiently large $n$ and $\kappa$ the following holds: 
		\begin{equation}\label{Delta-piccolo}
			Z_{i_1}\le  Z_0+ \sum_{i<i_1}  {\delta_i}<  Z_0+\tau_2 -\tau_1 
		\end{equation}	 
		where the latter inequality can be easily verified by direct inspection.
		For $ 0\le i\le {i_1}$,  define  $K_i:=N^\circ(Z_i)$. By construction,
		we have $T^\circ_{K_i}\le Z_i<T^\circ_{K_i+1} $.
		From this relationship, together with  \eqref{Delta-piccolo} and the  monotonicity of the paths $N_B^\circ$, 
		we deduce
		\begin{align}
			&\mathcal{A}_0 \subseteq \{ N^\circ_B(T^\circ_{K_{i_1}}) \le \primotheokappab \}  \subseteq\{  N^\circ_B [k]\le {\primotheokappab} \;\;\forall k \in [K_0,\,  K_{i_1}]\}  \quad  \text{and}\label{A0theo61}\\
				&\mathcal{A}_0\cap \{ k\in [K_i, \,K_{i+1}]\}\subseteq \mathcal{G}^{(k)}:=\{ N^\circ_B [k]\le  {\primotheokappab} \} \cap \left\{
				K^*\geq k
				\right\}.\quad  \text{for any } 0\le i<i_1.  \label{070223-bis}
			\end{align} 	
Indeed, since,  the original process is absorbed at $K^*$, we have
$K_{i+1}=N^\circ(Z_{i+1})\le K^*$. Hence, every
$k\in[K_i,K_{i+1}]$ satisfies $k\le K^*$.
Note that 	by construction $K_{i+1}\le 2K_i$,   implies $K_i \le 2^i\lfloor {\kappa} q\rfloor$. \\
			\noindent{\bf Step 4.}
 \noindent{\bf Step 4.}
	{
		On $\mathcal G^{(k)}$, we have $k\le K^*$; hence, the original
		and extended processes coincide up to activation index $k$.
		Moreover,
		\[
		N_R^\circ[k]
		=
		k-N_B^\circ[k]
		\geq k-h_B.
		\]
		Therefore, for any
		$k\in[2^i\lfloor\kappa q\rfloor,
		2^{i+1}\lfloor\kappa q\rfloor]$,
		Lemma~\ref{le:200523} applies with
		$\ell_R=k-h_B$ and $\ell_B=h_B$.
}
Hence for every $x>0$,
			\begin{align}
				\mathbb P\left(	|\mathcal S_R[k]|\le x,\, {\mathcal G}^{(k)}	\right)	
				&\le \mathbb P\left(		\text{Bin}\left( n_W,	\pi_R(k-h_B,h_B)\right)\le x\right)	\nonumber\\
				&\le			\mathbb P\left(\text{Bin}	\left(n_W,	\pi_R		\left(	2^i\lfloor\kappa q\rfloor-h_B,h_B\right)\right)\le x\right),
				\label{030223-2bis}
			\end{align}

			Note that,  for any $i$ such that $2^i\lfloor {\kappa} q \rfloor< p^{-1}$,  we have
			\begin{align}
				\pi_R(2^i\lfloor {\kappa} q \rfloor -{\primotheokappab},{\primotheokappab})\ge & \mathbb{P}( \text{Bin}(2^i\lfloor {\kappa} q \rfloor-{\primotheokappab}+a_R,p) = r ) \mathbb{P}(\text{Bin}({\primotheokappab}+a_B,p)=0) \nonumber \\
				> &\frac{[(2^i {\kappa} q -{\primotheokappab})p]^r}{r!} \text{e} ^{-1 } (1 +o(1)).\label{2404023-1}
			\end{align}				
			Moreover, by the definition of $i_0$, for every $i\ge i_0$,
			$2^i\lfloor {\kappa} q \rfloor p$ is bounded away from zero, uniformly in $i$, while
			$\primotheokappab p=o(1)$. Hence, for a sufficiently small
			constant $c_1>0$ and all sufficiently large $n$,
			\begin{equation}
				\pi_R(2^i\lfloor {\kappa} q \rfloor-\primotheokappab,\primotheokappab)>c_1,
				\qquad i_0\le i<i_1.
				\label{2404023-2}
			\end{equation}
			Therefore,  defining
			\begin{align}	
				\mathcal{K}_{i}:=&
				\left\{ |\mathcal{S}_R[ h ] |  > \gamma_i\,\, 
				\forall  h \in[  K_i, K_{i+1}]  \right\},  \quad  \text{for any } 0\le i<  {i_1},  \label{030825-3}\\
				\text{ with }\qquad&
				\gamma_i :=\left\{  \begin{array}{lc}
					\frac{\text{e}^{-1}}{2} n \frac{[(2^i\lfloor {\kappa}q \rfloor - \primotheokappab)p]^r}{r!},  \ &  0\le i < {i_0} \\
					c_1n/2,   & {i_0}  \le  i< {i_1}, 
				\end{array}\right.   \label{03022023-2}
			\end{align}
			we obtain
			\begin{align}\label{200525-2}
				\mathcal{K}_{i}\cap \{  k\in [ K_i, K_{i+1}] \} \subseteq \mathcal{K}_i^{(k)} :=\{ | \mathcal{S}_R[ k ]|   > \gamma_i  \}  \quad \text{and}
			\end{align}
			\begin{align}
				\mathcal{K}_{i}^c=& \left\{ \exists k\in[ K_i,K_{i+1} ]:  |\mathcal{S}_R[ k ]|   \le  \gamma_i
				\right\}= \bigcup_k \left[\left(\mathcal{K}_i^{(k)} \right)^c\cap \{k \in [K_i,\,K_{i+1}]\}\right].\label{160523}
			\end{align}

		}
		Exploiting \eqref{180523-1} and \eqref{200525-2},  it can be immediately checked that,  for $\kappa$ and $n$ sufficiently large,  
		in both cases $q=g$ and $g\ll q\ll p^{-1}$,  we have
		\begin{align}
			Q^R_{k+1} \ind_{\left\{ 
				\mathcal{K}_{i}\cap \{  k\in [ K_i, K_{i+1}] \}   \right\} }&\ge\lambda_i \ind_{\left\{ 
				\mathcal{K}_{i}\cap \{  k\in [ K_i, K_{i+1}] \}   \right\}, }
			\qquad  0\le k<n_W,\quad 0\le i<i_1.
			\label{0906723-2}
		\end{align}
		with $\lambda_i$ as in \eqref{03022023-1}. 
		Define
		{
		\begin{equation}
			\mathcal Z_i	:=			
			\{Z_{i+1}=T^\circ_{2K_i}\},	\qquad  	\mathcal D_i
			:=
			\{K_0=\lfloor\kappa q\rfloor\}
			\cap
			\bigcap_{j<i}\mathcal Z_j.
			\label{030825-4}
		\end{equation}
		}
		 we immediately obtain
		\begin{equation}\label{040223}
			\mathcal{D}_i\subseteq \{K_i= 2^i\lfloor {\kappa}q \rfloor \}  \text{   and, in particular    }
			\mathcal{D}_{i_1}\	\subseteq \{ K_{{i_1}}= 2^{{i_1}} \lfloor \kappa q \rfloor \ge \lfloor f(n)p^{-1}\rfloor  \}.
		\end{equation} 
On
$\bigcap_{i<i_1}(\mathcal K_i\cap\mathcal Z_i),$
relation~\eqref{0906723-2} holds on each interval
$[K_i,K_{i+1}]$, while \eqref{040223} implies that these
intervals successively cover $[K_0,K_{i_1}]$. Recalling that
$Q^R_{K^*+1}=0$, we necessarily have
		\begin{equation}\label{190525}
			K^*\ind_{ \left\{ \mathcal{A}_0 \cap (\cap_{i< i_1} (\mathcal{K}_i\cap \mathcal{Z}_i))\right\}}\ge 
			K_{i_1} \ind_{ \left\{ \mathcal{A}_0 \cap (\cap_{i< i_1} (\mathcal{K}_i\cap \mathcal{Z}_i))\right\}}
			\ge \lfloor f(n)p^{-1}\rfloor 
			\ind_{ \left\{ \mathcal{A}_0 \cap (\cap_{i< i_1} (\mathcal{K}_i\cap \mathcal{Z}_i))\right\}},
		\end{equation}
		from which,  applying \eqref{A0theo61},  we obtain
		\begin{align}\label{200525}
			\mathcal{A}_0 \cap [\cap_{i=0}^{{i_1}-1}(\mathcal{K}_i\cap  \mathcal{Z}_i)]\subseteq  \mathcal{B}&:=
			\{ N^\circ_B[ \lfloor f(n)p^{-1} \rfloor]\le {{\primotheokappab}}\}\cap\{  {K^*}\ge  \lfloor f(n)p^{-1} \rfloor \}. \nonumber \\
			&=	\{ N_B[ \lfloor f(n)p^{-1} \rfloor]\le {{\primotheokappab}}\}\cap\{  {K^*}\ge  \lfloor f(n)p^{-1} \rfloor \}.
		\end{align}
		\noindent{\bf Step 5.}		
		Applying the Borel--Cantelli lemma,
		we can conclude that
		\[
		\mathbb{P}(\limsup  \mathcal B^c)=0  
		\]
		provided that
		\begin{equation}\label{eq:07072025quattro}
			\sum_n \mathbb{P}(\mathcal B^c)\le \sum_n  \mathbb{P}(\mathcal{A}_0^c\cup (  \cup_{i=0}^{i_1-1}( \mathcal{K}_i^c\cup \mathcal{Z}_i^c))) 
			<\infty.
		\end{equation}
		
		To check this relation,  note that  by the definition of $\mathcal{D}_i$
		in \eqref{030825-4}, we immediately have
		{  
			$\mathcal A_0\cap
			\bigcap_{j<i}(\mathcal K_j\cap\mathcal Z_j)
			\subseteq
			\mathcal D_i.$ 
			}
		and\footnote{We conventionally set $\cap_{j<0}\mathcal{Z}_j=\cap_{j<0}\mathcal{K}_j= \cap_{j< 0}( \mathcal{K}_j \cap \mathcal{Z}_j) = \Omega $.}
		\begin{align}
			\mathbb{P}(\mathcal{A}^c_0\cup(  \cup_{i=0}^{{i_1}-1} (\mathcal{K}^c_i\cup 
			\mathcal{Z}_i^c)))
			= & \mathbb{P} (\mathcal{A}^c_0)+ \mathbb{P} \left( 
			\mathcal{A}_0\cap \Big[\cup_{i=0}^{{i_1}-1}  
			\big(  \mathcal{K}^c_i \cap [\cap_{j<i}(\mathcal{K}_j\cap \mathcal{Z}_j )]\big) \Big] \right) \nonumber\\
			+&\mathbb{P} \left( \mathcal{A}_0 \cap \Big[\cup_{i=0}^{{i_1}-1} \big( (\mathcal{Z}^c_i \cap \mathcal{K}_i)\cap [\cap_{j<i}(\mathcal{K}_j\cap \mathcal{Z}_j )] \big) \Big] \right)\nonumber\\
			\le &   \mathbb{P} (\mathcal{A}^c_0)+ \mathbb{P} ( \mathcal{A}_0\cap [\cup_i  ( \mathcal{K}^c_i \cap \mathcal{D}_{i})] ) +
			\mathbb{P} ( \cup_i (\mathcal{Z}^c_i  \cap \mathcal{K}_i \cap \mathcal{D}_{i} ) ) \nonumber \\
			\le &  \mathbb{P} (\mathcal{A}^c_0)+ \sum_i  \mathbb{P} (\mathcal{A}_0\cap   \mathcal{K}^c_i \cap \mathcal{D}_i)
			+\sum_i \mathbb{P}(\mathcal{Z}^c_i \cap \mathcal{K}_i \cap \mathcal{D}_{i}).
			\label{120923}
		\end{align}
		The infinite sum $\sum_n   \mathbb{P}(\mathcal{A}_0^c)$ converges by the (second) Borel--Cantelli lemma and \eqref{090223}, since the events $\mathcal{A}^{(n)}_0$ are independent.
		Hence \eqref{eq:07072025quattro} immediately follows if we prove that
		\begin{equation}\label{eq:07072025due}
			\sum_n\sum_i  \mathbb{P} (\mathcal{A}_0\cap   \mathcal{K}^c_i \cap \mathcal{D}_i)<\infty
		\end{equation} 
		and 
		\begin{equation}\label{eq:07072025tre}
			\sum_n  \sum_i \mathbb{P}(\mathcal{Z}^c_i\cap \mathcal{K}_i \cap \mathcal{D}_{i})<\infty.
		\end{equation}
		\label{eq:07072025uno}
		By  \eqref{160523}, \eqref{040223} and the fact that $K_{i+1}\le 2K_i$,  we have
			\begin{align}
				\mathbb{P}( \mathcal{A}_0 \cap \mathcal{K}^c_i \ \cap \mathcal{D}_{i})=&
				\mathbb{P}(  	\cup_k [(\mathcal{K}_i^{(k)} )^c \cap \{k\in [K_i,K_{i+1}]\} ] \cap \mathcal{A}_0 \cap \mathcal{D}_{i})\nonumber\\		
				=&\mathbb{P}(  	\cup_{k=2^i\lfloor {\kappa} q \rfloor}^{2^{i+1}\lfloor {\kappa} q \rfloor} [(\mathcal{K}_i^{(k)} )^c  \cap \{k\in [K_i, K_{i+1}]\}  \cap \mathcal{A}_0 \cap \mathcal{D}_{i}])\nonumber\\		 
				\le& \sum_{k=2^i\lfloor {\kappa} q \rfloor}^{2^{i+1}\lfloor {\kappa} q \rfloor}\mathbb{P}(  (\mathcal{K}_i^{(k)} )^c \cap \mathcal{G}^{(k)}),
				\label{eq:07072025sei}
			\end{align}
			where the last inequality follows from \eqref{070223-bis}.
			By  \eqref{030223-2bis} and \eqref{200525-2} it follows
			\begin{align}
				\mathbb{P}\left( (\mathcal{K}_i^{(k)})^c\cap \mathcal{G}^{(k)} \right)
				&
				\le \mathbb{P}\left( \text{Bin}( n_W, \pi_R(  2^i\lfloor {\kappa} q \rfloor - {\primotheokappab},{\primotheokappab}))\le \gamma_i \right)
				\nonumber\\
				&\le 
				\left\{ \begin{array}{cc}  
					\text{e}^{ 
						-n \text{e} ^{-1}  \frac{(2^i{\kappa}qp)^r}{3r!}\FH\left(\frac{1}{2}\right)},  & 0\le i<  {i_0}  \\
					\text{e}^{ -
						\frac{c_1n}{2} \FH\left(\frac{1}{2}\right)},  & {i_0}\le  i< {i_1}
				\end{array} \right. 
				\label{eq:07072025cinque}
			\end{align}
			where the latter inequality follows from \eqref{2404023-1},  \eqref{2404023-2} and the concentration inequality \eqref{Penrose-coda-sotto}.  Relation 
			\eqref{eq:07072025due} immediately follows from  \eqref{eq:07072025sei} and \eqref{eq:07072025cinque}.  

{  As for \eqref{eq:07072025tre}, on
	\(\mathcal Z_i^c\cap\mathcal K_i\cap\mathcal D_i\), fewer than \(K_i\) activations occur during \([Z_i,Z_i+\delta_i]\) .
	 Hence,
	\[
	\mathcal Z_i^c\cap\mathcal K_i\cap\mathcal D_i
	\subseteq
	\left\{
	N^\circ(Z_i+\delta_i)-K_i<K_i
	\right\}.
	\]}			
			
Recalling that on $\mathcal D_i$, we have
$	K_i=2^i\lfloor\kappa q\rfloor$ and  $	Z_i=T^\circ_{K_i}$,	set
$	J_i:=N^\circ(Z_i+\delta_i)-K_i. $
	Then, on
	$\mathcal Z_i^c\cap\mathcal K_i\cap\mathcal D_i$,
	we have $0\leq J_i<K_i$ and
	\[
	K_{i+1}=K_i+J_i \quad \text{and}\quad 
	T^\circ_{K_i+J_i}
	\leq Z_i+\delta_i
	<
	T^\circ_{K_i+J_i+1}.
	\]
	Therefore
	\begin{equation}\label{250826-bis}
		\delta_i
		<
		T^\circ_{K_i+J_i+1}-T^\circ_{K_i}.
		\end{equation}
			We claim that, on
			$\mathcal Z_i^c\cap\mathcal K_i\cap\mathcal D_i$,
			\[
			K_i+J_i<K^*.
			\]
			Indeed, $K_i+J_i=K_{i+1}=N^\circ(Z_{i+1})\leq K^*$. If equality held, then
			\eqref{0906723-2}, applied at $k=K^*$, would give
			\[
			Q^R_{K^*+1}\geq\lambda_i>0,
			\]
			contradicting $Q^R_{K^*+1}=0$.
			Therefore
			\begin{equation}\label{250826-tris}
				\delta_i
				<
				T^\circ_{K_i+J_i+1}-T^\circ_{K_i}
				=  	T_{K_i+J_i+1}-T_{K_i}=
				\sum_{h=0}^{J_i}W_{K_i+h}.
			\end{equation}

			For $0\le h<2^i\lfloor\kappa q\rfloor$, define
			\[
			\widehat W_h^{(i)}
			:=
			\frac{
				R_{2^i\lfloor\kappa q\rfloor+h+1}
			}{\lambda_i}
			W_{2^i\lfloor\kappa q\rfloor+h}.
			\]
			By Proposition~\ref{indeptau}, conditionally on the corresponding
			sequence of rates, the random variables
			$\{\widehat W_h^{(i)}\}_{0\le h<2^i\lfloor\kappa q\rfloor}$
			are independent and exponentially distributed with parameter
			$\lambda_i$. Since their conditional joint distribution does not
			depend on the values of the conditioning rates, they are
			unconditionally i.i.d. with law $\operatorname{Exp}(\lambda_i)$.
			
Furthermore, on
$\mathcal Z_i^c\cap\mathcal K_i\cap\mathcal D_i$,
relation~\eqref{0906723-2}, together with
\eqref{201023-3}, \eqref{eq:QSnew}, and \eqref{200524},
implies, for all sufficiently large $n$, that
\[
R_{K_i+h+1}\ge\lambda_i,
\qquad 0\le h\le J_i.
\]		
Indeed, before $K^*$ we have
\[
R_{k+1}=Q^R_{k+1}+Q^B_{k+1}\ge Q^R_{k+1}\ge\lambda_i.
\]
Since, uniformly over the indices considered above,
$k=o(n)$, while $\lambda_i\le c n$ for some constant $c<1$,
we also have $n_W-k>\lambda_i$ for all sufficiently large $n$.
Therefore
\[
 W_{K_i+h}=W_{2^i\lfloor\kappa q\rfloor+h}\ \le\widehat W_h^{(i)},	\qquad 0\le h\le J_i.
\]
Combining this relation with \eqref{250826-tris}, and recalling that
$J_i<K_i=2^i\lfloor\kappa q\rfloor$ on $\mathcal D_i$, we obtain
the pathwise inclusion
\begin{equation}
	\mathcal Z_i^c\cap\mathcal K_i\cap\mathcal D_i
	\subseteq
	\left\{
	\sum_{h=0}^{2^i\lfloor\kappa q\rfloor-1}
	\widehat W_h^{(i)}
	>\delta_i
	\right\}.
	\label{050825-bis}
\end{equation}

				From which  we  obtain:
					\begin{align}
						\mathbb{P}\left(  \mathcal{Z}_i^c \cap  \mathcal{K}_i\cap \mathcal{D}_{i} \right) 
						& \le  \mathbb{P}\left( \sum_{h=0}^{2^i\lfloor {\kappa} q\rfloor -1} \widehat W^{(i)}_h >\delta_i \right) \nonumber\\
						&=  \mathbb{P}\left( \text{Po}( \lambda_i \delta_i)<2^i\lfloor {\kappa} q\rfloor  \right) 
						< \text{e}^{-\lambda_i\delta_i\FH(1/2)}, 
						\label{160523-2}
				\end{align}
				where the latter inequality follows from \eqref{Penrose-coda-sotto-Po}.  Using \eqref{160523-2} one can immediately verify
				\eqref{eq:07072025tre}.
				
\subsection{Proof of Theorem \ref{thm:Lemma-new-2}}

The proof of Theorem \ref{thm:Lemma-new-2}  is divided into two parts.
We first establish an upper bound on the number of $B$-activated nodes at the stopping time of the process.
				Specifically,  we will show that
				\begin{equation}\label{eq:07072025dieci}
					\mathbb{P}\left(\liminf\{ N_B[ {K^*}]\le \lfloor (g_B(\kappa_{\bm g})+\varepsilon)q \rfloor \} \right)=1.
				\end{equation}
				In the second part, we will prove that the total number of activated nodes at the stopping time is large. Specifically, we   demonstrate that
				\begin{equation}\label{eq:07072025undici}
					\mathbb{P}\left( \liminf \{{K^*}\ge \lfloor c n\rfloor \} \right)=1.
				\end{equation}
				
				\subsubsection{Highlighting the main conceptual steps in the proof of \eqref{eq:07072025dieci}}
					The core idea of this part of the proof is to analyze the simpler dynamics of a specially defined stopped process, where the stopping time is set to  $Z^{\text{stop}}:=T^\circ_{\lfloor f(n) p^{-1}\rfloor}$. We break down the proof into three main steps:
					\begin{enumerate}
						\item[\bf{Step 1}.]  We prove that,  given the event 
						$\mathcal{E}_0$   specified in \eqref{defboco6.2},
						the number of \Bsupra{} nodes (for the stopped process) in a right neighborhood of $Z^{\text{stop}}$ is
						asymptotically negligible ($o_{\mathrm{a.s.}}(g)$) in a right neighborhood of $Z^{\text{stop}}$. 
						\item[\bf{Step 2}.] From the result of Step 1,
						we deduce  that necessarily,  for any $\varepsilon>0$, $K^{*, \text{stop}}\le  \lfloor f(n) p^{-1}\rfloor
						+\lfloor \varepsilon g\rfloor$,  a.s.,  given that   $ \liminf {\mathcal{E}_0} = 1$. 
						\item[\bf{Step 3}.] 
						Finally, we conclude the proof by showing that previous properties of the stopped process immediately carry over to the original, unrestricted process  by leveraging \eqref{stop-CB}.  
					\end{enumerate}
					
					\subsubsection{Detailed proof of \eqref{eq:07072025dieci}}
We prove Steps 1–3 outlined above.
				\\
\\
				\noindent{\bf Step 1.}				
				Let $f(n)$ be as in  Theorem \ref{lemma-nuovo}$(i)$
				and  consider the stopped process with $Z_{\text{stop}}=T^\circ_{\lfloor f(n) p^{-1}\rfloor}$.
				Define the following quantities: ${{\secondotheokappazero}}:= \lfloor f(n) p^{-1}\rfloor$,  ${\secondotheokappauno}:= {\secondotheokappazero}+ \lfloor \varepsilon g \rfloor$, 
				with $\varepsilon>0$ arbitrarily fixed. Similarly, set   ${\secondotheokappazerob}:= 	\lfloor {(g_B(\kappa_{\bm g})}+\varepsilon)g \rfloor $,  $\secondotheokappaunob:= 	h^{(0)}_B +\lfloor \varepsilon g \rfloor $.
				Define the events
				\begin{equation}\label{defboco6.2}
				\mathcal{B}_0:=
				\{ N_B[{\secondotheokappazero}]\le  {\secondotheokappazerob}  \}=
				\{ T_{{\secondotheokappazero}}< 	T^B_{{\secondotheokappazerob}+1}\},\qquad 
				\ \mathcal{C}_0:=\{  {{K^*}} \ge  {{\secondotheokappazero}} \}  	\quad \text{and} \quad  {\mathcal E}_0={\mathcal B}_0\cap{\mathcal C}_0.
				\end{equation}
				From \eqref{180523-1},  on ${\mathcal C}_0$ we have
				$	Q_{{\secondotheokappazero}+1}^{B, \text{stop}}  = Q_{{\secondotheokappazero}+1}^{B} 
				\le  |\mathcal{S}_B[{\secondotheokappazero}]|.$
				By Lemma \ref{le:200523} 
			{
					\begin{align}
				&
				|\mathcal{S}^{\text{stop}}_B[{\secondotheokappazero}]| \ind _{\mathcal{E}_0} = 
				|\mathcal{S}_B[{\secondotheokappazero}]| \ind_{\mathcal{E}_0}
				 \leq_{\mathrm{st}}   \text{Bin}(n_W,\pi_B(  {\secondotheokappazero}- {\secondotheokappazerob} ,  {\secondotheokappazerob}  ))  \qquad \text{and} \nonumber\\
				&\qquad \quad  \text{Bin}(n_W,\pi_B( {\secondotheokappazero}- {\secondotheokappazerob},  {\secondotheokappazerob}  )
				)/g\to 0,\quad\text{a.s.}  \label{eq:23102023}
			\end{align}
			}
				Indeed, it is straightforward to check 
				$n_W\pi_B({\secondotheokappazero}- {\secondotheokappazerob},  {\secondotheokappazerob} )/g\to 0.$
				Applying the concentration inequality \eqref{Penrose-coda-sopra}, we obtain
				 $\mathbb{P}(\text{Bin}(n_W,\pi_B ( {\secondotheokappazero}-{\secondotheokappazerob}, {\secondotheokappazerob}) )>\varepsilon  g)< \exp(- \frac{\varepsilon g}2)$,  for $n$ sufficiently large and $\varepsilon>0$.  The claim \eqref{eq:23102023} follows by a standard application of  the Borel--Cantelli lemma.
				{ 
					On $\mathcal E_0$, only $B$-activations occur in
					the stopped process after index $h_0$. Hence,
					\[
					N_B^{\mathrm{stop}}[h_1]
					=
					N_B[h_0]+h_1-h_0
					\le h_B^{(0)}+h_1-h_0
					=
					h_B^{(1)}.
					\]
					Therefore, Lemma~\ref{le:200523} yields
					\[
					|\mathcal S_B^{\mathrm{stop}}[h_1]|
					\ind_{{\mathcal E}_0}
					\le_{\mathrm{st}}
					\operatorname{Bin}\!\left(
					n_W,\pi_B(h_1-h_B^{(1)},h_B^{(1)})
					\right).
					\]
					{ Repeating the argument leading to \eqref{eq:23102023}, with
						$(h_0,h_B^{(0)})$ replaced by $(h_1,h_B^{(1)})$, yields
						\[
						\frac{\operatorname{Bin}\!\left(
							n_W,\pi_B(h_1-h_B^{(1)},h_B^{(1)})
							\right)}{g}
						\longrightarrow 0,
						\qquad\text{a.s.}
						\]}
				
				Therefore,  by Lemma \ref{lemma-140223}, recalling that the above random variables,  for different $n$, 
				are independent,  we conclude
				\begin{equation}\label{140223-2}
					|\mathcal{S}^{\text{stop}}_B[{\secondotheokappazero}]|  \ind_{\mathcal{E}_0}=  o_{\mathrm{a.s.}}(g)   \quad\text{and}\quad 
					|\mathcal{S}^{\text{stop}}_B[{\secondotheokappauno} ]|  \ind_{\mathcal{E}_0}=   o_{\mathrm{a.s.}}(g).   
				\end{equation}
				\noindent{\bf Step 2.}
				We begin by observing that on $\mathcal{C}_0$
				\begin{equation}\label{100323}
					\mathcal{S}_B[{\secondotheokappazero}]= \mathcal{S}^{\text{stop}}_B[{\secondotheokappazero}] \subseteq   \mathcal{S}^{\text{stop}}_B[{\secondotheokappazero}+k]\subseteq
					\mathcal{S}^{\text{stop}}_B[{\secondotheokappauno}],  \qquad \forall k\le \lfloor \varepsilon g\rfloor.
				\end{equation}
				Indeed, in the stopped process, no node becomes $R$-active after
				 $Z^{\text{stop}}=T_{\secondotheokappazero}=T^{\text{stop}}_{\secondotheokappazero}$, 
			and therefore the number of \Bsupra{}  nodes is 
				monotonically increasing,  for all times after $T^{\text{stop}}_{{\secondotheokappazero}}$. 
				Moreover, we clearly have
				\begin{equation}\label{100323-2}
					\mathcal{V}_{B}^{\text{stop}} [{\secondotheokappazero}] = 	\mathcal{V}_{B}[{\secondotheokappazero}]  \quad\text{and}\quad 
					\mathcal{V}_{R}^{\text{stop}}[{\secondotheokappazero}+k] =		\mathcal{V}_{R}^{\text{stop}}[{\secondotheokappazero}] = 	\mathcal{V}_{R}[{\secondotheokappazero}],  
					\qquad  \forall k\le \lfloor \varepsilon g\rfloor.
				\end{equation}
			Finally, recall the following facts: $(i)$ up to time $T_{K^*}$,  only \Ssupra{} nodes become $S$-active;
				$(ii)$ a node can be  \Ssupra{} only if it has collected at least $r$ $S$-marks, i.e., 
				$\{v\in S_S(t)\} \subseteq \{ D_S^{(v)}(t)\ge r \}$; $(iii)$ for each node $v$
,  the number of $S$-marks collected, $D_S^{(v)}[k]$,  
				 is non-decreasing in $k$.  
				Then
				\begin{align}
					&  	\ind_{\mathcal{C}_0}|(\mathcal{V}_W\setminus\mathcal{S}^{\text{stop}}_B[{\secondotheokappauno} ])\cap
					\mathcal{V}_B[{\secondotheokappazero}]\cap
					\{ v: D^{v,\text{stop}}_B[{\secondotheokappauno}  ] \ge r \} 
					| \nonumber\\
					\stackrel{(a)}{=}& \ind_{\mathcal{C}_0}|(\mathcal{V}_W\setminus\mathcal{S}^{\text{stop}}_B[{\secondotheokappauno}])\cap    \mathcal{V}_B[{\secondotheokappazero}] |   	
					\stackrel{(b)}{\le}   \ind_{\mathcal{C}_0}|(\mathcal{V}_W\setminus\mathcal{S}_B[{\secondotheokappazero}] )\cap \mathcal{V}_B[{\secondotheokappazero}]|.
					\label{290823}
				\end{align}
				Here: the equality $(a)$ 
				follows because, conditional on $ \mathcal{C}_0$, by properties $(i)$–$(iii$) stated earlier, we have
				\[
				\mathcal{V}_B[{\secondotheokappazero}]\subseteq  \{ v: D^{v,\text{stop}}_B[{\secondotheokappazero}]\ge r  \}\subseteq \{ v: D^{v,\text{stop}}_B[{\secondotheokappauno} ]\ge r  \};
				\] 
				the inequality $(b)$ follows from \eqref{100323}.
				Therefore,  recalling that $N^{\text{stop}}_B[{\secondotheokappauno} ]=N_B[{\secondotheokappazero}]+\lfloor \varepsilon g\rfloor$,  and that,  conditional on $\mathcal{C}_0$,  we have 
				\[  
				\mathcal{V}^{\text{stop}}_{R}[{\secondotheokappazero}] \subseteq    \{ v:  D_R^{v, \text{stop}}[{\secondotheokappazero} ]\geq r   \}=
				\{ v:  D_R^{v, \text{stop}}[{\secondotheokappauno} ]\geq r   \},
				\]
				by \eqref{eq:QSnew},   \eqref{100323}, \eqref{100323-2}  and \eqref{290823} it follows
					\begin{align}
						Q_{{\secondotheokappauno} +1 }^ {B,\text{stop}}\ind_{\mathcal{E}_0}
						=&	\Biggl[ |\mathcal{S}_B^{ \text{stop}}[{\secondotheokappauno}]|-N^{\text{stop}}_B[{\secondotheokappauno}] 
						-  |\mathcal{S}^{\text{stop}}_B[{\secondotheokappauno}]\cap \mathcal{V}^{\text{stop}}_{R}[{\secondotheokappauno}]\cap\{ v:  D_R^{v, \text{stop}}[{\secondotheokappauno} ]\geq r   \}| \nonumber\\
						&\hspace{ 0.5 cm} + 
						|(\mathcal{V}_W\setminus\mathcal{S}^{\text{stop}}_B[{\secondotheokappauno}])\cap \mathcal{V}^{ \text{stop}}_B[{\secondotheokappauno}]\cap\{ v:  D_B^{v, \text{stop}}[{\secondotheokappauno}]\geq r   \}|  
						\Biggl] \ind_{\mathcal{E}_0} \nonumber\\
						\le & \ind_{\mathcal{E}_0} \Biggl[ |\mathcal{S}^{\text{stop}}_B[{\secondotheokappauno}  ]|-|\mathcal{S}_B[{\secondotheokappazero} ]|+
						|\mathcal{S}_B[{\secondotheokappazero} ]|
						-N_B[{\secondotheokappazero}]-\lfloor \varepsilon g\rfloor\nonumber\\
						&\hspace{ 0.5 cm} - |\mathcal{S}_B[{\secondotheokappazero}]\cap\mathcal{V}_R[{\secondotheokappazero}] |
						+|(\mathcal{V}_W\setminus\mathcal{S}^{\text{stop}}_B[{\secondotheokappazero}])\cap \mathcal{V}_B[{\secondotheokappazero}]
						|\nonumber\\
						& \hspace{0.5 cm} +  |(\mathcal{V}_W\setminus\mathcal{S}^{\text{stop}}_B[{\secondotheokappauno} ] )\cap
						(\mathcal{V}_B^{\text{stop}}[{\secondotheokappauno}]\setminus \mathcal{V}_B[{\secondotheokappazero}]) 
						\cap \{ v:  D^{v,\text{stop}}_B[{\secondotheokappauno} ]\ge r  \}|
						\Biggr]\nonumber\\
						\le &   \ind_{\mathcal{E}_0}\Biggl[ Q_{{\secondotheokappazero}+1}^{B, \text{stop}}+  |\mathcal{S}^{\text{stop}}_B[{\secondotheokappauno} ]|
						-\lfloor \varepsilon g\rfloor \nonumber\\ 
						& \hspace{ 0.5 cm} + |(\mathcal{V}_W\setminus\mathcal{S}^{\text{stop}}_B[{\secondotheokappauno}])\cap (\mathcal{V}_B^{\text{stop}}[{\secondotheokappauno}]\setminus \mathcal{V}_B[{\secondotheokappazero}]) 
						\cap \{ v:  D^{v,\text{stop}}_B[{\secondotheokappauno} ]\ge r  \}
						| \Biggr].\label{eq:09032023}
					\end{align}
				We note that the last addend in \eqref{eq:09032023} is bounded above by  $ \lfloor \varepsilon g\rfloor$,  since
				\[
				| (\mathcal{V}_B^{\text{stop}}[{\secondotheokappauno}]\setminus \mathcal{V}_B[{\secondotheokappazero}])  | =N^{\text{stop}}_B[{\secondotheokappauno}] -  N_B[{\secondotheokappazero}] = 
				\lfloor \varepsilon g\rfloor.
				\]
				Moreover this term is nonzero only on the event  $\{K^{*,\text{stop}}<  {\secondotheokappauno} \}$.
				(Indeed,  for any  $k$ such that ${\secondotheokappazero} < k \le {\secondotheokappauno}$, 
				on $\{K^{*,\text{stop}}\ge  {\secondotheokappauno}\}$,  we have: $ V^{\text{stop}}_{k}\in \mathcal{S}^{\text{stop}}_B[k]$  with 
				$ \mathcal{S}^{\text{stop}}_B[k]\subseteq  \mathcal{S}^{\text{stop}}_B[{\secondotheokappauno} ]$.
				In other words, 	on the event $\{K^{*,\text{stop}}\ge h_1\}$,
				\[
				\mathcal V_B^{\text{stop}}[h_1]	\setminus		\mathcal V_B^{\text{stop}}[h_0]
				\subseteq		\mathcal S_B^{\text{stop}}[h_1].
				\]
				
				Consequently,  we have
				\[
				Q_{{\secondotheokappauno}+1}^{B, \text{stop}} \ind_{\mathcal{E}_0}
				\le  \left[ Q_{{\secondotheokappazero}+1}^{B, \text{stop}}+  |\mathcal{S}^{\text{stop}}_B[{\secondotheokappauno} ]|
				-\lfloor \varepsilon g\rfloor \ind_{ \{K^{ *,\text{stop}}\ge  {\secondotheokappauno}\} } \right]\ind_{\mathcal{E}_0},  	
				\qquad\text{a.s.}.
				\]
				Combining 
				\eqref{180523-1}  and \eqref{140223-2}, and recalling that
				$ \mathcal{C}_0 \subseteq \{  Q_{{\secondotheokappazero}+1}^{B, \text{stop}} \ge 0\}$
				we  obtain
				\begin{equation}  \label{300123}
					\ind_{\{ K^{*,\text{stop}}>  {\secondotheokappauno} \}}\ind_{\mathcal{E}_0}  	 \le 
					\frac{-Q_{{\secondotheokappauno} +1}^{B, \text{stop}}+o_{\mathrm{a.s.}} (g ) }{ \lfloor \varepsilon g\rfloor }\ind_{\mathcal{E}_0},\quad a.s.
				\end{equation}
				Since
				$
				\{Q_{{\secondotheokappauno} +1}^{B, \text{stop}}<0\}\subseteq \{K^{*,\text{stop}}\le {\secondotheokappauno}\},
				$
				it follows 
				\begin{equation}  \label{300123-1}
					\ind_{\{ K^{*,\text{stop}}>  {\secondotheokappauno}   \} }\le \ind_{\{  Q_{{\secondotheokappauno} +1}^{B, \text{stop}}\ge 0 \}}.
				\end{equation}
				Multiplying both sides of \eqref{300123}      by $\ind_{\{  Q_{{\secondotheokappauno} +1}^{B, \text{stop}}\ge 0 \}}$ and applying    \eqref{300123-1} we obtain
				\[
				\ind_{\{ K^{*,\text{stop}}> {\secondotheokappauno} \}}\ind_{\mathcal{E}_0}  	\le 
				\frac{-Q_{{\secondotheokappauno}+1}^{B, \text{stop}}+o_{\mathrm{a.s.}} (g ) }{ \lfloor \varepsilon g\rfloor} 
				\ind_{\{  Q_{{\secondotheokappauno} +1}^{B, \text{stop}}\ge 0 \}}\ind_{\mathcal{E}_0},\quad a.s.
				\]
				Now observe that
				\[
				\limsup  \frac{-Q_{{\secondotheokappauno}+1}^{B, \text{stop}}+o_{\mathrm{a.s.}} (g ) }{ \lfloor \varepsilon g\rfloor} 
				\ind_{\{  Q_{{\secondotheokappauno} +1}^{B, \text{stop}}\ge 0 \}} \ind_{\mathcal{E}_0}\le \limsup \frac{o_{\mathrm{a.s.}} (g )}{g}\ind_{\mathcal{E}_0}= 0, \qquad\text{a.s.}
				\]
			  We deduce  that
				\[
				\ind_{\{ K^{*,\text{stop}}>  {\secondotheokappauno}  \} }\ind_{\mathcal{E}_0}   
				\to 0,\qquad\text{a.s.}
				\]
				This implies
				\[
				\ind_{\{ K^{*,\text{stop}}>  {\secondotheokappauno}   \} }\to 0,\qquad\text{a.s.}
				\]
				since by Theorem \ref{lemma-nuovo} we have
				$\lim\ind_{\mathcal{E}_0} =\liminf   \ind_{\mathcal{E}_0} = 1$  a.s.
				\\
				{\bf Step 3.} Relation \eqref{eq:07072025dieci}  follows from the inclusion
				\[
				\mathcal{E}_0 \cap  \{ K^{*,\text{stop}}\le   {\secondotheokappauno}   \} \subseteq \{A_B^{*,\mathrm{stop}}
				\le {\secondotheokappazerob}+ \lfloor \varepsilon g\rfloor
				+ a_B \},
				\] 
				together with  the arbitrariness of $\varepsilon$ and \eqref{stop-CB}.   

\subsubsection{Proof of \eqref{eq:07072025undici}}
The proof of  \eqref{eq:07072025undici} is rather simple.  
We begin by defining the  event  $\mathcal{A}_0 := \{ {K^*} \ge  {\secondotheokappazero}\} \cap \{ N_B[K^*]\le  {\secondotheokappazerob} \}.$
By  Theorem \ref{lemma-nuovo} and \eqref{eq:07072025dieci} we have 
$\mathbb{P}\left(\limsup \mathcal{A}_0^c\right)=0.$
Next, we  analyze the dynamics over intervals $[T_{k_i}, T_{k_{i+1}})$,  where
\[
k_i:={\min\{2^i {\secondotheokappazero}, \lfloor cn \rfloor\}},\qquad 0\le i\le i_1:=\Big\lceil \log_2 \frac{\lfloor c n \rfloor +1 }{\lfloor f(n) p^{-1}\rfloor} \Big\rceil .
\]
For any $0\le i<i_1$ we show that $ \mathbb{P}( {K^*}\in [ k_{i}, k_{i+1}), \; \mathcal{A}_0)\to 0$ sufficiently fast.
More specifically,   recalling that 
$Q^R_{K^*+1}=0$,
for any  $0\le i< i_1$,  we have
\[
\begin{split}
	\{ {K^*}\in [k_{i}, k_{i+1}) \}\cap  \mathcal{A}_0\subseteq \{ \exists   k \in [k_i,k_{i+1}) \text{ s.t. }  Q^R_{k+1}=0, N_B[k]\le  {\secondotheokappazerob}\} \quad \text{and}\\
	\mathbb{P}(\exists   k \in [k_i,k_{i+1}) \text{ s.t. }  Q^R_{k+1}=0, N_B[k]\le {\secondotheokappazerob} )\le \sum_{k=k_i}^{k_{i+1}-1} 
	\mathbb{P}(Q^R_{k+1}=0, N_B[k]\le  
	{\secondotheokappazerob}).
\end{split}
\]
On the event $\{N_B[k]\le h_B^{(0)}\}$, since
$N_R[k]+N_B[k]=k$, we have
\[
N_R[k]\ge k-h_B^{(0)}.
\]
{
		Moreover, uniformly over $0\le i<i_1$,
		\[
		(k_i-h_B^{(0)})p\ge f(n)-o(1)\longrightarrow\infty,
		\qquad
		h_B^{(0)}p=O(gp)=o(1),
		\]
		and hence
		\[
		\inf_{0\le i<i_1}
		\pi_R(k_i-h_B^{(0)},h_B^{(0)})\longrightarrow1.
		\]
		Since $n_W/n\to1$, it follows that, for all sufficiently large $n$,
		\[
		\inf_{0\le i<i_1}
		n_W\pi_R(k_i-h_B^{(0)},h_B^{(0)})
		\ge \frac{1+c}{2}\,n.
		\]
		}
So by \eqref{180523-1},  Lemma \ref{le:200523}  and the concentration inequality \eqref{Penrose-coda-sotto},
it follows 

\begin{align*}
	&\mathbb{P}( {K^*}\in [ k_{i}, k_{i+1}), \; \mathcal{A}_0 ) \le 
	\hspace{-2mm} \sum_{k=k_i}^{k_{i+1}-1} \hspace{-2mm} 
	\mathbb{P}\left(  | \mathcal S_R[k]|\le k, N_B[k]\le   {\secondotheokappazerob} \right) 
	 \\		
	& \qquad \qquad \le 2^i {\secondotheokappazero}\mathbb{P}\left( \text{Bin} (n_W, \pi_R(k_i-{\secondotheokappazerob} , {\secondotheokappazerob} ) )< k_{i+1} \right)
	< \exp\left(-cn \FH\left(\frac{c}{\frac{1}{2}+ \frac{c}{2} } \right) \right),
\end{align*}
for any $0\le i< i_1$ and any $n$ large enough.  
{
	Since $i_1=O(\log n)$, the preceding exponential bound implies
	\[
	\sum_{n=1}^{\infty}
	\sum_{i=0}^{i_1-1}
	\mathbb P\left(
	K^*\in[k_i,k_{i+1}),\,\mathcal A_0
	\right)
	<\infty.
	\]
	Hence, by the first Borel--Cantelli lemma,
	$
	\mathbb P\left(
	\limsup
	\bigl(
	\{K^*<\lfloor cn\rfloor\}\cap\mathcal A_0
	\bigr)
	\right)=0.
	$
	Together with
	$
	\mathbb P(\limsup\mathcal A_0^c)=0,
	$
	this yields
	$
	\mathbb P\left(
	\liminf\{K^*\ge\lfloor cn\rfloor\}
	\right)=1,
	$
	which proves \eqref{eq:07072025undici}.
}
\subsection{Proof of Theorem  \ref{thm:Lemma-new-3}}
\subsubsection{Highlighting the main conceptual steps}
{The analysis is conducted recursively over the sequence of intervals $[Z_i,Z_{i+1})$,  where 
	\begin{align} 
		Z_0 :=\min\{T_{ \terzotheokappazero   }, T^B_{{\terzotheokappazerob }}\}\quad\text{and}\quad 
		Z_{i+1}:=\min\{T_{4^{i+1}\terzotheokappazero   }, T^B_{2^{i+1} {\terzotheokappazerob } }, T_{\lfloor cn \rfloor}\},
		\quad\text{$i\ge 0, $}\label{eq:TsuCo}
	\end{align}
	where the constants $\terzotheokappazero  $ and $\terzotheokappazerob $ are specified in  \eqref{100124}.  
	Informally, our arguments show that the $R$-activation process largely outpaces the $B$-activation process within each interval 
	$[Z_i,Z_{i+1})$.  This ensures that the events
	\begin{equation}\label{eq:16122022primo}
		\mathcal{A}_i:=
		\left\{T^B_{2^{i+1} {\terzotheokappazerob } }\ge  \min\left\{T_{4^{i+1}\terzotheokappazero }, T_{\lfloor cn\rfloor}\right\}\right\}\subseteq
		\left\{N_B[\min( 4^{i+1}\terzotheokappazero ,  \lfloor cn\rfloor)]\le  2^{i+1} {\terzotheokappazerob }\right\}
	\end{equation}
	occur with a probability that rapidly approaches $1$ for every meaningful $i$.
	   Furthermore, the number of \Ssupra{} nodes remains large enough to guarantee that the activation process never stops in the intervals defined above. More technically,  setting
	\begin{equation}\label{eq:20072022}
		I:=
		\min\{i:\,\, Z_{i+1}= T_{\lfloor c n\rfloor}\}=	\min\{i:\,\, T_{4^{i+1}\terzotheokappazero   }\ge  T_{\lfloor cn \rfloor},\,\,  T^B_{2^{i+1} {\terzotheokappazerob } }\ge T_{\lfloor cn \rfloor}\},
	\end{equation}
	Define
	\begin{equation}\label{eq:16122022secondo}
	\mathcal{B}_i:= \{Q^R_{h+1}>\lambda_i \text{ and } Q^B_{h+1}\le \phi_i
	\; \forall h \in [K_i,K_{i+1}), I\ge i\}\cup\{I<i\},
\end{equation}
where $\lambda_i$, $\phi_i$,  are suitable positive quantities and $K_i:=N(Z_i)$,
	we show that the probability that both $\mathcal A_i$ and $\mathcal B_i$
	occur tends rapidly to $1$ for every $i\geq0$.

	The claim then follows,  as in previous theorems,  by applying  the Borel--Cantelli lemma.
	The detailed proof is organized in three parts.  In the first part,  we prove the claim assuming that certain technical inequalities
	(i.e.,  \eqref{eq-new-b-q>>g} and \eqref{eq-new-a-q>>g}) are verified; in the second part we prove the first technical inequality (namely,  \eqref{eq-new-b-q>>g}); in the third part we prove the second technical inequality  (namely,  \eqref{eq-new-a-q>>g}).
	
	Before going through the details of the proof,  we introduce some notation.
	Let $f(n)$
	be the function considered in the statement of Theorem
	\ref{lemma-nuovo}$(ii)$ (i.e.,  for the case $g\ll q \ll p^{-1}$).
	Set
	\begin{equation}\label{100124}
		{\terzotheokappazero }:= \left\{ \begin{array}{ll}  \lfloor f(n)p^{-1} \rfloor  &   {\rm if}\quad g \ll q\ll p^{-1}\\
			\lfloor \kappa p^{-1} \rfloor  &  {\rm if}\quad q=  p^{-1}\\
			\lfloor \kappa q  \rfloor  &  {\rm if}\quad   p^{-1}\ll q \ll n 
		\end{array} \right.    \qquad {\terzotheokappazerob }:= \left\{ \begin{array}{ll}  \lfloor p^{-1} \rfloor  & {\rm if}\quad  g \ll q\ll p^{-1}\\
			\lfloor \overline{f}_B p^{-1} \rfloor  &  {\rm if}\quad q=  p^{-1}\\
			\lfloor q  \rfloor  &  {\rm if}\quad   p^{-1}\ll q \ll n  
		\end{array} \right.
	\end{equation}
	where { $\kappa$ } is a sufficiently large positive constant and $\overline{f}_B$ is defined in  Proposition \ref{le:coupledphi}. 
}
Due to the arbitrariness of $\kappa$ note that the ratio $\terzotheokappazero /\terzotheokappazerob$ can be assumed arbitrarily large for $n$ large enough.

\subsubsection{Part 1}

The proof unfolds over five steps.
\begin{enumerate}
	\item[{\bf Step 1}.] Preliminary relations are introduced, followed by the full definition of $\mathcal{B}_i$.
	\item[{\bf Step 2}.] We prove  that 
	\begin{equation}\label{eq:190822primo}
		(\cap_{  i\in \mathbb{J}} \mathcal{B}_i)\cap \mathcal{D}_0\cap \mathcal{C}_0   \subseteq  \{ T_{K^*}\ge	Z_{I+1} 
		= T_{\lfloor cn\rfloor} \}, 
	\end{equation}
	where $\mathbb{J}=\{0,1,\cdots \overline{i}-1\}$,  with $\overline{i}$  defined in \eqref{eq:19092022primo},
	\begin{equation}\label{09042024}
		\mathcal{C}_0:=\{ T_{\terzotheokappazero }< T^B_{{\terzotheokappazerob }}\}
		=\{ N_B[\terzotheokappazero]< {\terzotheokappazerob }\}
		\quad
		\text{and}
		\quad
		\mathcal{D}_0:=\{ T_{{K^*}}\ge  T_{\terzotheokappazero }\}=\{{{K^*}}\ge \terzotheokappazero \}.
	\end{equation}
	\item[{\bf Step 3}.] We show 	
	\begin{equation}\label{eq:12092022quarto}
		(\cap_{i\in \mathbb{J}}\mathcal{A}_i)\cap \mathcal{C}_0\subseteq \{ 
		N(Z_{I+1})\ge N(Z_{I}) = 4^{I} \terzotheokappazero ,  N_B(Z_{I+1})\le 2^{I+1}{\terzotheokappazerob }  \}.
	\end{equation}
	\item[{\bf Step 4}.] We prove 
	\begin{equation*}
		(\cap_{i\in \mathbb{J}} (\mathcal{A}_i\cap \mathcal{B}_i))\cap \mathcal{C}_0\cap \mathcal{D}_0 \subseteq 
		\left\{   T_{{K^*}}\ge T_{\lfloor cn\rfloor},  \frac{N_B\big[{\lfloor cn \rfloor }\big]} {n }\le  2^{-\underline{i}+1}c \right\},
	\end{equation*}
	from which  we get claim
	\eqref{eq:11082022-2-mod},  provided that 
		$		\sum_{n\geq 1}\mathbb{P}([(\cap_{i\in \mathbb{J}} (\mathcal{A}_i\cap \mathcal{B}_i))\cap \mathcal{C}_0\cap \mathcal{D}_0]^c)
		<\infty.$
		\item[{\bf Step 5}.] We show that the latter infinite sum converges exploiting \eqref{eq-new-b-q>>g} and \eqref{eq-new-a-q>>g}. 
	\end{enumerate}	
	\subsubsection{Detailed proof of Part 1}  
	
	We accomplish Steps 1-5 outlined above.
	\\
	\\
	\noindent{\bf Step 1.}
	Recalling the definitions of $Z_{i+1}$ in  
	\eqref{eq:TsuCo}  and $I$ in
	\eqref{eq:20072022},
	it is rather immediate to verify that $Z_{I+j}=T_{\lfloor cn \rfloor}$, $\forall  j\in \mathbb{N}$,
	and
	{
	\begin{equation}\label{eq:19092022primo}
		\begin{gathered}
			\underline i\leq I<\overline i,
			\qquad\text{a.s.},\\
			\underline i:=
			\left\lceil
			\log_4\frac{\lfloor cn\rfloor}{h_0}
			\right\rceil-1,\\
			\overline i:=
			\left\lceil
			\log_4\frac{\lfloor cn\rfloor}{h_0}
			\right\rceil
			+
			\left\lceil
			\log_2\frac{\lfloor cn\rfloor}{h_B^{(0)}}
			\right\rceil.
		\end{gathered}
	\end{equation}
}
	
	for all $n$ sufficiently large (in order to guarantee that all the involved quantities are meaningful).
		Recalling that $K_i:=N(Z_i)$, 
		we also have
		\begin{equation}\label{eq:22122022terzo}
			\begin{split}	
				K_i \le \min(4^{i}\terzotheokappazero , \lfloor cn \rfloor)\quad \text{and}  \quad
				N_B(Z_{i})\le 2^{i}{\terzotheokappazerob };  \quad  \forall i\ge  0  \qquad   
				Z_i= T_{\lfloor cn \rfloor }, \;\; \forall  i \ge {	\overline i}.
			\end{split}	
		\end{equation}
		Therefore,  we will limit our analysis to the intervals  $[Z_i, Z_{i+1})$ with $0\le i <  {	\overline i}$.  As far as the definition of the events $\mathcal B_i$ in \eqref{eq:16122022secondo} is concerned,
		fix $\delta\in(0,1-c)$ and $\varepsilon\in(0,1/18)$, and set
		$\rho_\varepsilon:=
		\min\left\{
		(1-\varepsilon)\FH(1/8),
		\frac{1}{36}\log\left(\frac{1}{18\varepsilon}\right)
		\right\}>0.
		$
		We define
		\[
		\lambda_i :=n(1-\delta)-\min\{4^{i+1}\terzotheokappazero +2^{i+1}{\terzotheokappazerob } , cn\},
	\quad\text{and}\quad
		\phi_i:=\max\left\{18n \text{e}^{-\rho_\varepsilon4^i\terzotheokappazero p}, g \right\}.
		\]
		\\
		{\bf Step 2.} 
		We have
		\begin{align}
			\cap_{  i\in \mathbb{J}}  \mathcal{B}_i &=
			\cup_{  j\in \mathbb{J}}  \Big((\cap_{  i\in \mathbb{J}} \mathcal{B}_i )\cap\{I=j\}\Big)
			=	\cup_{  j\in \mathbb{J}} \Big( \big( ( \cap_{i\le j}  \mathcal{B}_i )\cap
			( \cap_{j<i<{\overline i}} \mathcal{B}_i  ) \big) \cap\{I=j\} \Big) \nonumber \\
			&\supseteq \cup_{  j\in \mathbb{J}}
			\Big( \big( \cap_{i\le j} \mathcal{B}_i  \big)\cap \big( \cap_{j<i<\overline{i}}\{I<i\}  \big)  \cap  \{I=j\}\Big)
			=\cup_{  j\in \mathbb{J}}\Big(  \big( \cap_{i\le j} \mathcal{B}_i \big)\cap\{I=j\}\Big),\label{eq:12092022secondo}
		\end{align}
		where the inclusion follows from  the relation
		$\mathcal{B}_i\supseteq \{ I<i \}$, whereas the last equality follows from the inclusion $\{I=j\} \subseteq \cap_{j<i<\overline{i}}\{I<i\} $.  
		Comparing the second and the last expressions,
		in \eqref{eq:12092022secondo},  we immediately have
		\begin{equation}\label{eq:12192022terzo}
			\cap_{  i\in \mathbb{J}}  \mathcal{B}_i =	\cup_{  j\in \mathbb{J}} \Big( (\cap_{  i\in \mathbb{J}} \mathcal{B}_i )\cap\{I=j\}\Big)=	\cup_{  j\in \mathbb{J}}\Big(  \big( \cap_{i\le j} \mathcal{B}_i \big)\cap\{I=j\}\Big).
		\end{equation}
		By the definition of $\mathcal{B}_i$,  we obtain
		\[
		\mathcal{B}_i 
		\subseteq\{Q^R_{k+1}>0\,\, \forall k: k \in [K_{i}, K_{i+1}), I\ge i \} \cup \{ I<i\}.
		\]
		Therefore
		\[
		( \cap_{i\le j} \mathcal{B}_i )  \cap \{I=j\}\subseteq
		\{  Q^R_{k+1}>0\,\, \forall k:\; k \in [K_{0},K_{I+1} ),  I=j  \}.
		\]
		Combining this with 
		\eqref{eq:12192022terzo}, 
		we have
		\[
		\cap_{i\in \mathbb{J}}  \mathcal{B}_i\subseteq\{  Q^R_{k+1}>0\,\, \forall k:\; k \in [K_{0},K_{I+1} )\}.
		\]
		Similarly,   we obtain
		\begin{equation}\label{111022-1}
			\cap_{j\le i}  \mathcal{B}_j\subseteq\{  Q^R_{k+1}>0\;\; \forall k \in [K_{0},\min\{ K_{i+1},K_{I+1} \})\}.
		\end{equation}
		Considering the intersection with the set $\mathcal{D}_0\cap \mathcal{C}_0$, we finally have
		\eqref{eq:190822primo}, 
		since,  by construction, $Q^R_{K^*+1}=0$. 
		\\
		\noindent{\bf Step 3.} 
		By \eqref{eq:TsuCo},  the definitions of $\mathcal{C}_0$ and $\mathcal{A}_i$,  and \eqref{eq:20072022},
		for any $\omega \in \mathcal{A}_{i}\cap \mathcal{C}_0\cap\{I(\omega)=j\} $, with $i<j$,  we have:
		$Z_{i+1}(\omega)=T_{4^{i+1} \terzotheokappazero }(\omega)$,
		$T_{4^{i+1} {\terzotheokappazero }}(\omega) \le T^B_{2^{i+1}{\terzotheokappazerob }} (\omega)$ and  $T_{4^{i+1} \terzotheokappazero } (\omega)< T_{\lfloor cn \rfloor} (\omega)$. Similarly,  for any $\omega \in \mathcal{A}_j\cap \mathcal{C}_0\cap\{I(\omega)=j\}$,
		we have: $T_{\lfloor cn \rfloor} (\omega) \le T_{4^{j+1} \terzotheokappazero }(\omega)$ and 
		$T_{\lfloor cn \rfloor}(\omega)  \le T^B_{2^{j+1} {\terzotheokappazerob }}(\omega)$.  In particular,  for $\omega \in \mathcal{A}_{j-1}\cap  \mathcal{A}_{j} \cap \mathcal{C}_0\cap\{I(\omega)=j\}$,  we obtain:
		$Z_{j}(\omega) =Z_I(\omega)= T_{4^{I} \terzotheokappazero }(\omega)  \le $
		$Z_{j+1}(\omega)=Z_{I+1}(\omega)=T_{\lfloor cn \rfloor}(\omega)\le  T^B_{2^{I+1}h^{(0)}_B}(\omega)$. 
		The claim \eqref{eq:12092022quarto} then follows by taking the union over all values $j$ that $I$ can assume.\\
		\noindent{\bf Step 4.}
			Combining this estimate with \eqref{eq:190822primo}, we obtain
			{
				On
				$\bigcap_{i\in\mathbb J}\mathcal A_i\cap\mathcal C_0$,
				the definition of $I$ and the events $\mathcal A_i$
				give
				\[
				I=\underline i.
				\]
				Indeed, $\underline i$ is the first index $i$
				for which
				$4^{i+1}h_0\geq\lfloor cn\rfloor$,
				and $\mathcal A_{\underline i}$ ensures that
				the corresponding black threshold is not reached
				before $T_{\lfloor cn\rfloor}$.
				
				Moreover, for all sufficiently large $n$,
				$h_B^{(0)}\leq h_0$, and
				\[
				4^{\underline i}h_0
				<
				\lfloor cn\rfloor
				\leq
				4^{\underline i+1}h_0.
				\]
				Hence, by \eqref{eq:12092022quarto},
				\[
				\frac{N_B[\lfloor cn\rfloor]}{n}
				\leq
				\frac{2^{\underline i+1}h_B^{(0)}}{n}
				\leq
				2^{1-\underline i}
				\frac{4^{\underline i}h_0}{n}
				\leq
				2^{1-\underline i}c.
				\]
				Since $h_0=o(n)$, we have
				$\underline i\to\infty$.
				Combining this estimate with
				\eqref{eq:190822primo}, we obtain
			
			}
			
		\begin{equation}\label{eq:080822}
			\{(\cap_{i\in \mathbb{J} } (\mathcal{A}_i\cap \mathcal{B}_i))\cap \mathcal{C}_0\cap \mathcal{D}_0\} \subseteq 
			\mathcal{T} := \left\{   T_{{K^*}}\ge T_{\lfloor cn\rfloor},  \frac{N_B\big[{\lfloor cn \rfloor }\big]} {n }\le  2^{-\underline{i}+1}c \right\}.
		\end{equation}
		We shall show later that
		\begin{equation}\label{eq:28072022I}
			\sum_{n\geq 1}\mathbb{P}([(\cap_{i\in \mathbb{J}} (\mathcal{A}_i\cap \mathcal{B}_i))\cap \mathcal{C}_0\cap \mathcal{D}_0]^c)=
			\sum_{n\geq 1}	\left[\mathbb{P}  (\cup_{i\in \mathbb{J}} (\mathcal{A}_i^c\cup \mathcal{B}_i^c)\cap\mathcal{C}_0\cap\mathcal{D}_0)+\mathbb{P}(\mathcal{C}_0^c\cup \mathcal{D}_0^c)\right]
			<\infty.
		\end{equation}
		Therefore by the Borel--Cantelli lemma and \eqref{eq:080822},  we obtain: 
		$\mathbb{P}(\limsup \mathcal{T}^c)
		=0$, 
		which  implies \eqref{eq:11082022-2-mod}.  
		 Indeed by construction
		\[
		\limsup \frac{N_B[K^*]-N_B[\lfloor cn \rfloor ] }{n}\le \limsup \frac{\max(0, K^*-  cn )}{n} \le 1-c   \qquad \forall c \in (0,1).
		\]
		Applying the preceding estimate to the countable sequence
		$c_m:=1-m^{-1}\uparrow1$ yields $N_B[K^*]/n\to0$ a.s.\\
		\noindent{\bf Step 5.} To prove \eqref{eq:28072022I},  note that
		$	\mathbb{P}(\liminf\mathcal{C}_0\cap \mathcal{D}_0)=1$, which is 
		 an immediate consequence of Theorem \ref{lemma-nuovo}$(ii)$ for  the case $g \ll q \ll p^{-1}$ (recall that for $n$ sufficiently  large   $h^{(0)}_B := \lfloor p^{-1} \rfloor > \lfloor (g_B(\kappa_{\bm g})+\varepsilon)q\rfloor$),
		and of Theorem  \ref{prop:kminT} $(i)$  together with  Corollary \ref{cor:190823} for the remaining cases.
		Therefore,   by the second Borel--Cantelli lemma it follows that
		$\sum_n \mathbb{P}(\mathcal{C}_0^c\cup \mathcal{D}_0^c)<\infty$
		since the events $\mathcal{C}_0^c\cup \mathcal{D}_0^c$  are independent for different values of $n$. Thus  to establish  \eqref{eq:28072022I} it remains  to show
		\begin{align}\label{eq:280722II}
			\sum_{n\geq 1}\mathbb{P}  (\cup_{i\in \mathbb{J}} (\mathcal{A}_i^c\cup \mathcal{B}_i^c)\cap\mathcal{C}_0\cap\mathcal{D}_0) &
			<	\infty.
		\end{align}
		To this aim, note that proceeding similarly to \eqref{120923}, we have
		\begin{align*}
			\mathbb{P}\Big(\cup_{i \in \mathbb{J}} (\mathcal{A}_i^c\cup \mathcal{B}_i^c) \cap \mathcal{C}_0\cap \mathcal{D}_0\Big)
			&\le
			\sum_{i\in \mathbb{J} } 	\mathbb{P}\left( \Big(   \mathcal{B}_i^c\cap\big(\cap_{j<i} \mathcal{A}_j \big)\Big) \cap \mathcal{C}_0\right) \nonumber \\
			&
			+\sum_{i\in \mathbb{J}} \mathbb{P}\left( \big(  \mathcal{A}_i^c\cap  \mathcal{B}_i\big) \cap 
			\big( 
			\cap_{j<i} (\mathcal{A}_j\cap \mathcal{B}_j) \big) \Big) \cap \mathcal{C}_0\cap \mathcal{D}_0\right).
			\label{eq-new-qggp-28June-2}
		\end{align*}
		Since ${\overline i}=O(\log_2 (np))$, relation  \eqref{eq:280722II} follows from  \eqref{eq:bootcond}, provided that we can show
		\begin{equation} 	\label{eq-new-b-q>>g} 
			\sup_{ i\in \mathbb{J} }\mathbb{P} (\mathcal{B}^c_{i}\cap (\cap_{ j<i} \mathcal{A}_j) \cap \mathcal{C}_0)
			\leq
			n^3 \left(\text{e}^{-n\left(1- \frac{\delta}{2}\right)\FH\left( \frac{1-\delta}{1-\delta/2 } \right)}
			+ \text{e}^{-\frac{g}{2}\log 8 }\right),  
		\end{equation}
		
		and
		\begin{equation}
			\sup_{ i\in \mathbb{J}}\mathbb{P}\left(\mathcal{A}^c_{i}\cap \mathcal{B}_i\cap  (\cap_{ j< i}  \mathcal{A}_j\cap\mathcal{B}_j)\cap \mathcal{C}_0\cap \mathcal{D}_0 \right)
			\leq
			\text{e}^{ -\frac{{\terzotheokappazerob }}{2} \log(10)}.
			\label{eq-new-a-q>>g}
		\end{equation}
		for all $n$ large enough.
		 For $i=0$,  we conventionally  set: $ (\cap_{0\le j\le -1}\, \mathcal{A}_j):=\Omega$.
		To conclude the proof of the theorem,  it remains to verify \eqref{eq-new-b-q>>g} and  \eqref{eq-new-a-q>>g}.  This  will be accomplished  in the Parts 2 and 3 of the proof.
		
		\subsubsection{Part 2} 		
		We break down the proof of this part in three steps.
		
		\begin{enumerate}
			\item[{\bf Step 1}.] Let $\terzotheokappazero $ and ${\terzotheokappazerob }$ be as  in \eqref{100124}.
			We define the sets  $	\mathcal{E}_{i-1}$, 	$\mathcal{M}_{i}$, 	$\mathcal{E}^{(k)}_{i-1}$ and 
			$\mathcal{M}^{(k)}_{i}$
			and establish some set inclusions concerning these events,  namely  relations
			\eqref{081022-5} and \eqref{081022-3}.
			\item[{\bf Step 2}.]
			{	We apply Lemma~\ref{le:200523} using deterministic lower and upper
			bounds for $N_R[k]$ and $N_B[k]$, respectively. This yields
			stochastic bounds for the corresponding suprathreshold counts in
			terms of binomial random variables evaluated at a fixed configuration.
			These bounds are then used to derive the joint tail estimates
			\eqref{081022-1} and~\eqref{081022-2}. 
		}
			
			\item[{\bf Step 3}.] We provide an upper bound on the probability $\mathbb{P}(\mathcal{B}_i^c\cap\mathcal{E}_{i-1})$,  from which the claim immediately follows.
		\end{enumerate}	
		
		\subsubsection{Detailed proof of Part 2}
		
	We now proceed to carry out Steps 1–3 as outlined above.\\
		
		\noindent{\bf Step 1.}
		Setting for $i>0$
		\begin{equation}\label{defEi}
			\mathcal E_{i-1}
			:=
			\left(\bigcap_{0\le j<i}\mathcal A_j\right)
			\cap\mathcal C_0
			\cap\{I\ge i\},  \quad \text{with} \quad 	\mathcal E_{-1}=\mathcal{C}_0
		\end{equation}	
		We have
		\begin{equation}\label{081022-5}
			\mathcal{E}_{i-1}\subseteq\mathcal{M}_{i}:=	\{N_R[h]\geq 4^{i}\terzotheokappazero  -2^{i+1}{\terzotheokappazerob }\,\, \text{and} \,\, N_B[h]\leq 2^{i+1}{\terzotheokappazerob }\,\,\, \forall h\in [K_i, K_{i+1})\}.
		\end{equation}
		Indeed,  if $\omega\in\mathcal{E}_{i-1}$,  then 
		\begin{equation}\label{27102022}
			\omega\in\mathcal{A}_{i-1}\cap \mathcal{C}_0\cap\{I \geq i\}\subseteq\{   
			Z_{i}=T_{4^{i} \terzotheokappazero  } \}= \{K_i = 4^{i} \terzotheokappazero \}, 
		\end{equation}
		which implies
		$N[h](\omega)=N_R[h](\omega)+N_B[h](\omega)\ge 4^{i} \terzotheokappazero  $, for any $h \in [K_i(\omega), K_{i+1}(\omega))$.  Furthermore,  by definition (see \eqref{eq:TsuCo})  $Z_{i+1}(\omega)\le T^B_{ 2^{i+1}{\terzotheokappazerob } }(\omega)$,  
		which yields
		$N_B[h](\omega)
		\le 2^{i+1}  {\terzotheokappazerob } $, for any $h\in [K_i(\omega), K_{i+1}(\omega) )$. 
		Moreover, by  \eqref{081022-5} it follows 
		\begin{align}
			\mathcal{E}_{i-1}^{(k)}&:=\mathcal{E}_{i-1}\cap \{ k\in [K_i, K_{i+1}) \}\subseteq  
			\mathcal{M}_{i}\cap \{ k\in [K_i, K_{i+1}) \}\nonumber \\ 
			&\subseteq \mathcal{M}_{i}^{(k)} := 
			\{  N_R[k]    \ge 4^{i}\terzotheokappazero  - 2^{i+1}{\terzotheokappazerob }, N_B[k]  \le 2^{i+1} {\terzotheokappazerob } \}.
			\label{081022-3}
		\end{align}
		\noindent{\bf Step 2.}
		Setting ${ \underline{k}^{(i)}_R}:=  4^{i}\terzotheokappazero  - 2^{i+1}  {\terzotheokappazerob }$ and $
		{ \overline{k}^{(i)}_B}= 2^{i+1}  {\terzotheokappazerob }$,
		and applying  Lemma \ref{le:200523}  with $\ell_R=\underline{k}^{(i)}_R$ and $\ell_B= \overline{k}^{(i)}_B $, we have  $\forall\,x>0$
			{ 
			\begin{equation}\label{eq:10102022primo}
				\begin{split}
			&\mathbb{P}\left(|\mathcal{S}_R[k]|\le x,  \mathcal{M}_{i}^{(k)}\right)\le \mathbb{P}\left( \text{Bin}(n_W,\pi_R(\underline{k}^{(i)}_R ,  \overline{k}^{(i)}_B  ))\le x\right),  \\ 	
		&\mathbb{P}\left(	|\mathcal{S}_B[k ]|\ge x, \mathcal{M}_{i}^{(k)}   \right) \le \mathbb{P}\left(\text{Bin}(n_W,\pi_B(\underline{k}^{(i)}_R , \overline{k}^{(i)}_B))\ge x \right).
		\end{split}
				\end{equation}
			}
		Note that,  for any $z\ge r$ and any $S\in\{R,B\}$,  it holds
		\begin{align}
			\pi_S( k_R,     k_{B}   )&= 
			\mathbb{P}( \text{Bin}(k_S+a_S,p)-\text{Bin}(k_{\overline{S}}+a_{\overline{S}},p)\ge r)\nonumber\\ 
			&\ge 
			\mathbb{P}( \text{Bin}(k_S+a_S,p)\ge z,  \text{Bin}(k_{\overline{S}}+a_{\overline{S}},p)\le z-r)\nonumber\\
			&	\ge 1- \mathbb{P}( \text{Bin}(k_{S}+a_{S},p)< z) -\mathbb{P}( \text{Bin}(k_{\overline{S}}+a_{\overline{S}},p)> z-r). 
			\label{eq:10102022quarto}
		\end{align}
		Moreover,  for $n$ sufficiently large,  assuming ${\terzotheokappazerob }/\terzotheokappazero < ({\terzotheokappazerob }+a_B)/\terzotheokappazero  <\varepsilon/2$, we have
		\[
		\begin{split}
			&\mathbb{E}[\text{Bin}(\underline{k}^{(i)}_R+a_R,p)]\ge \mathbb{E}[\text{Bin}(\underline{k}^{(i)}_R,p)]\ge
			4^{i}\terzotheokappazero p \left(1- \frac{2{\terzotheokappazerob }}{2^i\terzotheokappazero }\right)\ge 4^{i}\terzotheokappazero p \left(1-\varepsilon \right),\quad\text{and} \\
			&\mathbb{E}[\text{Bin}(\overline{k}^{(i)}_B+a_B,p)]
			\le  2^{i+1} ({\terzotheokappazerob }+a_B) p.
		\end{split}
		\]
		Therefore, taking $z=4^i\terzotheokappazero p/9$ in
		\eqref{eq:10102022quarto} and applying the concentration inequalities
		in Appendix~\ref{Penrose}, we obtain, uniformly in $i\ge0$,
		\begin{equation}\label{eq:20092022primo}
			\begin{split}
				\pi_R(\underline{k}^{(i)}_R,\overline{k}^{(i)}_B)
				&\ge 1-2\text{e}^{-\rho_\varepsilon4^i\terzotheokappazero p},\\
				\pi_B(\underline{k}^{(i)}_R,\overline{k}^{(i)}_B)
				&\le 2\text{e}^{-\rho_\varepsilon4^i\terzotheokappazero p}.
			\end{split}	
		\end{equation}
		Here we used
		$\log(2^{i-1}/(9\varepsilon))\ge\log(1/(18\varepsilon))$.
		{ Since the error term in \eqref{eq:20092022primo} is non-increasing
			in $i$, by choosing $\kappa$ sufficiently large when $q=p^{-1}$ and
			using $n_W/n\to1$, we obtain, uniformly in $i\ge0$ and for all
			sufficiently large $n$,
			\[
			n_W\pi_R(\underline{k}^{(i)}_R,\overline{k}^{(i)}_B)
			\ge n(1-\delta/2).
			\]}
		Now , from \eqref{eq:10102022primo}, by applying again the concentration inequalities reported in Appendix \ref{Penrose},
		for any $i$ and all $n$ large enough,  we obtain
			{
			\begin{equation}\label{081022-1}
				\mathbb{P}(|\mathcal{S}_R[k]|\le (1-\delta)n, \mathcal{M}_{i}^{(k)} )\le 
				\mathbb{P}\left( \text{Bin}(n_W,\pi_R(\underline{k}^{(i)}_R ,  \overline{k}^{(i)}_B  ))\le   (1-\delta)n, \right)\le 				
				 \text{e}^{-n\left(1- \frac{\delta}{2}\right)\FH\left( \frac{1-\delta}{1-\frac{\delta}{2}} \right)}.
			\end{equation}
		}

		Similarly,  exploiting  \eqref{eq:10102022primo} and \eqref{eq:20092022primo}, for any $i$ and all $n$ sufficiently large, we have
{
		\begin{equation}\label{081022-2}
			\mathbb{P}\left(|\mathcal{S}_B[k]|\ge \phi_i, \mathcal{M}_{i}^{(k)} \right) \le 
			\mathbb{P}\left(\text{Bin}(n_W,\pi_B(\underline{k}^{(i)}_R , \overline{k}^{(i)}_B))\ge \phi_i   \right)
			\leq\text{e}^{-
				\frac{\phi_i}{2}\log 8}.
		\end{equation}
}
where
$\overline{\mu}^B_i:=2n\text{e}^{-\rho_\varepsilon4^i\terzotheokappazero p}
\ge n_W\pi_B(\underline{k}_R^{(i)},\overline{k}_B^{(i)})$
and $\phi_i:=\max(9\overline{\mu}^B_i,g)$.\\
		\noindent{\bf Step 3.}
		By  \eqref{eq:16122022secondo},
		for any $i$ and all $n$ large enough,  we obtain
		
		\begin{align}
			\mathbb{P}(\mathcal{B}_i^c\cap&\mathcal{E}_{i-1})
			=\mathbb{P}\left(\bigcup_{k} \Big(\Big\{k \in [K_i, \, K_{i+1}),\,\, Q^R_{k+1}\le \lambda_i \text{ or } 
			Q^B_{k+1}> \phi_i\Big\}
			\cap \mathcal{E}_{i-1}\Big) \right)\nonumber\\
			&=
			\mathbb{P}\left(\bigcup_{k}\Big(\Big\{ Q^R_{k+1}\le \lambda_i \text{ or } 
			Q^B_{k+1}> \phi_i\Big\}\cap \mathcal{E}_{i-1}^{(k)}\Big) \right)\nonumber\\
			& 
			=\sum_{k_i, k_{i+1}}
			\mathbb{P}\Big( 	    \{K_i =k_i, K_{i+1}=k_{i+1} \} 
			\bigcap \big( \bigcup_{k=k_i}^{k_{i+1}-1} \big[ (\{ Q^R_{k+1}\le \lambda_i\}\cup\{ 
			Q^B_{k+1}> \phi_i \})  \cap \mathcal{E}_{i-1}^{(k)} \big] \big) \Big) \nonumber\\
			&\stackrel{(a)}{\le}\sum_{k_i,k_{i+1}}
			\sum_{k=k_i}^{k_{i+1}-1} \mathbb{P}
			\left(
			\{  K_i =k_i,   K_{i+1}=k_{i+1} \}  \bigcap\Big(\{Q^R_{k+1}\le \lambda_i\} \cup\{ 
			Q^B_{k+1}> \phi_i\}\Big)\cap \mathcal{M}_{i}^{(k)}\right) \nonumber\\ 
			&\le \sum_{k_i,k_{i+1}} \sum_{k=k_i}^{k_{i+1}-1}\mathbb{P}(
			(\{Q^R_{k+1}\le \lambda_i\}\cup\{ 
			Q^B_{k+1}> \phi_i \})\cap \mathcal{M}_{i}^{(k)})\nonumber\\     
			&\stackrel{(b)}{\le}\sum_{k_i,k_{i+1}}\sum_{k=k_i}^{k_{i+1}-1}
			\mathbb{P}(\left(\{|\mathcal{S}_R[k]| \le \lambda_i +k\}\cup\{ 
			|\mathcal{S}_B[k]|> \phi_i \}\right)\cap  \mathcal{M}_{i}^{(k)}  )
			\nonumber\\
			& \stackrel{(c)}{\le} n^3 \left(e^{-n\left(1- \frac{\delta}{2}\right)\FH\left( \frac{1-\delta}{1-\frac{\delta}{2} } \right)}
			+ \text{e}^{-\frac{g}{2}\log 8 }\right),\label{eq:10102022sei}
		\end{align}
		where the indices $k_i$ and $k_{i+1}$ in the sums range over the support of $K_i$ and $K_{i+1}$, respectively.
		Here, inequality $(a)$ follows from \eqref{081022-3},   $(b)$  from \eqref{180523-1}, and  $(c)$ combines 
		\eqref{081022-1} and \eqref{081022-2}
		(using $\lambda_i+k\le (1-\delta)n$),  the union bound, the fact that  $K_i $, for every $ i$,  takes values  in $\{0,\ldots, n_W \}$, and the monotonicity of $\phi_i$ in $i$.
		
		Finally,  we note that relation \eqref{eq-new-b-q>>g} follows immediately. Indeed, after recalling \eqref{defEi} and observing that
		$\mathcal{B}_i^c\cap \{I<i \}=\emptyset$
		(by \eqref{eq:16122022secondo}), we have
		\[
		\mathbb{P}(\mathcal{B}_i^c\cap  (\cap_{j< i} \mathcal{A}_j)\cap \mathcal{C}_0)=\mathbb{P}(\mathcal{B}_i^c\cap (\cap_{ j< i} \mathcal{A}_j)\cap\mathcal{C}_0\cap \{I \ge i\})=\mathbb{P}(\mathcal{B}_i^c\cap \mathcal{E}_{i-1} ). 
		\]

		{
		\subsubsection{Part 3}
		
		We now prove~\eqref{eq-new-a-q>>g}. The proof is organized into four steps.
		
		In Step~1, we establish the preliminary relations~\eqref{270525}
		and~\eqref{270525-2}.
		
		In Step~2, we introduce suitable truncated versions of the
		$B$-activation increments. Their conditional success probabilities
		are uniformly bounded from above by a deterministic quantity, which
		allows us to stochastically dominate their partial sums by suitable
		binomial random variables.
		
		In Step~3, we introduce the events $\mathcal G_i$ and
		$\widetilde{\mathcal B}_i$ and show that an early termination of the
		$i$-th interval requires a sufficiently large number of
		$B$-activations.
		
		Finally, in Step~4, we combine the previous relations with the
		binomial concentration inequalities in Appendix~\ref{Penrose} to prove~\eqref{eq-new-a-q>>g}.

		\subsubsection{Detailed proof of Part 3} We now proceed to carry out Steps 1–4 as outlined above.\\
			\noindent{\bf Step 1.} 
		Since  $Q^R_{{K^*+1}}=0$,   it follows from \eqref{09042024},
		\eqref{111022-1} and \eqref{defEi}  that
			\begin{equation}
			\label{270525}
			\mathcal D_0\cap\left(\bigcap_{j\le i}\mathcal B_j\right)\cap\mathcal E_{i-1}	\subseteq
			\mathcal C_0\cap\mathcal D_0\cap \left(\bigcap_{j\le i}\mathcal B_j\right)\cap\{I\ge i\}
			\subseteq		\{T_{K^*}\ge Z_{i+1}\}.
		\end{equation}
		Moreover, by \eqref{26102022}, we obtain
		\begin{equation}\label{270525-2}
			N_S[K_{i+1}]=N_S[K_{i}]+ \sum_{k=1}^{K_{i+1}-K_i}M^S_{K_i+k}.
		\end{equation}

		\noindent
		{\bf Step 2.}
		Set
		\begin{equation}
			\label{eq:theta-i}
			\theta_i
			:=
			\frac{\phi_i}{\lambda_i+\phi_i}.
		\end{equation}
		For every $m\ge1$, define
		\begin{equation}
			\label{eq:truncated-M}
			\widehat M^B_{i,m}
			:=
			M^B_{4^ih_0+m}
			\ind_{\{U^B_{4^ih_0+m}\le\theta_i\}}.
		\end{equation}
		
		By Proposition~\ref{M-U},
		\[
		\mathbb P
		\left(
		M_k^B=1
		\,\middle|\,
		\mathcal H_{k-1}
		\right)
		=
		U_k^B.
		\]
		Since $U_k^B$ is $\mathcal H_{k-1}$-measurable, it follows that
		\begin{align}
						\mathbb P\left(\widehat M^B_{i,m}=1		\,\mid\,		\mathcal H_{4^ih_0+m-1}		\right)	
			& =
			\mathbb E\left [\widehat M^B_{i,m}\, \mid \,   	\mathcal H_{4^ih_0+m-1}	  \right]
			=
			\ind_{\{U^B_{4^ih_0+m}\le\theta_i\}}
			U^B_{4^ih_0+m}
			\le
			\theta_i .
			\label{240826}
		\end{align}
		
		Consequently, for every deterministic $h\ge1$,
		\begin{equation}
			\label{eq:stoch-dom-M}
			\left.
			\sum_{m=1}^{h}\widehat M^B_{i,m}
			\,\right|\,
			\mathcal H_{4^ih_0}
			\; \leq_{\mathrm{st}}\;
			\text{Bin}(h,\theta_i).
		\end{equation}
		{
		To prove~\eqref{eq:stoch-dom-M}, let
		\[
		S_m:=\sum_{\ell=1}^m \widehat M^B_{i,\ell},
		\qquad S_0:=0,
		\]
		and let $\varphi:\mathbb{N}\to\mathbb{R}$ be an arbitrary bounded
		non-decreasing function. Set
		\[
		P_m:=
		\mathbb P\left(
		\widehat M^B_{i,m}=1
		\,\middle|\,
		\mathcal H_{4^ih_0+m-1}
		\right).
		\]
		By~\eqref{240826}, $P_m\le\theta_i$. Moreover, since $S_{m-1}$ is
		$\mathcal H_{4^ih_0+m-1}$-measurable, we have
		\begin{align*}
			\mathbb E\left[
			\varphi(S_m)
			\,\middle|\,
			\mathcal H_{4^ih_0+m-1}
			\right]
			&=
			(1-P_m)\varphi(S_{m-1})
			+P_m\varphi (S_{m-1}+1)\\
			&\le
			(1-\theta_i)\varphi(S_{m-1})
			+\theta_i\varphi(S_{m-1}+1),
		\end{align*}
		where the inequality follows because $\varphi$ is non-decreasing and
		$P_m\le\theta_i$.
		
		Define the operator
		\[
		(\mathsf T_{\theta_i}\varphi)(s)
		:=
		(1-\theta_i)\varphi(s)
		+\theta_i\varphi(s+1).
		\]
		Since $\mathsf T_{\theta_i}\varphi$ is non-decreasing whenever
		$\varphi$ is non-decreasing, iterating the preceding inequality and
		conditioning successively yields
		\[
		\mathbb E\left[
		\varphi(S_h)
		\,\middle|\,
		\mathcal H_{4^ih_0}
		\right]
		\le
		(\mathsf T_{\theta_i}^{\circ h}\varphi)(0)
		=
		\mathbb E\left[
		\varphi\bigl(\text{Bin}(h,\theta_i)\bigr)
		\right].
		\]
		Since this holds for every bounded non-decreasing function $\varphi$,
		it proves~\eqref{eq:stoch-dom-M}.

		In particular, for every real $x\ge0$,
		\begin{equation}
			\label{eq:tail-dom-M}
			\mathbb P
			\left(
			\sum_{m=1}^{h}\widehat M^B_{i,m}
			\ge x
			\,\middle|\,
			\mathcal H_{4^ih_0}
			\right)
			\le
			\mathbb P
			\left(
			\text{Bin}(h,\theta_i)\ge x
			\right).
		\end{equation}

		\noindent
		{\bf Step 3.}
		From	\eqref{defEi} and \eqref{27102022}, we have
		\[
		\mathcal E_{i-1}	\subseteq		\mathcal A_{i-1}	\cap		\mathcal C_0 \cap \{I\ge i\}
		\subseteq \{K_i=4^ih_0\}.
		\]
		Therefore, define
		\begin{equation}
			\label{240826-2}
			\mathcal G_i	:=			\mathcal D_0		\cap		\left(\bigcap_{j<i}\mathcal B_j\right)			\cap
			\mathcal E_{i-1}			\subseteq	\{K_i=4^ih_0\}.
		\end{equation}
		
		Recalling the definition of $\mathcal B_i$ and observing that,
		by~\eqref{270525},
		\[
		\mathcal G_i\cap\mathcal B_i
		\subseteq
		\{T_{K^*}\ge Z_{i+1}\}
		\subseteq
		\left\{
		Q^B_{k+1}\ge0,\quad
		\forall k\in[K_i,K_{i+1})
		\right\},
		\]
		we obtain
		\begin{equation}
			\label{240826-3}
			\mathcal G_i\cap\mathcal B_i
			=
			\mathcal G_i\cap\mathcal B_i\cap\{I\ge i\}
			\subseteq
			\mathcal G_i\cap\widetilde{\mathcal B}_i,
		\end{equation}
		where
		\begin{equation}
			\label{eq:Btilde}
			\widetilde{\mathcal B}_i
			:=
			\left\{
			U^B_{k+1}
			<
			\frac{\phi_i}{\lambda_i+\phi_i},
			\quad
			\forall k\in[K_i,K_{i+1})
			\right\}.
		\end{equation}
	{ For indices such that $\mathcal G_i\neq\emptyset$, 
		\eqref{240826-2} and $\mathcal G_i\subseteq\mathcal E_{i-1}$ imply
		\[
		K_i=4^ih_0<\lfloor cn\rfloor,
		\]
		where the strict inequality follows from $I\ge i$ and the definition
		of $I$. Hence,
		\[
		\delta_{\max}^{(i)}
		:=
		\min\{4^{i+1}h_0,\lfloor cn\rfloor\}-4^ih_0
		\]
		is a strictly positive integer. If $\mathcal G_i=\emptyset$, the
		required probability estimate is trivial.}
		
		Let $0\le h<\delta_{\max}^{(i)}$, we have. 
		\begin{align}
			&
			\{K_{i+1}=h+4^ih_0\}
			\cap\mathcal G_i
			\nonumber\\
			&\qquad=
			\{K_{i+1}-K_i=h\}
			\cap\mathcal G_i
			\nonumber\\
			&\qquad\subseteq
			\left\{
			N_B[K_{i+1}]
			=
			N_B[K_i+h]
			=
			2^{i+1}h_B^{(0)}
			\right\}
			\cap\mathcal G_i
			\nonumber\\
			&\qquad=
			\left\{
			N_B[K_i+h]-N_B[K_i]
			=
			2^{i+1}h_B^{(0)}-N_B[K_i],
			\;
			N_B[K_i]\le2^ih_B^{(0)}
			\right\}
			\cap\mathcal G_i
			\nonumber\\
			&\qquad\subseteq
			\left\{
			N_B[K_i+h]-N_B[K_i]
			\ge
			2^ih_B^{(0)}
			\right\}
			\cap\mathcal G_i
			\nonumber\\
			&\qquad=
			\left\{
			\sum_{m=1}^{h}
			M^B_{4^ih_0+m}
			\ge
			2^ih_B^{(0)}
			\right\}
			\cap\mathcal G_i.
			\label{240826-4}
		\end{align}
		Here, the first inclusion follows from the definition of $Z_{i+1}$,
		whereas the last equality follows from~\eqref{270525-2}.
		
		On the event
		\[
		\widetilde{\mathcal B}_i
		\cap
		\{K_{i+1}=h+4^ih_0\}
		\cap
		\mathcal G_i,
		\]
		we have
		\[
		U^B_{4^ih_0+m}<\theta_i,
		\qquad
		1\le m\le h,
		\]
		and therefore, by~\eqref{eq:truncated-M},
		\[
		\widehat M^B_{i,m}
		=
		M^B_{4^ih_0+m},
		\qquad
		1\le m\le h.
		\]
		Combining this observation with~\eqref{240826-3} and~\eqref{240826-4},
		we obtain
		\begin{align}
			&
			\{K_{i+1}<\delta_{\max}^{(i)}+K_i\}
			\cap
			\mathcal B_i
			\cap
			\mathcal G_i
			\nonumber\\
			&\qquad\subseteq
			\bigcup_{h=1}^{\delta_{\max}^{(i)}-1}
			\left\{
			\sum_{m=1}^{h}
			\widehat M^B_{i,m}
			\ge
			2^ih_B^{(0)}
			\right\}
			\cap
			\mathcal G_i
			\nonumber\\
			&\qquad\subseteq
			\left\{
			\sum_{m=1}^{\delta_{\max}^{(i)}-1}
			\widehat M^B_{i,m}
			\ge
			2^ih_B^{(0)}
			\right\}
			\cap
			\mathcal G_i,
			\label{240826-5}
		\end{align}
		where the last inclusion follows from the non-negativity of the
		variables $\widehat M^B_{i,m}$.

		\noindent
		{\bf Step 4.}
		From the definitions of $Z_{i+1}$ and $\mathcal A_i$, it follows
		that
		\begin{equation}\label{240826-6}
		\mathcal A_i^c	=	\left\{		K_{i+1}	<	\min\{4^{i+1}h_0,\lfloor cn\rfloor\}\right\}		\subseteq
		\{I\ge i\}.
		\end{equation}
		Therefore, by the definitions of $\mathcal E_{i-1}$ and
		$\mathcal G_i$,
		\begin{align}
			&
			\mathcal A_i^c
			\cap\mathcal B_i
			\cap
			\left(
			\bigcap_{j<i}
			(\mathcal A_j\cap\mathcal B_j)
			\right)
			\cap\mathcal C_0
			\cap\mathcal D_0
			\nonumber\\
			&\qquad=
			\mathcal A_i^c
			\cap\mathcal B_i
			\cap\mathcal G_i.
			\label{eq:part3-id}
		\end{align}
		{In view of \eqref{eq:part3-id}, it remains to bound
			$\mathbb P(\mathcal A_i^c\cap\mathcal B_i\cap\mathcal G_i)$.
			To apply the conditional estimate \eqref{eq:tail-dom-M}, observe that
			$
			\mathcal G_i\in\mathcal H_{4^ih_0}.
			$
			Indeed, on $\mathcal E_{i-1}$, the events
			$\mathcal A_0,\ldots,\mathcal A_{i-1}$ and $\{I\ge i\}$ imply
			recursively that $K_j=4^jh_0$ for $0\le j\le i$. Consequently, for
			$j<i$, the event $\mathcal B_j$ is determined by the variables
			$Q^R_{h+1}$ and $Q^B_{h+1}$ with
			$h<4^{j+1}h_0\le4^ih_0$. The remaining events entering the definition
			of $\mathcal G_i$ are also determined by the embedded chain up to
			index $4^ih_0$, proving the claimed measurability.
			
			Therefore, \eqref{240826-5}, \eqref{240826-6}, and
			\eqref{eq:tail-dom-M} yield}
		\begin{align}
			&
			\mathbb P	\left(		\mathcal A_i^c		\cap		\mathcal B_i			\cap	\mathcal G_i\right)		\nonumber\\
			&\qquad\le \mathbb P\left(\left\{			\sum_{m=1}^{\delta_{\max}^{(i)}-1}			\widehat M^B_{i,m}
			\ge 	2^ih_B^{(0)} \right\}	\cap			\mathcal G_i	\right)			\nonumber\\
			&\qquad\le
			\mathbb P	\left(\text{Bin}\left(		\delta_{\max}^{(i)}-1,\theta_i	\right)
			\ge		2^ih_B^{(0)}	\right)			\mathbb P(\mathcal G_i).
			\label{240826-7}
		\end{align}{
		Since $\delta_{\max}^{(i)}\le3\cdot4^ih_0$, stochastic monotonicity
		in the number of trials gives
		\[
		\mathbb P\left(
		\operatorname{Bin}(\delta_{\max}^{(i)}-1,\theta_i)
		\ge2^ih_B^{(0)}
		\right)
		\le
		\mathbb P\left(
		\operatorname{Bin}(3\cdot4^ih_0,\theta_i)
		\ge2^ih_B^{(0)}
		\right).
		\]
		For every $i$ such that $\mathcal G_i\neq\emptyset$, we have
		$4^ih_0<cn$ and $\lambda_i\ge n(1-\delta-c)$. Hence, denoting by
		\[
		\mu_i:=3\cdot4^ih_0\theta_i
		\]
		the mean of the binomial random variable above, the definitions of
		$\theta_i$ and $\phi_i$ yield, uniformly over such $i$,
		\[
		\frac{\mu_i}{2^ih_B^{(0)}}
		\le
		C\frac{2^ih_0}{h_B^{(0)}}
		\left(
		e^{-\rho_\varepsilon4^ih_0p}+\frac gn
		\right).
		\]
		The first term is uniformly negligible since
		\[
		\frac{2^ih_0}{h_B^{(0)}}e^{-\rho_\varepsilon4^ih_0p}
		\le
		\frac{4^ih_0p}{ph_B^{(0)}}
		e^{-\rho_\varepsilon4^ih_0p},
		\]
		where $ph_B^{(0)}$ is bounded away from zero. In the regimes
		$g\ll q\ll p^{-1}$ and $q\gg p^{-1}$, we have $h_0p\to\infty$;
		when $q=p^{-1}$, the same term can be made arbitrarily small by
		choosing $\kappa$ sufficiently large. Moreover,
		\[
		\frac{2^ih_0}{h_B^{(0)}}\frac gn
		\le
		\sqrt c\,\frac{g\sqrt{h_0}}{h_B^{(0)}\sqrt n}
		=o(1)
		\]
		by \eqref{100124} and $g\ll q\ll n$. Therefore,
		$2^ih_B^{(0)}\ge10\mu_i$ for all sufficiently large $n$, uniformly
		in $i$, and the binomial upper-tail bound gives
		\[
		\mathbb P\left(
		\operatorname{Bin}(3\cdot4^ih_0,\theta_i)
		\ge2^ih_B^{(0)}
		\right)
		\le
		e^{-2^{i-1}h_B^{(0)}\log(10)}.
		\]
		Combining this estimate with \eqref{eq:part3-id} and
		\eqref{240826-7}, and using $\mathbb P(\mathcal G_i)\le1$, yields
		\[
		\mathbb P\left(
		\mathcal A_i^c\cap\mathcal B_i
		\cap\bigcap_{j<i}(\mathcal A_j\cap\mathcal B_j)
		\cap\mathcal C_0\cap\mathcal D_0
		\right)
		\le
		e^{-\frac12h_B^{(0)}\log(10)}.
		\]
		Taking the supremum over $i\in\mathbb J$ proves
		\eqref{eq-new-a-q>>g}.
		}

	}

\subsection{Proof of Theorem \ref{prop:supercritical}}
Theorem \ref{prop:supercritical} follows directly from Theorems \ref{thm:Lemma-new-2} (case $q=g$)
and \ref{thm:Lemma-new-3} (case $q\gg g$).
\section{Conclusions}

We introduced and analyzed a competitive extension of bootstrap percolation
on the Erd\H{o}s--R\'enyi random graph, in which two irreversible activation
processes evolve concurrently. Our results reveal a sharp qualitative
dichotomy. In the sub-critical regime, competition is asymptotically
negligible at first order, whereas in the super-critical regime the initial
advantage of one process is amplified: the dominant process almost
percolates the graph, while the competing process remains confined to a
smaller scale.

Natural extensions include the symmetric case $\alpha_R=\alpha_B$, the
qualitatively different regime $r=1$, and the analysis of competing
bootstrap processes on more heterogeneous random graph models.

\appendix
\section{Proof of Lemma \ref{lemma:deterministic-marks} }\label{App-new}

\begin{proof}
	Fix a deterministic vector
	$\bm{k}=(k_R,k_B)$.
	For every $v\in \mathcal V_W$, the two random variables
	\[
	\widehat D^{(v)}_S(\bm{k})
	=
	\sum_{i=1}^{k_S+a_S} E_i^{S,(v)}
	\qquad
	\text{and}
	\qquad 
	\widehat D^{(v)}_{{\overline S}}(\bm{k})
	=
	\sum_{i=1}^{k_{\overline S}+a_{\overline S}} E_i^{\overline S,(v)}
	\]
	are independent binomial random variables with parameters
	$(k_S+a_S,p)$ and $(k_{\overline S}+a_{\overline S},p)$, respectively.
	Moreover, the corresponding pairs are independent for different
	vertices $v\in \mathcal V_W$. Therefore, the indicators
	\[
	\ind_
	{\left\{
	\widehat D^{(v)}_S(\bm{k})
	-
	\widehat D^{(v)}_{{\overline S}}(\bm{k})
	\ge r
	\right\}},
	\qquad v\in \mathcal V_W,
	\]
	are independent and identically distributed Bernoulli random
	variables with mean $\pi_S(\bm{k})$. This proves
	\eqref{static-binomial}.
	
	Next, by the definitions of $D_S^{(v)}[k]$ and
	$\bm N[k]$, we have, pathwise,
	\[
	D_S^{(v)}[k]
	=
	\sum_{i=1}^{N_S[k]+a_S} E_i^{S,(v)}
	=
	\widehat D_S^{(v)}(\bm N[k]).
	\]
	Hence
	\[
	\ind_{
	\left\{
	D_S^{(v)}[k]-D_{{\overline S}}^{(v)}[k]\ge r
	\right\}}
	=
	\ind_{
	\left\{
	\widehat D_S^{(v)}(\bm N[k])
	-
	\widehat D_{{\overline S}}^{(v)}(\bm N[k])
	\ge r
	\right\}},
	\]
	for every $v\in \mathcal V_W$, which proves
	\eqref{static-representation}.
	
	Finally, assume that
	$k_S\ge\ell_S$ and
	$k_{{\overline S}}\le\ell_{\overline S}$.
	Since all marks are non-negative,
	\[
	\widehat D_S^{(v)}(\bm{k})
	\ge
	\widehat D_S^{(v)}(\bm{\ell})
\qquad 
	\text{and}
\qquad
	\widehat D_{\overline S}^{(v)}(\bm{k})
	\le
	\widehat D_{\overline S}^{(v)}(\bm{\ell})
	\]
	for every $v\in \mathcal V_W$. Therefore,
	\[
	\widehat D_S^{(v)}(\bm{k})
	-
	\widehat D_{\overline S}^{(v)}(\bm{k})
	\ge
	\widehat D_S^{(v)}(\bm{\ell})
	-
	\widehat D_{\overline S}^{(v)}(\bm{\ell}),
	\]
	which immediately yields
	\eqref{static-monotonicity}.
\end{proof}

\section{ Further Consequences of Markovianity} \label{appendix-MC}

{
\begin{Proposition}\label{indeptau}
	{\it Define
		\[
		{\mathbb{S}}_m:=\{\bm z\in\mathbb{S}:\,\, R(\bm z)=m\},
		\quad
		\text{$m\in \{0,1,\ldots,n_W\}$.}
		\]
		For $k\in\mathbb{N}\cup\{0\}$ and
		$\{m_h\}_{0\leq h\leq k}\subseteq\{1,\ldots,n_W\}$,
		given the event
		\[
		\bigcap_{0\leq h\leq k}
		\{\bm Z_h\in\mathbb{S}_{m_h}\},
		\]
		provided that this event has positive probability, it holds that:
		
		\noindent $(i)$ The sojourn times
		$\{W_h\}_{0\leq h\leq k}$ of the Markov chain $\bm Z$
		are independent.
		
		\noindent $(ii)$ Each random variable $W_h$, $0\leq h\leq k$,
		is exponentially distributed with parameter $m_h$.}
\end{Proposition}

\begin{proof}
	By the Markov property of the process $\bm Z$, for any arbitrary
	finite sequence of states
	$\{\bm z_h\}_{0\leq h\leq k}\subset
	\mathbb{S}\setminus\mathbb{S}_0$
	and any arbitrary finite sequence of positive numbers
	$\{a_h\}_{0\leq h\leq k}\subset(0,\infty)$, the following
	well-known identity holds:
	\begin{align}
		\label{131123}
		&\mathbb{P}\left(
		\bigcap_{0\leq h\leq k}
		\left(
		\{\bm Z_h=\bm z_h\}
		\cap
		\{W_h>a_h\}
		\right)
		\right)
		\nonumber\\
		&\qquad =
		\mathbb{P}(\bm Z_0=\bm z_0)
		\left(
		\prod_{0\leq h<k}
		p_{\bm z_h\bm z_{h+1}}
		\right)
		\left(
		\prod_{0\leq h\leq k}
		{\rm e}^{-R(\bm z_h)a_h}
		\right),
	\end{align}
	where $(p_{\bm z \bm y} )$ denotes the transition matrix
	of the embedded Markov chain $\{\bm Z_k\}_k$.
	
	Moreover,
	\[
	\bigcap_{0\leq h\leq k}
	\{\bm Z_h\in\mathbb{S}_{m_h}\}
	=
	\bigcup_{
		\substack{
			\bm z_0\in\mathbb{S}_{m_0},\ldots,
			\bm z_k\in\mathbb{S}_{m_k}
	}}
	\bigcap_{0\leq h\leq k}
	\{\bm Z_h=\bm z_h\},
	\]
	where the union on the right-hand side is over disjoint events.
	Hence,
	\begin{align*}
		&\mathbb{P}\left(
		\bigcap_{0\leq h\leq k}\{W_h>a_h\}
		\,\Big|\,
		\bigcap_{0\leq h\leq k}
		\{\bm Z_h\in\mathbb{S}_{m_h}\}
		\right)
		\\
		&=
		\frac{
			\displaystyle
			\sum_{
				\substack{
					\bm z_0\in\mathbb{S}_{m_0},\ldots,
					\bm z_k\in\mathbb{S}_{m_k}
			}}
			\mathbb{P}(\bm Z_0=\bm z_0)
			\left(
			\prod_{0\leq h<k}		p_{\bm z_h\bm z_{h+1}}
			\right)
			\left(
			\prod_{0\leq h\leq k}
			{\rm e}^{-R(\bm z_h)a_h}
			\right)
		}{
			\displaystyle	\sum_{	\substack{
					\bm z_0\in\mathbb{S}_{m_0},\ldots,
					\bm z_k\in\mathbb{S}_{m_k}
			}}	\mathbb{P}(\bm Z_0=\bm z_0) \prod_{0\leq h<k}
			p_{\bm z_h\bm z_{h+1}}
		}.
	\end{align*}
	Since $\bm z_h\in\mathbb{S}_{m_h}$ implies
	$R(\bm z_h)=m_h$, we obtain
	\begin{align*}
		&\mathbb{P}\left(
		\bigcap_{0\leq h\leq k}\{W_h>a_h\}		\,\Big|\,
		\bigcap_{0\leq h\leq k}	\{\bm Z_h\in\mathbb{S}_{m_h}\}		\right)
		\\
		&\qquad =
		\prod_{0\leq h\leq k}
		{\rm e}^{-m_h a_h}.
	\end{align*}
	The right-hand side is the product of the survival functions of
	independent exponential random variables with parameters
	$m_0,\ldots,m_k$. This proves both claims.
\end{proof}
}

\section{Concentration inequalities}\label{Penrose}
Throughout the paper,  we extensively employ classical deviation bounds for binomial and Poisson distributions.  These results can be found e.g.  in \cite{P},  and are reported here for the reader's convenience.  Hereafter,  {$\FH$} denotes the function defined in \eqref{eq:H}.  

Let $X\sim\mathrm{Bin}(m,q)$ and let $\mu:=mq$, with
$m\in\mathbb N$ and $q\in(0,1)$. For any real $x\in(0,m)$,
the following inequalities hold:

\begin{equation}\label{Penrose-coda-sopra}
	\mathbb P(X\geq x)\leq
	\begin{cases}
		\exp\left\{-\mu\FH\left(\frac{x}{\mu}\right)\right\},
		& x\geq\mu,\\[1mm]
		\exp\left\{-\frac{x}{2}\log\left(\frac{x}{\mu}\right)\right\},
		& x\geq \mathrm e^2\mu.
	\end{cases}
\end{equation}

and

\begin{equation}\label{Penrose-coda-sotto}
	\mathbb P(X\leq x)
	\leq
	\exp\left\{-\mu\FH\left(\frac{x}{\mu}\right)\right\},
	\qquad x\leq\mu.
\end{equation}
Let $\lambda>0$. For every real $x\in[0,\lambda]$,
\begin{equation}\label{Penrose-coda-sotto-Po}
	\mathbb P\big(\mathrm{Po}(\lambda)\leq x\big)
	\leq
	\exp\left\{
	-\lambda\FH\left(\frac{x}{\lambda}\right)
	\right\}.
\end{equation}
Although $X$ is integer-valued, the bounds are stated for real
thresholds; this follows immediately by replacing $x$ with
$\lceil x\rceil$ or $\lfloor x\rfloor$ and using the monotonicity
of $\FH$ on either side of $1$.
\section{Properties of the solutions of Cauchy's problem \ref{eq:CP},  and proof of Proposition \ref{le:coupledphi}}\label{Appendix-CP}

We begin  by stating a lemma which establishes a relationship between $\bm f$ and $\bm g$,  i.e.,  the maximal solutions of Cauchy's problems \eqref{eq:CP} and \eqref{simplified:CPcoupled},  respectively.  This relationship holds when $q\ll p^{-1}$,  which entails $\beta_S(x_R,x_B)=\beta_S(x_S)$.
The proof is omitted since the claim follows directly by inspection.

\begin{Lemma}\label{le:transfer}
	Assume $\beta_S(x_R,x_B)=\beta_S(x_S)$,  $S\in\{R,B\}$,  and that the Cauchy problem \eqref{simplified:CPcoupled} has a unique maximal solution $\bm g$ on $(0,\kappa_{\bm g})$ {with $g_R$ and $g_B$
		strictly increasing.} Then the Cauchy problem
	\eqref{eq:CP} has a unique maximal solution $\bm f$ on $(0,\kappa_{\bm f})$, with $\kappa_{\bm f}:=\lim_{y \uparrow \kappa_{\bm g}} z(y)$ and $z:=g_R+g_B$,  provided by
	\begin{equation*}
		\bm f(x)=\bm g(z^{-1}(x)).
	\end{equation*}
\end{Lemma}

Under the assumption
$\beta_S(x_R,x_B)=\beta_S(x_S)$,  $\bm g$ can be written 
in terms of the maximal solutions of the following one-dimensional Cauchy problems:
\begin{equation}\label{simplified:CP}
h_S'(y)=\beta_S(h_S(y)),
\quad y\in (0,\kappa_{h_S}),\quad h_S(0)=0,
\end{equation}
i.e.,
\begin{equation} \label{eq:identityCP}
	\bm g(y)\equiv (h_R(y), h_B(y)),\quad  y\in (0, \kappa_{\bm g}),\quad
	\kappa_{\bm g}:=\min\{\kappa_{h_R}, \kappa_{h_B}\}
\end{equation}
As a consequence,  for $g=q$ or $g \ll q \ll p^{-1}$, 
we can compute  $ \lim_{x\uparrow \kappa_{\bm f}}\bm f(x)$  first evaluating   
$\lim_{y\uparrow \kappa_{\bm g}}(h_R(y),h_B(y) )$ and then  invoking both Lemma \ref{le:transfer} and  identity \eqref{eq:identityCP}. 
{\color{blue}
\begin{Remark}
	Assume $q=g$ or $g\ll q\ll p^{-1}$.
	The functions $\bm f$ and $\bm g$ describe the same
	limiting evolution under two different parameterizations.
	
	By Theorem~\ref{prop:equadiff},
	$\bm f(x)$ describes the normalized trajectory
	$\widetilde{\bm N}(xq)/q$, where $xq$ parameterizes
	the total number of non-seed activations.
	In contrast, $\bm g(y)$ describes the evolution in
	rescaled physical time
	\[
	y=\eta t,
	\]
	where $\eta$ is defined in \eqref{eta}.
	More precisely, writing
	$\bm N^\circ=(N_R^\circ,N_B^\circ)$ for the original
	physical process, for every $y_0\in(0,\kappa_{\bm g})$
	we have
	\[
	\sup_{0\le y\le y_0}
	\left\|
	\frac{\bm N^\circ(y/\eta)}{q}-\bm g(y)
	\right\|
	\longrightarrow0,
	\qquad\text{a.s.}
	\]
	This follows from
	Theorems~\ref{prop:equadiff}, \ref{prop:kminT},
	and~\ref{propVk}, by inverting the limiting
	activation-time map and using
	Lemma~\ref{le:transfer}.
	
	In particular, when $q=g$, we have $\eta=1$,
	so the argument of $\bm g$ is unscaled physical time.
	When $g\ll q\ll p^{-1}$, physical time must instead
	be accelerated by the factor $\eta$.
	
	Since $\beta_S(x_R,x_B)=\beta_S(x_S)$ in these regimes,
	the scalar equations in \eqref{simplified:CP}
	show that the two activation processes are decoupled
	at first order in rescaled physical time, even though
	the finite systems remain coupled.
\end{Remark}

}
The interpretation of the identities \eqref{simplified:CP} and \eqref{eq:identityCP} is that the two activation processes evolve largely independently  over \lq\lq physical \rq\rq\,time,  exhibiting a negligible dependence.

The next two lemmas provide some properties of $h_S$ when $g=q$ and $g \ll q \ll p^{-1}$,  respectively.

\begin{Lemma}\label{le:subcrit-supercrit}
	Assume $q=g$. \\
	(i) If $\alpha_S<1$,  then the Cauchy problem \eqref{simplified:CP} has a unique (strictly increasing) solution $h_S$ 
	on $(0,\infty)$ and $h_S(y)\uparrow z_S$, as $y\uparrow+\infty$.\\
	(ii) If $\alpha_S>1$,  then the Cauchy problem \eqref{simplified:CP} has a unique (strictly increasing) solution $h_S$ on $(0,\kappa_{h_S})$,
	where
	\begin{equation*}
		\kappa_{h_S}:=\int_{0}^{\infty}\frac{\dd y}{\beta_S(y)}<\infty
	\end{equation*}
	and $h_S(y)\uparrow+\infty$, as $y\uparrow\kappa_{h_S}$.  Moreover $\kappa_{h_R}<\kappa_{h_B}$.\\
	(iii) If $\alpha_S=1$, then the Cauchy problem
	\eqref{simplified:CP} has a unique strictly increasing solution
	$h_S$ on $(0,\infty)$ and
	\[
	h_S(y)\uparrow z_S,
	\qquad y\to\infty.
	\]
\end{Lemma}
\begin{proof}
	\noindent{\it Proof of $(i)$.} By Remark \ref{re:29Ott} the function $\beta_S(x_S)$ has two strictly positive zeros, say $z_S<z'_{S}$,
	which represent two equilibrium points for the dynamical system.  Furthermore $\beta_S(x_S)$ is positive  and decreasing for $x_S<z_S$.  
	Since  $h_S(0)=0<z_S$,  necessarily $h_S(y)\le z_S$ for every $y\in [0,\infty)$, 
	then $\kappa_{h_S}=+\infty$,  $h'_S(y)=\beta_S( h_S(y))\ge 0 $ for any $y\in [0,\infty)$, 
	and $\lim_{y \to +\infty}h_S(y)=\sup_{y\in [0,\infty)} h_S(y)$.  Since
	\[
	z_S\ge  \lim_{y\to \infty} h_S(y)=\lim_{y\to \infty} \int _0^y h'_S(u) \dd  u= \lim_{y\to \infty} \int _0^y 
	\beta_S(h_S(u)) \dd  u\ge \lim_{y\to \infty} y \beta_S(h_S(y)),
	\]
	we finally have $\lim_{y\to \infty }\beta_S(h_S(y))=\beta_S( \lim_{y\to \infty }h_S(y) ) =0$.\\
	\noindent{\it Proof of $(ii)$.} By Remark \ref{re:29Ott} the function $\beta_S(x_S)$ is strictly positive for $x_S\ge 0$. Moreover,
	$\lim_{x_S\to+\infty}\beta_S(x_S)=+\infty$, as can be checked directly. Therefore $\inf_{x_S\in [0,\infty)} \beta_S(x_S)>0$.
	So the unique solution $h_S$ is strictly increasing, and $h'_S(y)$ is bounded away from zero for all $y$. In particular, this property guarantees that $h_S$ has no horizontal asymptotes.
	Therefore there are only two possible cases: $(i)$ $h_S$ is defined on the whole non-negative half-line $[0,\infty)$ and $h_S(y)\uparrow+\infty$,
	as $y\uparrow+\infty$; $(ii)$ $h_S$ is defined on a finite interval of the form $[0,\kappa_{h_S})$, for some $\kappa_{h_S}\in (0,\infty)$
	and $h_S(y)\uparrow+\infty$, as $y\uparrow\kappa_{h_S}$. We now verify that case $(ii)$ holds.
	Let $\mathcal{D}_{h_S}$ be the domain of $h_S$.  From \eqref{simplified:CP}, we have
	\begin{equation}\label{160124}
		\frac{h'_S(y)}{\beta_S(h_S(y))}= 1,\quad\text{$\forall$ $y\in \mathcal{D}_{h_S}$}.
	\end{equation}
	 Integrating both sides yields
	\begin{equation}\label{sep-var}
		\int_{h_S(0)}^{h_S(y)}\frac{1}{\beta_S(u)}\,\dd u=\int_{0}^{y}\frac{h'_S(u)}{\beta_S(h_S(u))}\,\dd u=\int_0^y\dd u=y,
		\quad\text{$\forall$ $y\in \mathcal{D}_{h_S}$}.
	\end{equation}
	Now observe that
	\[
	\int_{0}^{\infty}\frac{1}{\beta_S(u)}\,\dd  u=\int_{0}^{\infty}\frac{\dd u}{-u + r^{-1}[(1- r^{-1})]^{r-1}(\alpha_S+ u)^r}=\kappa_{h_S}<\infty.
	\]
	Therefore by \eqref{sep-var} we conclude that $\mathcal{D}_{h_S}=[0,\kappa_{h_S})$ and $h_S(y)\uparrow+\infty$,
	as $y\uparrow\kappa_{h_S}$. 
	Finally, if $\alpha_B>1$, since
	$\alpha_R>\alpha_B$ implies
	$\beta_R(x)>\beta_B(x)$ for every $x\geq0$, we have
	\[
	\kappa_{h_R}
	=
	\int_0^\infty\frac{dx}{\beta_R(x)}
	<
	\int_0^\infty\frac{dx}{\beta_B(x)}
	=
	\kappa_{h_B}.
	\]
	If instead $\alpha_B\leq1$, parts $(i)$ and $(iii)$ give
	$\kappa_{h_B}=\infty$, whereas $\kappa_{h_R}<\infty$.
	Therefore, in all cases,
	\[
	\kappa_{\bm g}=\kappa_{h_R}.
	\]

\noindent{\it Proof of $(iii)$.}
The proof is analogous to that of part $(i)$. Since $z_S$ is the
unique positive zero of $\beta_S$, $\beta_S(x)>0$ for
$x\in[0,z_S)$, and $h_S(0)=0$, uniqueness prevents the solution
from crossing the equilibrium point $z_S$. Hence $h_S$ is
increasing and bounded by $z_S$, and therefore it is defined on
$[0,\infty)$. Its limit must be a zero of $\beta_S$, and therefore
\[
\lim_{y\to\infty}h_S(y)=z_S.
\]

\end{proof}

When $g \ll q\ll p^{-1}$ we have the analytic expression of $h_S$.  Indeed,  the next lemma holds.

\begin{Lemma}\label{le:betapos}
	Let $g\ll q\ll p^{-1}$.  Then the Cauchy problem \eqref{simplified:CP} has a unique solution $h_S$ on $(0,\kappa_{h_S})$,  with
	\begin{equation*}
		h_S(y):=\frac{1}{\left(\alpha_S^{1-r}-\frac{r-1}{r!}y\right)^{1/(r-1)}}-\alpha_S,\qquad	\kappa_{h_S}:=\frac{r!}{(r-1)\alpha_S^{r-1}}.
	\end{equation*}
\end{Lemma}
The claim follows by direct inspection, so the proof is omitted.

Now,  we direct our attention to the case $q=p^{-1}$.  In this scenario,  the identity  $\beta_S(x_R,x_B) =\beta_S(x_S)$  no longer holds and so the previous methodology 
is no longer applicable.   Nevertheless, a comparative analysis is possible by examining the solution of the Cauchy problem \eqref{eq:CP} in relation to the solution of an auxiliary Cauchy problem where the aforementioned identity holds true.

\begin{Lemma}\label{le-CPq=p-1}
	Assume $q=p^{-1}$, 
	let $\bm f$ be the solution of the Cauchy problem \eqref{eq:CP}  and 
	let $\widetilde{\bm f}$ be the solution of the Cauchy problem
	\begin{equation}\label{eq:CP2}
		\widetilde{\bm f}^{'}(x)=\frac{\widetilde{\bm{\beta}} (\widetilde{\bm f}(x))}{\widetilde \beta_R(\widetilde{f}_{R}(x))+
			\widetilde{\beta}_{B}(\widetilde f_{B}(x)) },
		\quad x\in (0,\kappa_{\bm f}),\quad\text{$\widetilde{\bm f}(0)=(0,0)$}
	\end{equation}
	where $\widetilde \beta_S(\bm x)=\frac{1}{r!}(x_S+\alpha_S)^r$,  $S\in\{R,B\}$.
	Then
	$f_R(x)> \widetilde f_R(x)$ and $f_B(x)< \widetilde f_B(x)$,  for every $x\in(0,\kappa_{\bm f})$. 
\end{Lemma}
\begin{proof}
	First,  note that
	\begin{equation}
		\frac{\beta_B( x_R, x_B )} {\beta_R( x_R,  x_B)  }<\frac{\widetilde{\beta}_B( x_B )} {\widetilde{\beta}_R( x_R)  }=
		\left( \frac{ x_B+\alpha_B}{ x_R+\alpha_R}\right)^r,  \qquad \text{for }  x_B+\alpha_B< x_R+\alpha_R.  \label{15-5-25}
	\end{equation}
	Indeed, grouping the terms in the series defining
	$\beta_R$ and $\beta_B$ according to
	$d=r'-r''\geq r$, the corresponding terms in $\beta_B$
	are obtained from those in $\beta_R$ by multiplication by
	\[
	\left(
	\frac{x_B+\alpha_B}{x_R+\alpha_R}
	\right)^d.
	\]
	Since $x_B+\alpha_B<x_R+\alpha_R$ and terms with $d>r$
	are strictly positive, \eqref{15-5-25} follows.
	
	Second,  note that $\widetilde{\beta}_R(\cdot)$ and  $\widetilde{\beta}_B(\cdot)$ have the same expression of 
	$\beta_R(\cdot)$ and  $\beta_B(\cdot)$ for the case $g\ll q\ll p^{-1}$, therefore 
	we can apply  Lemma \ref{le:betapos},  identity \eqref{eq:identityCP}  and  Lemma \ref{le:transfer} to obtain the analytical expression of $\widetilde {\bm f}(x)$,  from which we infer that $\widetilde{f}_R(x)\ge \widetilde{f}_B(x)$,  for any $x\in [0, \infty)$.
	By \eqref{15-5-25} we have
	\[
	f'_R(0)=\frac{\beta_R (0,0)} {  \beta_R(0,0)+
		{\beta}_{B}( 0,0)}> \frac{\widetilde{\beta}_R (0)}{\widetilde \beta_R(0)+
		\widetilde{\beta}_{B}(0)}=\widetilde f'_R(0)
	\]
	and similarly $f'_B(0) < \widetilde f'_B(0)$. 
	Therefore $f_R(x) > \widetilde f_R(x)$ and $f_B(x) < \widetilde f_B(x)$
	in a right-neighborhood of $0$.
	Reasoning by contradiction,  assume $x_0<\kappa_{\bm f}$,  where
	\[
	x_0:=\inf\{x>0:  f_R(x) \le \widetilde f_R(x) \text{ or } f_B(x) \ge  \widetilde f_B(x)\}.
	\]
	Then
	\begin{align*}
		f_R(x_0) &=f_R(0)+ \int_0^{x_0} \frac{{{\beta}_R} (\bm f(x))}{ \beta_R(\bm f(x))+\beta_B(\bm f(x))}\dd  x\\
		&> f_R(0)+\int_0^{x_0} \frac{{\widetilde {\beta}_R} ({f}_R(x))}{ \widetilde \beta_R(f_R(x))+
			\widetilde \beta_B({f}_B(x))}\dd  x\\
		& > \widetilde f_R(0)+\int_0^{x_0} 
		\frac{ {\widetilde {\beta}_R} (\widetilde{f}_R(x))} { \widetilde \beta_R(\widetilde{f}_R(x))+
			\widetilde{\beta}_B(\widetilde{f}_B(x))}\dd  x=\widetilde{f}_R(x_0)
	\end{align*}
	and similarly $f_B(x_0)<\widetilde f_B(x_0)$.  This contradicts the definition of $x_0$,
	and thus concludes the proof of the lemma.
\end{proof}

\subsection{Proof of Proposition \ref{le:coupledphi}}
Cases $(i)$ and $(ii)$ of Proposition \ref{le:coupledphi} follow directly by Lemmas  \ref{le:transfer},  \ref{le:subcrit-supercrit},   and the identity \eqref{eq:identityCP}.
Case $(iii)$   follows from Lemmas \ref{le:transfer} and 
\ref{le:betapos} and the identity \eqref{eq:identityCP}.
Case $(iv)$ 
easily follows from Lemma \ref{le-CPq=p-1},  since,  as already mentioned,
$\widetilde {\bm f}$
coincides with the solution of the Cauchy problem \eqref{eq:CP} for  $g\ll q\ll p^{-1}$.
Finally case $(v)$ is of immediate verification.

{
\section{Proof of Proposition \ref{cor:BC}}\label{sec:conc}
{
	
	\subsection{Proof strategy and reduction}
	
	In this appendix we prove Proposition~\ref{cor:BC}.
	The proof is organized in a top-down fashion.
	We first reduce the claim to two uniform concentration estimates and
	then establish the auxiliary results needed to prove them.
	
	Throughout this appendix, we set
	\begin{equation}
		\xi_n:=\eta q,
		\label{D:xi}
	\end{equation}
	where $\eta$ is defined in~\eqref{eta}. Hence
	\begin{equation}
		\xi_n=
		\begin{cases}
			q, & q=g,\\[1mm]
			n(qp)^r, & g\ll q\ll p^{-1},\\[1mm]
			n, & q=p^{-1}\ \text{or}\ p^{-1}\ll q\ll n.
		\end{cases}
		\label{D:xi-regimes}
	\end{equation}
	
	For $S\in\{R,B\}$, define
	\begin{equation}
		\Upsilon_S(\kappa)
		:=
		\sup_{\bm{k}\in{\mathbb T}(\kappa)}
		Y_S(\bm{k}),
		\qquad
		\widehat\Upsilon_S(\kappa)
		:=
		\sup_{\bm{k}\in{\mathbb T}(\kappa)}
		\widehat Y_S(\bm{k}),
		\label{D:upsilon}
	\end{equation}
	where $Y_S(\bm{k})$ and $\widehat Y_S(\bm{k})$
	are defined in~\eqref{defY}--\eqref{defhatY}, and
	\begin{equation}
		\Psi_S(\kappa)
		:=
		\sup_{j\le \kappa q}
		\frac{|\widehat N_S[j]|}{q}.
		\label{D:Psi}
	\end{equation}
	
	Recalling the definition of $\Gamma_S(\kappa)$ in
	\eqref{eq:zero6nov}, we have
	\begin{equation}
		\Gamma_S(\kappa)
		=
		\max\left\{
		\Upsilon_S(\kappa),
		\frac{\widehat\Upsilon_S(\kappa)}{\xi_n},
		\Psi_S(\kappa)
		\right\}.
		\label{D:Gamma-decomp}
	\end{equation}
	Thus, Proposition~\ref{cor:BC} follows once uniform exponential
	bounds are established for the three quantities appearing on the
	right-hand side of \eqref{D:Gamma-decomp}.
	
	The first two quantities are controlled by the following proposition.
	
	\begin{Proposition}
		\label{prop:D1-new}
		Let $\kappa\in(0,\kappa_{\bm f})$, $\delta>0$, and
		$S\in\{R,B\}$. There exist positive constants
		$c=c(\kappa,\delta)$ and $C=C(\kappa,\delta)$ such that,
		for all sufficiently large $n$,
		\begin{equation}
			\max
			\left\{
			\mathbb P
			\left(
			\Upsilon_S(\kappa)>\delta
			\right),
			\,
			\mathbb P
			\left(
			\widehat\Upsilon_S(\kappa)>\delta\xi_n
			\right)
			\right\}
			\le
			Ce^{-c\xi_n}.
			\label{D:prop1}
		\end{equation}
	\end{Proposition}
	
	The martingale term is controlled by the following proposition.
	
	\begin{Proposition}
		\label{prop:D2-new}
		Let $S\in\{R,B\}$.
		For every $\kappa>0$ and $\delta>0$,
		\begin{equation}
			\mathbb P
			\left(
			\Psi_S(\kappa)>\delta
			\right)
			\le
			2\kappa q
			\exp
			\left\{
			-\frac{\delta^2q}{2\kappa}
			\right\}
			\label{D:Psi-bound}
		\end{equation}
		for all sufficiently large $n$.
		In particular, there exists
		$c=c(\kappa,\delta)>0$ such that
		\begin{equation}
			\mathbb P
			\left(
			\Psi_S(\kappa)>\delta
			\right)
			\le
			e^{-cq}
			\label{D:Psi-exp}
		\end{equation}
		up to a multiplicative subexponential factor.
	\end{Proposition}
	
	We first show that these two propositions are sufficient to conclude
	the proof of Proposition~\ref{cor:BC}.
	
	\begin{proof}[Reduction of Proposition~\ref{cor:BC}]
		Fix $\kappa\in(0,\kappa_{\bm f})$ and $S\in\{R,B\}$.
		For every $\delta>0$, Propositions~\ref{prop:D1-new}
		and~\ref{prop:D2-new} yield positive constants
		$c_\delta,C_\delta$ such that
		\begin{align}
			\mathbb P
			\left(
			\Gamma_S(\kappa)>\delta
			\right)
			&\le
			C_\delta e^{-c_\delta\xi_n}
			+
			2\kappa q
			\exp
			\left\{
			-c_\delta q
			\right\}
			\nonumber\\
			&\le
			C'_\delta q\,e^{-c'_\delta q},
			\label{D:Gamma-tail}
		\end{align}
		for suitable positive constants
		$c'_\delta,C'_\delta$, where we have used
		$\xi_n\ge q$ for all sufficiently large $n$.
		
		Under assumption~\eqref{eq:bootcond}, the critical scale $g$
		satisfies
		\begin{equation}
			\frac{g}{\log n}\longrightarrow\infty.
			\label{D:g-logn}
		\end{equation}
		Since $q\ge g$, it follows that
		\[
		\sum_{n=1}^{\infty}
		q\,e^{-cq}
		<
		\infty
		\qquad
		\text{for every }c>0.
		\]
		Hence~\eqref{D:Gamma-tail} implies
		\[
		\sum_{n=1}^{\infty}
		\mathbb P
		\left(
		\Gamma_S(\kappa)>\delta
		\right)
		<
		\infty.
		\]
		By the Borel--Cantelli lemma,
		\[
		\mathbb P
		\left(
		\Gamma_S(\kappa)>\delta
		\ \text{infinitely often}
		\right)
		=
		0.
		\]
		Since this holds for every rational $\delta>0$, we obtain
		\[
		\Gamma_S(\kappa)
		\longrightarrow0,
		\qquad\text{a.s.},
		\]
		for each $S\in\{R,B\}$.
		Therefore, once Propositions~\ref{prop:D1-new}
		and~\ref{prop:D2-new} are established,
		Proposition~\ref{cor:BC} follows.
	\end{proof}

	\subsection{Proof of Proposition \ref{prop:D1-new}}
	
	The proof of Proposition~\ref{prop:D1-new} is based on a uniform
	approximation of $Q^S_{k+1}$.
	
	\begin{Lemma}
		\label{lemma:D-Q}
		Let $\kappa\in(0,\kappa_{\bm f})$, $\delta>0$, and
		$S\in\{R,B\}$. There exist positive constants
		$c=c(\kappa,\delta)$ and $C=C(\kappa,\delta)$ such that,
		for all sufficiently large $n$,
		\begin{equation}
			\mathbb P
			\left(
			\sup_{\bm{k}\in{\mathbb T}(\kappa)}
			\ind_{\{\bm N[k]=\bm{k}\}}
			\left|
			Q^S_{k+1}
			-
			\xi_n\beta_S(\bm{k}/q)
			\right|
			>
			\delta\xi_n
			\right)
			\le
			Ce^{-c\xi_n},
			\label{D:Q-concentration}
		\end{equation}
		where $k=k_R+k_B$.
	\end{Lemma}
	
	Assuming Lemma~\ref{lemma:D-Q} for the moment, we can immediately
	prove Proposition~\ref{prop:D1-new}.
	
	\begin{proof}[Proof of Proposition~\ref{prop:D1-new}]
		The estimate concerning
		$\widehat\Upsilon_S(\kappa)$ follows directly from
		Lemma~\ref{lemma:D-Q}.
		
		It remains to control $\Upsilon_S(\kappa)$.
		Set
		\begin{equation}
			b_\kappa
			:=
			\inf_{\bm{x}\in{\mathbb T}'(\kappa)}
			\left(
			|\beta_R(\bm{x})|
			+
			|\beta_B(\bm{x})|
			\right).
			\label{D:b-kappa}
		\end{equation}
		Since $\kappa<\kappa_{\bm f}$, the definition of
		$\kappa_{\bm f}$ and the properties of
		$\beta_R$ and $\beta_B$ imply
		\begin{equation}
			b_\kappa>0.
			\label{D:b-positive}
		\end{equation}
		
		Recall that, before the stopping time, $Q^S_{k+1}$ equals the
		number of jointly white and $S$-\supra{} nodes.
		After the stopping time, $U^S_{k+1}$ is defined by
		\eqref{eq:redefineU}. Thus, along the extended process,
		with the convention $0/0:=1/2$,
		\begin{equation}
			U^S_{k+1}
			=
			\frac{|Q^S_{k+1}|}
			{|Q^R_{k+1}|+|Q^B_{k+1}|}.
			\label{D:UQ}
		\end{equation}
		
		Fix
		\[
		0<\varepsilon<\frac{b_\kappa}{4}.
		\]
		Suppose that, for $T\in\{R,B\}$,
		\begin{equation}
			\left|
			Q^T_{k+1}
			-
			\xi_n\beta_T(\bm{k}/q)
			\right|
			\le
			\varepsilon\xi_n.
			\label{D:good-Q}
		\end{equation}
		Then
		\[
		\left|
		|Q^T_{k+1}|
		-
		\xi_n|\beta_T(\bm{k}/q)|
		\right|
		\le
		\varepsilon\xi_n.
		\]
		Consequently,
		\begin{align}
			|Q^R_{k+1}|+|Q^B_{k+1}|
			&\ge
			\xi_n
			\left(
			|\beta_R(\bm{k}/q)|
			+
			|\beta_B(\bm{k}/q)|
			-
			2\varepsilon
			\right)
			\nonumber\\
			&\ge
			\frac{b_\kappa}{2}\xi_n.
			\label{D:den-lower}
		\end{align}
		
		A direct computation then gives
		\begin{align}
			&
			\left|
			\frac{|Q^S_{k+1}|}
			{|Q^R_{k+1}|+|Q^B_{k+1}|}
			-
			\frac{|\beta_S(\bm{k}/q)|}
			{|\beta_R(\bm{k}/q)|
				+
				|\beta_B(\bm{k}/q)|}
			\right|
			\nonumber\\
			&\qquad
			\le
			\frac{6\varepsilon}{b_\kappa}.
			\label{D:ratio-bound}
		\end{align}
		
		Choose $\varepsilon$ sufficiently small so that
$	\frac{6\varepsilon}{b_\kappa}<\delta.$
		It follows from~\eqref{D:ratio-bound} that
		\begin{align}
			\{\Upsilon_S(\kappa)>\delta\}
			\subseteq
			\bigcup_{T\in\{R,B\}}
			\Bigg\{
			\sup_{\bm{k}\in\mathbb T(\kappa)}
			\ind_{\{\bm N[k]=\bm{k}\}}
			\left|
			Q^T_{k+1}
			-
			\xi_n\beta_T(\bm{k}/q)
			\right|
			>
			\varepsilon\xi_n
			\Bigg\}.
			\label{D:U-inclusion}
		\end{align}
		Lemma~\ref{lemma:D-Q} therefore yields the desired
		exponential estimate for $\Upsilon_S(\kappa)$ and completes
		the proof.
	\end{proof}

	\subsection{Proof of Lemma \ref{lemma:D-Q}}
	
	We now prove the uniform approximation of $Q^S_{k+1}$.
	The key idea is to separate the error into three contributions:
	a fluctuation term, a deterministic mean-approximation term, and
	a residual term.
	
	To this end, for a deterministic vector
	$\bm{k}=(k_R,k_B)\in(\mathbb N\cup\{0\})^2$, define
	\begin{equation}
		X_S(\bm{k})
		:=
		\left|
		\widehat{\mathcal S}_S(\bm{k})
		\right|,
		\qquad S\in\{R,B\},
		\label{D:X}
	\end{equation}
	where $\widehat{\mathcal S}_S(\bm{k})$ is the
	fixed-configuration \Ssupra{} set introduced in
	Lemma~\ref{lemma:deterministic-marks}.
	
	For every deterministic $\bm{k}$,
	\begin{equation}
		X_S(\bm{k})
		\overset{\mathrm L}{=}
		\text{Bin}
		\left(
		n_W,\pi_S(\bm{k})
		\right),
		\label{D:X-bin}
	\end{equation}
	where
	\begin{equation}
		\pi_S(\bm{k})
		:=
		\mathbb P
		\left(
		\text{Bin}(k_S+a_S,p)
		-
		\text{Bin}(k_{\overline S}+a_{\overline S},p)
		\ge r
		\right),
		\label{D:pi}
	\end{equation}
	and the two binomial random variables appearing on the
	right-hand side are independent.
	
	Moreover, pathwise,
	\begin{equation}
		\mathcal S_S[k]
		=
		\widehat{\mathcal S}_S(\bm N[k]),
		\qquad k\le n_W.
		\label{D:static-path}
	\end{equation}
	
	For later use, define also
	\begin{equation}
		C_S(\bm{k})
		:=
		\sum_{v\in\mathcal V_W}
		\ind_{
		\left\{
		\widehat D^{(v)}_S(\bm{k})\ge r,\,
		\widehat D^{(v)}_{\overline S}(\bm{k})\ge1
		\right\}}.
		\label{D:C}
	\end{equation}
	For every deterministic $\bm{k}$,
	\begin{equation}
		C_S(\bm{k})
		\overset{\mathrm L}{=}
		\text{Bin}
		\left(
		n_W,\widetilde\pi_S(\bm{k})
		\right),
		\label{D:C-bin}
	\end{equation}
	where
	\begin{equation}
		\widetilde\pi_S(\bm{k})
		:=
		\mathbb P
		\left(
		\text{Bin}(k_S+a_S,p)\ge r
		\right)
		\mathbb P
		\left(
		\text{Bin}(k_{\overline S}+a_{\overline S},p)\ge1
		\right).
		\label{D:pi-tilde}
	\end{equation}
	
	We shall use the following two fixed-configuration estimates, whose
	proofs are postponed to the next subsection.
	
	\begin{Lemma}
		\label{lemma:D-mean}
		Let $\kappa\in(0,\kappa_{\bm f})$ and $S\in\{R,B\}$.
		Then
		\begin{equation}
			\sup_{\bm{k}\in{\mathbb T}(\kappa)}
			\left|
			\frac{n_W\pi_S(\bm{k})-k_S}{\xi_n}
			-
			\beta_S(\bm{k}/q)
			\right|
			\longrightarrow0.
			\label{D:mean}
		\end{equation}
		Moreover, if $q=g$,
		\begin{equation}
			\sup_{\bm{k}\in{\mathbb T}(\kappa)}
			\frac{
				n_W\widetilde\pi_S(\bm{k})
			}{q}
			\longrightarrow0.
			\label{D:mixed-mean}
		\end{equation}
	\end{Lemma}
	
	\begin{Lemma}
		\label{lemma:D-concentration}
		Let $\kappa\in(0,\kappa_{\bm f})$, $\delta>0$, and
		$S\in\{R,B\}$. There exist positive constants
		$c=c(\kappa,\delta)$ and $C=C(\kappa,\delta)$ such that,
		for all sufficiently large $n$,
		\begin{equation}
			\mathbb P
			\left(
			\sup_{\bm{k}\in{\mathbb T}(\kappa)}
			\left|
			X_S(\bm{k})
			-
			n_W\pi_S(\bm{k})
			\right|
			>
			\delta\xi_n
			\right)
			\le
			C e^{-c\xi_n}.
			\label{D:X-concentration}
		\end{equation}
		
		If $q=g$, then there exist positive constants
		$c'=c'(\kappa,\delta)$ and
		$C'=C'(\kappa,\delta)$ such that
		\begin{equation}
			\mathbb P
			\left(
			\sup_{\bm{k}\in{\mathbb T}(\kappa)}
			C_S(\bm{k})
			>
			\delta q
			\right)
			\le
			C'e^{-c'q}.
			\label{D:C-concentration}
		\end{equation}
	\end{Lemma}
	
	We now show how these two lemmas imply
	Lemma~\ref{lemma:D-Q}.
	
	\begin{proof}[Proof of Lemma~\ref{lemma:D-Q}]
		Fix
		$\bm{k}=(k_R,k_B)\in{\mathbb T}(\kappa)$
		and let $k=k_R+k_B$.
		
		On the event
		$\{\bm N[k]=\bm{k}\}$,
		the fixed-configuration representation gives
		\begin{equation}
			|\mathcal S_S[k]|=X_S(\bm{k}),
			\qquad
			N_S[k]=k_S,
			\label{D:on-N}
		\end{equation}
		and, for every $v\in\mathcal V_W$,
		\[
		D^{(v)}_R[k]=\widehat D^{(v)}_R(\bm{k}),
		\qquad
		D^{(v)}_B[k]=\widehat D^{(v)}_B(\bm{k}).
		\]
		
		Using~\eqref{eq:QSnew}, on this event we can write
		\begin{equation}
			Q^S_{k+1}
			=
			X_S(\bm{k})
			-
			k_S
			+
			R_S[k],
			\label{D:Q-decomposition}
		\end{equation}
		where
		\begin{align}
			R_S[k]
			:={}&
			\left|
			(\mathcal V_W\setminus\mathcal S_S[k])
			\cap \mathcal V_S[k]
			\cap
			\{v:D^{(v)}_S[k]\ge r\}
			\right|
			\nonumber\\
			&-
			\left|
			\mathcal S_S[k]
			\cap \mathcal V_{\overline S}[k]
			\cap
			\{v:D^{(v)}_{\overline S}[k]\ge r\}
			\right|.
			\label{D:R}
		\end{align}
		
		Therefore,
		\begin{align}
			&
			\left|
			Q^S_{k+1}
			-
			\xi_n\beta_S(\bm{k}/q)
			\right|
			\nonumber\\
			&\qquad\le
			\underbrace{
				\left|
				X_S(\bm{k})
				-
				n_W\pi_S(\bm{k})
				\right|
			}_{\text{fluctuation term}}
			+
			\underbrace{
				\left|
				n_W\pi_S(\bm{k})
				-
				k_S
				-
				\xi_n\beta_S(\bm{k}/q)
				\right|
			}_{\text{mean-approximation term}}
			+
			\underbrace{|R_S[k]|}_{\text{residual term}}.
			\label{D:Q-bound}
		\end{align}
		
		The fluctuation term is controlled uniformly by
		Lemma~\ref{lemma:D-concentration}, whereas the mean-approximation
		term is $o(\xi_n)$ uniformly over $\mathbb T(\kappa)$ by
		Lemma~\ref{lemma:D-mean}.
		It remains to control the residual term.
		
		\medskip
		\noindent
		{\bf Case $q=g$.}
		For the first term in~\eqref{D:R}, if
		\[
		v\notin\mathcal S_S[k]
		\qquad\text{and}\qquad
		D^{(v)}_S[k]\ge r,
		\]
		then
		\[
		D^{(v)}_{\overline S}[k]\ge1.
		\]
		Hence this term is bounded by $C_S(\bm{k})$.
		
		For the second term, if
		\[
		v\in\mathcal S_S[k]
		\qquad\text{and}\qquad
		D^{(v)}_{\overline S}[k]\ge r,
		\]
		then necessarily
		$D^{(v)}_S[k]\ge r$.
		Therefore this term is also bounded by
		$C_S(\bm{k})$. Consequently,
		\begin{equation}
			|R_S[k]|
			\le
			2C_S(\bm{k})
			\qquad
			\text{on }
			\{\bm N[k]=\bm{k}\}.
			\label{D:R-qg}
		\end{equation}
		By~\eqref{D:C-concentration}, this residual is
		$o(q)=o(\xi_n)$ uniformly with exponentially high probability.
		
		\medskip
		\noindent
		{\bf Case $g\ll q$.}
		Here a deterministic estimate is sufficient.
		The first cardinality in~\eqref{D:R} is bounded by
		$N_S[k]$, whereas the second is bounded by
		$N_{\overline S}[k]$. Thus, on
		$\{\bm N[k]=\bm{k}\}$,
		\begin{equation}
			|R_S[k]|
			\le
			k_S+k_{\overline S}
			=
			k
			\le
			\kappa q.
			\label{D:R-largeq}
		\end{equation}
		Furthermore, in every regime with $g\ll q$,
		\begin{equation}
			\frac{q}{\xi_n}\longrightarrow0.
			\label{D:q-xi}
		\end{equation}
		Hence the residual term is again $o(\xi_n)$ uniformly.
		
		Combining the three estimates in~\eqref{D:Q-bound} yields
		\eqref{D:Q-concentration}.
	\end{proof}

	\subsection{Fixed-configuration estimates}
	
	It remains to prove Lemmas~\ref{lemma:D-mean}
	and~\ref{lemma:D-concentration}.
	The key advantage of working with deterministic configurations is
	that concentration can first be established simultaneously over
	$\bm{k}\in\mathbb T(\kappa)$ and only afterwards evaluated at
	the random configuration $\bm N[k]$. In particular, no
	independence between $\bm N[k]$ and the underlying mark field
	is required.

	\subsubsection{Proof of Lemma \ref{lemma:D-mean}}
	
	\begin{proof}
		We consider separately the different scaling regimes.
		
		\medskip
		\noindent
		{\bf Case $q=g$.}
		In this regime, uniformly over
		$\bm{k}\in{\mathbb T}(\kappa)$,
		\[
		(k_S+a_S)p=O(qp)=o(1).
		\]
		Using the standard small-probability expansion for the binomial
		distribution,
		\[
		\mathbb P(\text{Bin}(m,p)\ge r)
		=
		\frac{(mp)^r}{r!}
		\left(
		1+O(mp+m^{-1})
		\right),
		\]
		we obtain, uniformly over
		$\bm{k}\in{\mathbb T}(\kappa)$,
		\begin{equation}
			\pi_S(\bm{k})
			=
			\frac{[(k_S+a_S)p]^r}{r!}
			(1+o(1)).
			\label{D:pi-small}
		\end{equation}
		Indeed, the contribution corresponding to one or more
		$\overline S$-marks is smaller by a factor
		$O(qp)=o(1)$.
		
		Recalling the definition of the critical scale $g=q$ and
		$a_S/q\to\alpha_S$, relation~\eqref{D:pi-small} yields
		\begin{equation}
			\sup_{\bm{k}\in{\mathbb T}(\kappa)}
			\left|
			\frac{n_W\pi_S(\bm{k})}{q}
			-
			\left[
			\beta_S(\bm{k}/q)
			+
			\frac{k_S}{q}
			\right]
			\right|
			\longrightarrow0.
			\label{D:qg-mean}
		\end{equation}
		Since $\xi_n=q$, this proves~\eqref{D:mean}.
		
		Furthermore, uniformly over
		$\bm{k}\in{\mathbb T}(\kappa)$,
		\[
		\mathbb P
		\left(
		\text{Bin}(k_S+a_S,p)\ge r
		\right)
		=
		O((qp)^r),
		\]
		and
		\[
		\mathbb P
		\left(
		\text{Bin}
		(k_{\overline S}+a_{\overline S},p)\ge1
		\right)
		=
		O(qp).
		\]
		Consequently,
		\[
		\widetilde\pi_S(\bm{k})
		=
		O((qp)^{r+1}).
		\]
		Since, by the definition of $g$,
		\[
		n(qp)^r=\Theta(q),
		\]
		we obtain
		\[
		n_W\widetilde\pi_S(\bm{k})
		=
		O(q^2p)
		=
		o(q),
		\]
		because $qp\to0$. This proves~\eqref{D:mixed-mean}.
		
		\medskip
		\noindent
		{\bf Case $g\ll q\ll p^{-1}$.}
		Again $qp\to0$, and therefore the expansion
		\eqref{D:pi-small} holds uniformly over
		$\mathbb T(\kappa)$. Hence
		\[
		\sup_{\bm{k}\in{\mathbb T}(\kappa)}
		\left|
		\frac{n_W\pi_S(\bm{k})}
		{n(qp)^r}
		-
		\frac{1}{r!}
		\left(
		\frac{k_S}{q}+\alpha_S
		\right)^r
		\right|
		\longrightarrow0.
		\]
		In this regime
	$	\xi_n=n(qp)^r$
		and, by the assumptions on $q$,
		$ \frac{q}{n(qp)^r}\longrightarrow0. $
		Since $k_S\le\kappa q$,
		\[
		\sup_{\bm{k}\in{\mathbb T}(\kappa)}
		\frac{k_S}{\xi_n}
		\longrightarrow0.
		\]
		Recalling that
		\[
		\beta_S(\bm{x})
		=
		\frac{1}{r!}(x_S+\alpha_S)^r
		\]
		in this regime, relation~\eqref{D:mean} follows.
		
		\medskip
		\noindent
		{\bf Case $q=p^{-1}$.}
		Set
		\[
		\lambda_S(\bm{k})
		:=
		(k_S+a_S)p.
		\]
		Uniformly over
		$\bm{k}\in{\mathbb T}(\kappa)$,
		$\lambda_S(\bm{k})$ remains in a compact subset of
		$(0,\infty)$ and
		\[
		\lambda_S(\bm{k})
		=
		\frac{k_S}{q}+\alpha_S+o(1).
		\]
		
		By the standard Poisson approximation for the binomial law,
		uniformly over $\bm{k}\in{\mathbb T}(\kappa)$,
		$\text{Bin}(k_S+a_S,p)$ can be replaced, with an
		$o(1)$ error in total variation, by a Poisson random variable
		with mean $\lambda_S(\bm{k})$. The same holds for color
		$\overline S$. {Indeed, if $Y_S$ and $Y_{\overline S}$ are independent Poisson random
			variables with means $\lambda_S(\bm{k})$ and
			$\lambda_{\overline S}(\bm{k})$, respectively, then
			\[
			\mathbb P\big(Y_S-Y_{\overline S}\ge r\big)
			=
			\sum_{r'=r}^{\infty}\sum_{r''=0}^{r'-r}
			\frac{\lambda_S(\bm{k})^{r'}}{r'!}
			\frac{\lambda_{\overline S}(\bm{k})^{r''}}{r''!}
			e^{-\lambda_S(\bm{k})-\lambda_{\overline S}(\bm{k})}.
			\]
			Since, uniformly over $\bm{k}\in\mathbb T(\kappa)$,
			\[
			\lambda_S(\bm{k})
			=
			\frac{k_S}{q}+\alpha_S+o(1),
			\qquad S\in\{R,B\},
			\]
			the right-hand side converges uniformly to
			$\beta_S(\bm{k}/q)$.}
		It follows that
		\begin{equation}
			\sup_{\bm{k}\in{\mathbb T}(\kappa)}
			\left|
			\pi_S(\bm{k})
			-
			\beta_S(\bm{k}/q)
			\right|
			\longrightarrow0.
			\label{D:poisson-regime}
		\end{equation}
		This immediately yields the claim.
	
		\medskip
		\noindent
		{\bf Case $p^{-1}\ll q\ll n$.}
		In this regime $qp\to\infty$.
		By the definition of ${\mathbb T}(\kappa)$, the configurations
		under consideration remain uniformly on the $R$-dominant side
		of the switching surface. More precisely, since
		$\alpha_R>\alpha_B$, there exists a constant
		$\varepsilon_\kappa>0$ such that, for all sufficiently large $n$,
		\begin{equation}
			k_R+a_R
			\ge
			(1+\varepsilon_\kappa)(k_B+a_B),
			\qquad
			\forall\bm{k}\in{\mathbb T}(\kappa).
			\label{D:separation}
		\end{equation}

		{
		Set
		\[
		\mu_S:=(k_S+a_S)p,\qquad S\in\{R,B\}.
		\]
		By \eqref{D:separation},
		$	\mu_R\ge (1+\varepsilon_\kappa)\mu_B	$
		uniformly over $\bm k\in\mathbb T(\kappa)$.
		Choose $\delta_\kappa>0$ sufficiently small so that
		$
		(1-\delta_\kappa)(1+\varepsilon_\kappa)>1+\delta_\kappa.
		$
		Since $\mu_B=\Theta(qp)$ uniformly over $\mathbb T(\kappa)$,
		$qp\to\infty$, and $r$ is fixed, for all sufficiently large $n$,
		$
		(1-\delta_\kappa)\mu_R-(1+\delta_\kappa)\mu_B>r.
		$
		Hence, by a union bound and the binomial tail inequalities
		\eqref{Penrose-coda-sopra}--\eqref{Penrose-coda-sotto},
		\[
		1-\pi_R(\bm k),\;\pi_B(\bm k)
		\le
		\exp\{-\mu_R\FH(1-\delta_\kappa)\}
		+
		\exp\{-\mu_B\FH(1+\delta_\kappa)\}.
		\]
		Since $\mu_R,\mu_B=\Theta(qp)$ uniformly over
		$\mathbb T(\kappa)$, there exists $c_\kappa>0$ such that
		\begin{equation}
			\sup_{\bm{k}\in\mathbb T(\kappa)}
			\left[1-\pi_R(\bm{k})\right]
			\le e^{-c_\kappa qp},
			\qquad
			\sup_{\bm{k}\in\mathbb T(\kappa)}
			\pi_B(\bm{k})
			\le e^{-c_\kappa qp}.
			\label{D:strong-regime}
		\end{equation}
		}
		Therefore, uniformly over $\mathbb T(\kappa)$,
		\[
		\pi_R(\bm{k})\longrightarrow1,
		\qquad
		\pi_B(\bm{k})\longrightarrow0.
		\]
		These limits coincide with
		$\beta_R(\bm{k}/q)=1$ and
		$\beta_B(\bm{k}/q)=0$
		on $\mathbb T(\kappa)$.
		
		Since $\xi_n=n$, $n_W/n\to1$, and $q/n\to0$,
		relation~\eqref{D:mean} follows also in this regime.
	\end{proof}

	\subsubsection{Proof of Lemma \ref{lemma:D-concentration}}
	
	\begin{proof}
		By Lemma~\ref{lemma:D-mean}, there exists a finite constant
		$K_\kappa>0$ such that, for all sufficiently large $n$,
		\begin{equation}
			n_W\pi_S(\bm{k})
			\le
			K_\kappa\xi_n,
			\qquad
			\forall\bm{k}\in{\mathbb T}(\kappa).
			\label{D:mean-upper}
		\end{equation}
		
		For a binomial random variable $X$ with mean $\mu$,
		Bernstein's inequality gives
		\[
		\mathbb P(|X-\mu|>t)
		\le
		2\exp
		\left\{
		-\frac{t^2}{2(\mu+t/3)}
		\right\}.
		\]
		Taking $t=\delta\xi_n$ and using~\eqref{D:mean-upper},
		we obtain
		\begin{equation}
			\sup_{\bm{k}\in{\mathbb T}(\kappa)}
			\mathbb P
			\left(
			\left|
			X_S(\bm{k})
			-
			n_W\pi_S(\bm{k})
			\right|
			>
			\delta\xi_n
			\right)
			\le
			2e^{-c_1\xi_n}
			\label{D:single-X}
		\end{equation}
		for some $c_1=c_1(\kappa,\delta)>0$.
		
		Moreover,
		\[
		|{\mathbb T}(\kappa)|
		\le
		(\lfloor\kappa q\rfloor+1)^2.
		\]
		Thus, by the union bound,
		\begin{align}
			&\mathbb P
			\left(
			\sup_{\bm{k}\in\mathbb T(\kappa)}
			\left|
			X_S(\bm{k})
			-
			n_W\pi_S(\bm{k})
			\right|
			>
			\delta\xi_n
			\right)
			\nonumber\\
			&\qquad
			\le
			2(\lfloor\kappa q\rfloor+1)^2
			e^{-c_1\xi_n}.
			\label{D:union-X}
		\end{align}
		Under assumptions~\eqref{eq:bootcond}--\eqref{eq:grelq},
		\[
		\log q=o(\xi_n),
		\]
		and hence~\eqref{D:X-concentration} follows.
		
		Suppose now that $q=g$.
		By~\eqref{D:mixed-mean}, for all sufficiently large $n$,
		\[
		\sup_{\bm{k}\in{\mathbb T}(\kappa)}
		n_W\widetilde\pi_S(\bm{k})
		\le
		\frac{\delta q}{2}.
		\]
		Since $C_S(\bm{k})$ is binomial with mean
		$n_W\widetilde\pi_S(\bm{k})$, a standard binomial
		upper-tail bound yields
		\[
		\sup_{\bm{k}\in{\mathbb T}(\kappa)}
		\mathbb P
		\left(
		C_S(\bm{k})>\delta q
		\right)
		\le
		e^{-c_2q}
		\]
		for some $c_2>0$.
		A union bound over ${\mathbb T}(\kappa)$ and
		$\log q=o(q)$ prove~\eqref{D:C-concentration}.
	\end{proof}

	\subsection{Proof of Proposition \ref{prop:D2-new}}
	
	\begin{proof}
		Recall that
		\[
		\widehat N_S[j]=N_S[j]-J_S[j],
		\qquad
		J_S[j]
		=
		\sum_{h=1}^{\min\{j,n_W-1\}}U_h^S.
		\]
		Since $q=o(n)$, for every fixed $\kappa>0$ and all sufficiently
		large $n$,
		\[
		\lfloor\kappa q\rfloor<n_W-1.
		\]
		Therefore, for
		$0\le k<\lfloor\kappa q\rfloor$,
		\begin{equation}
			\widehat N_S[k+1]-\widehat N_S[k]
			=
			M^S_{k+1}-U^S_{k+1}.
			\label{D:mart-increment}
		\end{equation}
		By Proposition~\ref{M-U},
		\[
		\mathbb E
		\left[
		M^S_{k+1}\mid\mathcal H_k
		\right]
		=
		U^S_{k+1}.
		\]
		Hence
		\[
		\{\widehat N_S[k]\}_{0\le k\le\lfloor\kappa q\rfloor}
		\]
		is an $\{\mathcal H_k\}$-martingale.
		
		Moreover,
		\[
		\left|
		\widehat N_S[k+1]-\widehat N_S[k]
		\right|
		\le1
		\qquad\text{a.s.}
		\]
		Therefore, by Azuma--Hoeffding's inequality, for every
		$j\le\kappa q$,
		\[
		\mathbb P
		\left(
		|\widehat N_S[j]|>\delta q
		\right)
		\le
		2\exp
		\left\{
		-\frac{\delta^2q^2}{2j}
		\right\}
		\le
		2\exp
		\left\{
		-\frac{\delta^2q}{2\kappa}
		\right\}.
		\]
		Using the union bound,
		\begin{align}
			\mathbb P
			\left(
			\sup_{j\le\kappa q}
			|\widehat N_S[j]|>\delta q
			\right)
			&\le
			\sum_{j=1}^{\lfloor\kappa q\rfloor}
			\mathbb P
			\left(
			|\widehat N_S[j]|>\delta q
			\right)
			\nonumber\\
			&\le
			2\kappa q
			\exp
			\left\{
			-\frac{\delta^2q}{2\kappa}
			\right\}.
		\end{align}
		This proves~\eqref{D:Psi-bound} and hence
		Proposition~\ref{prop:D2-new}.
	\end{proof}
	
	Together with the reduction established at the beginning of the
	appendix, Propositions~\ref{prop:D1-new} and~\ref{prop:D2-new}
	complete the proof of Proposition~\ref{cor:BC}.
	
}

}
{

}

\section{Proof of Proposition \ref{le:inequnifbis}}\label{sec:prop5.3}

Let $i\in\mathbb N\cup\{0\}$,  $\bm k=(k_R,k_B)\in(\mathbb{N}\cup\{0\})^2$ and $k=k_R+k_B< n_W$.  By construction	
$\bm N[k+i]-\bm N[k]$ takes values on $\mathbb{I}_i$  (defined in \eqref{170124}). Hence
\begin{equation}\label{eq:ind1}
\sum_{\bm i\in \mathbb{I}_i} \ind_{\mathcal{E}_{\bm i}(k,i)}=1, \quad \text{where}\quad \mathcal{E}_{\bm i}(k,i):=\left\{\omega\in\Omega:\,\,\bm N[k+i]-\bm N[k]=\bm i\right\},\quad\bm i\in \mathbb{I}_i
\end{equation}
and
\begin{equation}\label{eq:doubleind}
 \ind_{\mathcal{E}_{\bm i}(k,i)} \ind_{\{\bm N[k]=   
	\bm k\}}=
 \ind_{\mathcal{E}_{\bm i}(k,i)} \ind_{\{\bm N[k+i]=\bm k+\bm i\}},  \quad  \text{for any } \bm i\in \mathbb{I}_i.
\end{equation}
Since $q=o(n)$, for the values of $k$ and $i$ considered below,
all indices are smaller than $n_W-1$ for all sufficiently large $n$.
So, for any $z>0$,  recalling the definition of $J_S[k]$ in  \eqref{170124-3},  we have 
\begin{align*}
&\big(J_S[k+\lfloor zq\rfloor]-J_S[k]\big)\ind_{\{\bm N[k]=\bm k\}}
=\left[\sum_{i=0}^{\lfloor zq\rfloor-1}\big(J_S[k+i+1]-J_S[k+i]\big)\right]\ind_{\{\bm N[k]=\bm k\}}\\
&\stackrel{(a)}{=}\sum_{i=0}^{\lfloor zq\rfloor-1}\sum_{\bm i\in \mathbb{I}_i}\big(J_S[k+i+1]-J_S[k+i]\big)
 \ind_{\mathcal{E}_{\bm i}(k,i)}\ind_{\{\bm N[k]=\bm k\}}\\
&=\sum_{i=0}^{\lfloor zq\rfloor-1}\sum_{\bm i\in \mathbb{I}_i}U_{k+i+1}^{S}
 \ind_{\mathcal{E}_{\bm i}(k,i)}\ind_{\{\bm N[k]=\bm k\}} 
\stackrel{(b)}{=} \sum_{i=0}^{\lfloor zq\rfloor-1}\sum_{\bm i\in \mathbb{I}_i}U_{k+i+1}^{S}
 \ind_{\mathcal{E}_{\bm i}(k,i)}\ind_{\{\bm N[k+i]=\bm k+\bm i\}},
\end{align*}
where identity $(a)$  follows from \eqref{eq:ind1} and $(b)$ from \eqref{eq:doubleind}.
Therefore, for any $y,z>0$,
\begin{align}
J_S[\lfloor yq\rfloor+\lfloor zq\rfloor]&-J_S[\lfloor yq\rfloor]=\sum_{\bm k\in \mathbb{I}_{\lfloor yq\rfloor}}\big(J_S[\lfloor yq\rfloor+\lfloor zq\rfloor]-J_S[\lfloor yq\rfloor]\big)
\ind_{\{\bm N[\lfloor yq\rfloor]=\bm k\}}\nonumber\\
&=\sum_{\bm k\in \mathbb{I}_{\lfloor yq\rfloor}}\sum_{i=0}^{\lfloor zq\rfloor-1}\sum_{\bm i\in \mathbb{I}_i}U_{\lfloor yq\rfloor+i+1}^S \ind_{\mathcal{E}_{\bm i}(k,i)}\ind_{\{\bm N[\lfloor yq\rfloor+i]=\bm k+\bm i\}}.\label{eq:relsulleZbis}
\end{align}
{Fix $\kappa<\kappa_{\bm f}$ 
and assume 
$y+2z\le\kappa$.}
For any $\bm k\in \mathbb{I}_{\lfloor yq\rfloor}$ and $i\in\{1,\ldots,\lfloor zq\rfloor-1\}$, we have
\[
k_R+k_B+i=\lfloor yq\rfloor+i\leq (y+z)q\leq(\kappa -z)q.
\]
{
Therefore, if $q\ll p^{-1}$ or $q=p^{-1}$, for any
$\bm k\in\mathbb I_{\lfloor yq\rfloor}$,
$i=0,\ldots,\lfloor zq\rfloor-1$, and
$\bm i\in\mathbb I_i$, we have
$\bm k+\bm i\in\mathbb T(\kappa)$.
If instead $q\gg p^{-1}$, we distinguish two cases.
If
$\mathbb L_{\bm k/q}(\kappa,z)\subseteq\mathbb T'(\kappa),$
then, since
$
\frac{\bm k+\bm i}{q}
\in\mathbb L_{\bm k/q}(\kappa,z),
$
we have $\bm k+\bm i\in\mathbb T(\kappa)$, and therefore,
by the definition of $\Omega_\kappa$, for every $\varepsilon>0$
and all sufficiently large $n$,
\[
 \ind_{\{\bm N[\lfloor yq\rfloor+i]
	=\bm k+\bm i\}}
\left|
U^S_{\lfloor yq\rfloor+i+1}
-
\frac{
	|\beta_S((\bm k+\bm i)/q)|
}{
	|\beta_R((\bm k+\bm i)/q)|
	+
	|\beta_B((\bm k+\bm i)/q)|
}
\right|
<\varepsilon.
\]
If
$
\mathbb L_{\bm k/q}(\kappa,z)
\not\subseteq\mathbb T'(\kappa),
$
the desired upper and lower bounds follow directly from
$0\le U^S_{\lfloor yq\rfloor+i+1}\le1$ and from
definitions \eqref{170124-8}--\eqref{170124-9}.}

Using  this relation,  the fact that
$q^{-1}(k_R+i_R,k_B+i_B)\in\mathbb{L}_{{\bm k/q}}(\kappa,z)$ (with $\mathbb{L}_{{\bm k/q}}(\kappa,z)$ defined in \eqref{170124-6}),
the definitions of
$\overline{\beta}_{S,\mathbb{L}_{{\bm k/q}}(\kappa,z)}$ and $\underline{\beta}_{S,\mathbb{L}_{{\bm k/q}}(\kappa,z)}$ (in \eqref{170124-7}  or  \eqref{170124-8} and \eqref{170124-9}) 
and the fact that $0\le U_{\lfloor yq\rfloor+i+1}  \le 1$,
it follows that,
for all $\omega\in\Omega_\kappa$ and any $\varepsilon>0$,  there exists $n(\omega,\varepsilon)$
such that,  for all $n\geq n(\omega,\varepsilon)$
\begin{align}
\ind_{\{\bm N[\lfloor yq\rfloor +i]=\bm k+\bm i\}}(\omega)(\underline{\beta}_{S,\mathbb{L}_{{\bm k/q}}(\kappa,z)}-\varepsilon)&\leq
\ind_{\{\bm N[\lfloor yq\rfloor+i]=\bm k+\bm i\}}(\omega)U_{\lfloor yq\rfloor+i+1}^S(\omega)\nonumber\\
&\leq\ind_{\{\bm N[\lfloor yq\rfloor+i]=\bm k+\bm i\}}(\omega)(\overline{\beta}_{S,\mathbb{L}_{{\bm k/q}}(\kappa,z) }+\varepsilon).\nonumber
\end{align}
Combining this relation with \eqref{eq:relsulleZbis}, we have that,
for all $\omega\in\Omega_\kappa$ and any $\varepsilon>0$,  there exists $n(\omega,\varepsilon)$
such that,  for all $n\geq n(\omega,\varepsilon)$,
\begin{align}
&\sum_{\bm k\in\mathbb{ I}_{\lfloor yq \rfloor}}\sum_{i=0}^{\lfloor zq\rfloor-1}\sum_{\bm i\in \mathbb{I}_i}  \ind_{\mathcal{E}_{\bm i}(k,i)}(\omega)\ind_{\{\bm N[\lfloor y q\rfloor+i]=\bm k+\bm i\}}(\omega)(\underline{\beta}_{S,  \mathbb{L}_{{\bm k/q}}(\kappa,z) }-\varepsilon)\nonumber \\
&\qquad\qquad\qquad\leq J_S[\lfloor yq\rfloor+\lfloor zq\rfloor](\omega)-J_S[\lfloor yq\rfloor](\omega)
\nonumber\\
&\qquad\qquad\qquad\leq\sum_{\bm k\in \mathbb{I}_{\lfloor yq\rfloor}}\sum_{i=0}^{\lfloor zq\rfloor-1}\sum_{\bm i\in \mathbb{I}_i}
 \ind_{\mathcal{E}_{\bm i}(k,i)}(\omega)\ind_{\{\bm N[\lfloor yq\rfloor+i]=\bm k+\bm i\}}(\omega)(\overline{\beta}_{S,\mathbb{L}_{{\bm k/q}}(\kappa,z)}+\varepsilon),
\nonumber
\end{align}
i.e. ,  (using \eqref{eq:ind1} and \eqref{eq:doubleind})
\begin{align}
\lfloor zq\rfloor\sum_{\bm k\in \mathbb{I}_{\lfloor yq\rfloor}}&\ind_{\{\bm N[\lfloor yq\rfloor]=\bm k\}}(\omega)(\underline{\beta}_{S,\mathbb{L}_{{\bm k/q}}(\kappa,z)}-\varepsilon)
\leq J_S[\lfloor yq\rfloor+\lfloor zq\rfloor](\omega)-J_S[\lfloor yq\rfloor](\omega)\nonumber\\
&\qquad\qquad\qquad
\leq\lfloor zq\rfloor\sum_{\bm k\in \mathbb I_{\lfloor yq\rfloor}}\ind_{\{\bm N[\lfloor yq\rfloor]=\bm k\}}(\omega)(\overline{\beta}_{S,\mathbb{L}_{{\bm k/q}}(\kappa,z)}+\varepsilon).\label{eq:semifinal}
\end{align}
We note that, for any $\omega\in\Omega_\kappa$,
\begin{align}
&N_S[\lfloor yq\rfloor+\lfloor zq\rfloor](\omega)-N_S[\lfloor yq\rfloor](\omega)\nonumber\\
&=J_S[\lfloor yq\rfloor+\lfloor zq\rfloor](\omega)-J_S[\lfloor yq\rfloor](\omega)+\widehat{N}_S[\lfloor yq\rfloor+\lfloor zq\rfloor](\omega)-\widehat{N}_S[\lfloor yq\rfloor](\omega)].\nonumber
\end{align}
Since $\lfloor yq\rfloor\le  \lfloor yq\rfloor+2\lfloor zq\rfloor\leq\kappa q$,  by the definition of $\Omega_\kappa$ 
(in \eqref{170124-4} ),  $\widehat{N}_S[k]$ (in \eqref{170124-3}) and $\Psi_S(\kappa)$ (in \eqref{D:Psi}),  we have that, for any $\omega\in\Omega_\kappa$ and any $\varepsilon>0$, there exists
$n'(\omega,\varepsilon)$ such that,  for any $n\geq n'(\omega,\varepsilon)$, we have
\[
-\varepsilon q<\widehat{N}_S[\lfloor yq\rfloor+\lfloor zq\rfloor](\omega)-\widehat{N}_S[\lfloor yq\rfloor](\omega)<\varepsilon q
\]
and so
\begin{align}
-\varepsilon q+J_S[\lfloor yq\rfloor+\lfloor zq\rfloor](\omega)-J_S[\lfloor yq\rfloor](\omega)&<
N_S[\lfloor yq\rfloor+\lfloor zq\rfloor](\omega)-N_S[\lfloor yq\rfloor](\omega)\nonumber\\
&<\varepsilon q+J_S[\lfloor yq\rfloor+\lfloor zq\rfloor](\omega)-J_S[\lfloor yq\rfloor](\omega).\nonumber
\end{align}
Combining this inequality with \eqref{eq:semifinal}, we have that, for all $\omega\in\Omega_\kappa$ and any $\varepsilon>0$,
there exists $n''(\omega,\varepsilon)$ such that, for all $n\geq n''(\omega,\varepsilon)$,
\begin{align}
&-\varepsilon q+\lfloor zq\rfloor\sum_{\bm k\in \mathbb{I}_{\lfloor yq \rfloor}}\ind_{\{\bm N[\lfloor yq\rfloor](\omega)=\bm k\}}(\underline{\beta}_{S,\mathbb{L}_{{\bm k/q}}(\kappa,z)}-\varepsilon)
\nonumber\\
&\qquad\qquad\qquad
<N_S[\lfloor yq\rfloor+\lfloor zq\rfloor](\omega)-N_S[\lfloor yq\rfloor](\omega)\nonumber\\
&\qquad\qquad\qquad
<\varepsilon q+\lfloor zq\rfloor \sum_{\bm k\in\mathbb{I}_{\lfloor yq \rfloor}}
\ind_{\{\bm N[\lfloor yq\rfloor](\omega)=\bm k\}}(\overline{\beta}_{S,\mathbb{L}_{{\bm k/q}}(\kappa,z)}+\varepsilon).
\nonumber
\end{align}

The claim follows by dividing by $q$, since the discrepancies due to the linear interpolation and to replacing $\lfloor (y+z)q\rfloor$ with $\lfloor yq\rfloor+\lfloor zq\rfloor$ are uniformly $O(1)$, while $\lfloor zq\rfloor/q\to z$.

\section{Proof of Proposition \ref{lemma-compare}}\label{sec:lemmacompare}
We divide the proof in two steps. In the first step we prove the proposition assuming $a_{R,1}=a_{R,2}$. In the second step we consider the general case.\\
\noindent{\bf \,\,Case\,\,$a_{R,1}=a_{R,2}$.}
Let $\mathcal{V}_{S,h}$, $S\in\{R,B\}$, $h\in\{1,2\}$, denote the set of $S$-seeds for the process $h$.
Note that $|\mathcal{V}_{S,h}|=a_{S,h}$. Since $a_{R,1}=a_{R,2}$ and $a_{B,1}\ge  a_{B,2}$,
we can, without loss of generality,  assume that	
$\mathcal{V}_{R,1}\equiv\mathcal{V}_{R,2}$ and $\mathcal{V}_{B,1}\supseteq\mathcal{V}_{B,2}$.  Consequently $\mathcal{V}_{W,2}\supseteq\mathcal{V}_{W,1}$
and 
\begin{equation}\label{eq:20092021I}
	\mathcal{V}_{W,2}\setminus\mathcal{V}_{W,1}=\mathcal{V}_{B,1}\setminus\mathcal{V}_{B,2}.
\end{equation}
Let $\mathcal{V}_{S,h}(t)$ and $\mathcal{W}_{S,h}(t)$ denote, respectively, the random subsets of $\mathcal{V}_{W,1}$ and $\mathcal{V}_{W,2}$,  defined on $\Omega$, consisting of $S$-active nodes at time $t$ for the process $h$.  We denote by $\mathcal{V}_{S,h}(\infty)$ and $\mathcal{W}_{S,h}(\infty)$ the corresponding random subsets of $\mathcal{V}_{W,1}$ and $\mathcal{V}_{W,2}$  formed by $S$-active nodes when the process $h$ terminates.
We will show later on that
\begin{equation}\label{eq:15092021I}
	|\mathcal{V}_{R,1}(\infty)|  \leq_{\mathrm{st}}|\mathcal{V}_{R,2}(\infty)| \quad \text{and}\quad   	|\mathcal{V}_{B,2}(\infty)| \leq_{\mathrm{st}}|\mathcal {V}_{B,1}(\infty)|;
\end{equation}
the claim then  follows immediately by observing that
$|\mathcal{V}_{S,1}(\infty)|=N_{S,1}([0,\infty)\times \mathcal{V}_{W,1})$ and
$|\mathcal{W}_{S,2}(\infty)|=N_{S,2}([0,\infty)\times \mathcal{V}_{W,2})$,
$S\in\{R,B\}$.
For instance, regarding the $B$-active nodes, we have:  
\begin{align}
	A^*_{B,1} &=| {\mathcal{V}}_{B,1}(\infty) |+a_{B,1}
	 \geq_{\mathrm{st}} | {\mathcal{V}}_{B,2}(\infty) |+a_{B,1}\nonumber\\
	&	=| {\mathcal{V}}_{B,2}(\infty) |+|  \mathcal{V}_{W,2}\setminus\mathcal{V}_{W,1} |+a_{B,2}
	 \geq_{\mathrm{st}} |{\mathcal{W}}_{B,2}(\infty) |+a_{B,2}= A^*_{B,2}, \nonumber
\end{align}
where the second equality follows from \eqref{eq:20092021I}.  The  final inequality holds because 
by construction ${\mathcal{W}}_{B,2}(\infty)\subseteq {\mathcal{V}}_{B,2}(\infty)\cup (\mathcal{V}_{W,2}\setminus
\mathcal{V}_{W,1})$.

It remains to prove  \eqref{eq:15092021I}.   We will establish \eqref{eq:15092021I}  through a coupling argument;  that is,
we will consider a probability space $(\widetilde{\Omega},\widetilde{\mathcal{F}},\widetilde{\mathbb{P}})$ and random subsets defined on it, say $\widetilde{\mathcal{V}}_{S,h}(\infty)$, $h\in\{1,2\}$, such that:
{
	\begin{equation}\label{190825}
		\begin{split}
			(i)\;& \widetilde{\mathcal{V}}_{S,h}(\infty)
			\overset{\mathrm L}{=}\mathcal{V}_{S,h}(\infty),
			\qquad S\in\{R,B\},\ h\in\{1,2\},\\
			(ii)\;&
			\widetilde{\mathcal{V}}_{R,1}(\infty)
			\subseteq\widetilde{\mathcal{V}}_{R,2}(\infty),
			\qquad
			\widetilde{\mathcal{V}}_{B,2}(\infty)
			\subseteq\widetilde{\mathcal{V}}_{B,1}(\infty),
			\quad\widetilde{\mathbb P}\text{-a.s.}
		\end{split}
	\end{equation}
	
}	
Then, \eqref{eq:15092021I}  follows immediately.
To verify \eqref{190825},
we begin defining the processes $\widetilde{N}^{'(h)}:=\sum_{v\in\mathcal{V}_{W,h}}\widetilde{N}_v'$   for $h\in\{1,2\}$,
where $\{N_v'\}_{v\in\mathcal{V}_W}$ are independent Poisson processes on 
$\widetilde{\Omega} \times [0,\infty)\times\mathcal{V}_{W_2}$ with $N_v'$ having mean measure $\dd t\delta_v(\dd \ell)$. 
Since $\mathcal{V}_{W,1}\subseteq \mathcal{V}_{W,2}$, it follows that 
$\widetilde{N}^{(1)}\subseteq\widetilde{N}^{(2)}=\widetilde{N}^{(1)}\cup(\widetilde{N}^{(2)}\setminus\widetilde{N}^{(1)})$,
$\widetilde{\mathbb P}$-almost surely.
We denote  the points of $\widetilde{N}^{'(h)}$ by $\{(\widetilde{T'}^{(h)}_k,\widetilde{V'}^{(h)}_k)\}_{k\in\mathbb N}$.
For each $v\in\mathcal{V}_{W,2}$,  we consider $\{\widetilde E_i^{(v)}\}_{i\in\mathbb N}$ and $\{\widetilde E_i^{'(v)}\}_{i\in\mathbb N}$,  which are independent sequences of independent random variables defined on $\widetilde{\Omega}$ with the Bernoulli law of mean $p$.  These sequences are assumed to be independent of $\widetilde N^{(2)}$. 

Our focus here is on the resulting coupled versions of the competing bootstrap percolation processes, which are defined on
$\widetilde\Omega$.  We
denote by $\widetilde{\mathcal {V}}_{S,h}(t)$ and  $\widetilde{\mathcal {W}}_{S,h}(t)$
the random subset of $\mathcal{V}_{W,1}$ and  $\mathcal{V}_{W,2}$, defined on $\widetilde{\Omega}$, consisting of  $S$-active nodes at time $t$. 

Observe that  the coupled processes,  namely $\widetilde{N}^{'(h)}$, $\{\widetilde E_i^{(v)}\}_{i\in\mathbb N}$ and $\{\widetilde E_i^{'(v)}\}_{i\in\mathbb N}$, are constructed to follow the same law as their original counterparts defined on  $\Omega$. Consequently, the derived quantities $\widetilde{\mathcal{V}}_{S,h}(t)$ and $\widetilde{\mathcal{V}}_{S,h}(\infty)$ are distributed  identically to  ${\mathcal{V}}_{S,h}(t)$ and ${\mathcal{V}}_{S,h}(\infty)$, respectively. 
This establishes, \eqref{190825}-$(i)$.

Moreover, by construction,  for an arbitrarily fixed $k\in\mathbb N$,  the set $\widetilde{\mathcal{V}}_{S,h}(t)$ remains  constant 
for $\widetilde T'^{(1)}_k\leq t<\widetilde T'^{(1)}_{k+1}$, and may increase (with respect to the set inclusion) by the addition of a new node of color $S$,
at time $t=\widetilde{T'}^{(1)}_{k+1}$.  Relation 
\eqref{190825}-$(ii)$  
 follows if we prove that,  for any $k\in\mathbb N$,
\begin{equation}\label{eq:15092021III}
	\widetilde{\mathcal{V}}_{R,1}\left(\widetilde{T'}^{(1),-}_{k}\right)\subseteq\widetilde{\mathcal{V}}_{R,2}\left(\widetilde{T'}_{k}^{(1),-}\right)\qquad\text{and}\qquad
	\widetilde{\mathcal{V}}_{B,2}\left(\widetilde{T'}^{(1),-}_{k}\right)\subseteq\widetilde{\mathcal{V}}_{B,1}\left(\widetilde{T'}^{(1),-}_{k}\right),
	\qquad \text{$\widetilde{\mathbb P}$-a.s.}
\end{equation}
Indeed  
 for $S\in\{R,B\}$ and $h\in\{1,2\}$,  by construction it holds 
\begin{equation}\label{eq:20092021II}
	\widetilde{\mathcal{V}}_{S,h}(\infty)
	=\bigcup_{k\in\mathbb N}\widetilde{\mathcal{V}}_{S,h}\left(\widetilde{T'}^{(1),-}_{k}\right) , 
	\quad\text{$\widetilde{\mathbb P}$-a.s.}
\end{equation}
We prove \eqref{eq:15092021III}  by induction on $k\geq 1$. 
First, observe that the base case   $k=1$  holds trivially. Indeed, for any  $h\in\{1,2\}$ and $S \in \{R,B\}$ , we have 
 $\widetilde{\mathcal{V}}_{S,h}\left(\widetilde{T'}_{1}^{(1),-}\right)=\emptyset$.
 Now, assume that \eqref{eq:15092021III} holds for $k=j$  with $j \in \mathbb{N}$.  We aim to prove that the statement also holds for $k=j+1$.
 By the inductive hypothesis  observe that the following  relations  hold $\widetilde{\mathbb P}$-almost surely:
\begin{align} 
	\widetilde N_{R,2}\left(\big[0,\widetilde{T'}^{(1)}_j\big)\times\mathcal{V}_{W,2}\right)&\geq
	\widetilde N_{R,2}\left([0, \widetilde{T'}^{(1)}_j))\times\mathcal{V}_{W,1}\right)
	\nonumber\\
	&
	=\left|\widetilde{\mathcal{V}}_{R,2}\left( \widetilde{T'}^{(1),-}_{j}\right) \right|\geq \left|\widetilde{\mathcal{V}}_{R,1}\left( \widetilde {T'}^{(1),-}_{j}\right)\right|
	= \widetilde N_{R,1}\left(\big[0, \widetilde{T'}^{(1)}_j\big)\times\mathcal{V}_{W,1}\right)\nonumber
\end{align}

and
\begin{align}
	\widetilde{N}_{B,2}\left(\big[0, \widetilde{T'}^{(1)}_{j}\big)\times\mathcal{V}_{W,1}\right)=|\widetilde{\mathcal{V}}_{B,2}\left(\widetilde{T'}^{(1),-}_{j} \right) | 
	\leq \left|\widetilde{\mathcal{V}}_{B,1}(\widetilde{T'}^{(1),-}_{j}) \right|=\widetilde N_{B,1}\left( \big[0, \widetilde{T'}^{(1)}_{j}\big) \times\mathcal{V}_{W,1}\right).\nonumber
\end{align}
From the relations established above, it follows that  $\widetilde{\mathbb P}$-a.s. for every  $v\in\mathcal{V}_{W,1}$ we have
\begin{align*}
	\widetilde{D}_{R,1}^{(v)}\left( \widetilde{T'}_{j}^{(1),-}\right) &:=\sum_{i=1}^{\widetilde N_{R,1}([0,
		\widetilde{T'}_{j}^{(1)} )\times\mathcal{V}_{W,1} )+a_{R,1}} \widetilde E_i^{(v)}
	&\leq\sum_{i=1}^{{\widetilde N}_{R,2}([0, \widetilde{T'}^{(1)}_{j})\times\mathcal{V}_{W,2})+a_{R,2}}{ \widetilde E}_i^{(v)}
	=: {\widetilde D}_{R,2}^{(v)}\left(  \widetilde{T'}^{(1),-}_{j}\right)
\end{align*}
and  
\begin{align*}
	\widetilde D_{B,2}^{(v)}( \widetilde{T'}_{j}^{(1),-}):= 
	&\sum_{i=1}^{\widetilde N_{B,2}([0,\widetilde{T'}_{j} ^{(1)})\times\mathcal{V}_{W,2})+a_{B,2}} \widetilde {E}_i^{'(v)} \le
	\sum_{i=1}^{\widetilde N_{B,2}([0,\widetilde{T'}_{j}^{(1)}  )\times\mathcal{V}_{W,1})+a_{B,1}} \widetilde {E}_i^{'(v)} \\
	&\leq\sum_{i=1}^{\widetilde N_{B,1}([0,\widetilde{T'}_{j}^{(1)} )\times\mathcal{V}_{W,1} )+a_{B,1}} \widetilde{E}_i^{'(v)}=:{\widetilde D}_{B,1}^{(v)}\left(\widetilde{T'}_{j}^{(1),-}\right).
\end{align*}
Indeed since
$a_{B,1} = |\mathcal{V}_{W,2}\setminus\mathcal{V}_{W,1}|+a_{B,2}$,  we have,
$\widetilde{\mathbb P}$-a.s.: 
\begin{align}
	\widetilde{\mathcal{S}}_{R,1}\left( \widetilde{T'}^{(1),-}_{j}\right)
	&:=\left\{v\in\mathcal{V}_{W,1}:\,\widetilde D_{R,1}^{(v)}\left(\widetilde {T'}^{(1),-}_{j} \right)-\widetilde D_{B,1}^{(v)}\left( \widetilde {T'}^{(1),-}_{j}\right)\geq r\right\} 
	\nonumber\\
	& \subseteq \left\{v\in\mathcal{V}_{W,1}:\,\widetilde D_{R,2}^{(v)}\left(\widetilde{T'}^{(1),-}_{j}\right)-\widetilde D_{B,2}^{(v)}\left( \widetilde {T'}^{(1),-}_{j}\right) \geq r\right\}	\nonumber \\
	&\subseteq  \left\{  v\in\mathcal{V}_{W,2}:\,\widetilde D_{R,2}^{(v)}\left(\widetilde{T'}^{(1),-}_{j}\right) -
	\widetilde D_{B,2}^{(v)} \left( \widetilde {T'}^{(1),-}_{j}\right)\geq r\right\}  \nonumber\\
	&=:\widetilde{\mathcal{S}}_{R,2}\left( \widetilde{T'}^{(1),-}_{j}\right).\label{eq:16092021V}
\end{align}
{ Note that  
\[
v\in\widetilde{\mathcal{V}}_{R,h}\left(\widetilde{T'}^{(1),-}_{j+1} \right) \setminus\widetilde{\mathcal{V}}_{R,h}\left(\widetilde{T'}^{(1),-}_{j}\right) ,\quad h\in\{1,2\}
\]
implies  that  
\[
v\in\widetilde{\mathcal{S}}_{R,h}\left({\widetilde T}^{(1),-}_{j}\right) \setminus\widetilde{\mathcal{V}}_{R,h}\left(\widetilde{T'}^{(1),-}_{j}\right),\quad h\in\{1,2\}.
\]
Therefore,  if $v\in\widetilde{\mathcal{V}}_{R,1}(( \widetilde{T'}^{(1),-}_{j+1}))\setminus\widetilde{\mathcal{V}}_{R,1}(( \widetilde {T'}^{(1),-}_{j}))$,  then it must be that
$v\in\widetilde{\mathcal{S}}_{R,1}\left(\widetilde{T'}^{(1),-}_{j}\right)$. By \eqref{eq:16092021V} this implies
$v\in \widetilde{\mathcal{S}}_{R,2}\left(\widetilde{T'}^{(1),-}_{j}\right)$. 
Indeed, by the inductive hypothesis, $v$ is either already $R$-active
in process~2 or is still white; in the latter case it becomes $R$-active
at the common wake-up time $\widetilde T_j^{\prime(1)}$.
Therefore  we have
$v\in\widetilde{\mathcal{V}}_{R,2}\left( \widetilde{T'}^{(1),-}_{j+1} \right)$. 
}

 This completes the proof of the first relation in \eqref{eq:15092021III}.
Observe indeed,  that   the claim follows directly from the inductive hypothesis  when $v\in\widetilde{\mathcal{V}}_{R,1}( \widetilde{T'}^{(1),-}_{j})$.
The second relation in \eqref{eq:15092021III} follows along similar lines, observing that

\begin{align*}
	\widetilde{\mathcal{S}}_{B,2}
	\left(\widetilde{T'}^{(1),-}_{j}\right)\cap\mathcal V_{W,1}
	&=
	\left\{v\in\mathcal V_{W,1}:
	\widetilde D_{B,2}^{(v)}
	\left(\widetilde{T'}^{(1),-}_{j}\right)
	-
	\widetilde D_{R,2}^{(v)}
	\left(\widetilde{T'}^{(1),-}_{j}\right)\ge r
	\right\}\\
	&\subseteq
	\left\{v\in\mathcal V_{W,1}:
	\widetilde D_{B,1}^{(v)}
	\left(\widetilde{T'}^{(1),-}_{j}\right)
	-
	\widetilde D_{R,1}^{(v)}
	\left(\widetilde{T'}^{(1),-}_{j}\right)\ge r
	\right\}\\
	&=
	\widetilde{\mathcal{S}}_{B,1}
	\left(\widetilde{T'}^{(1),-}_{j}\right).
\end{align*}

\noindent{\bf Case\,\,$a_{R,1}\leq a_{R,2}$.}
{
	To prove the general case, introduce a third
	activation process with seed counts
	\[
	(a_{R,3},a_{B,3})
	=
	(a_{R,1},a_{B,2}).
	\]
	This configuration is always admissible, since
	\[
	a_{R,3}+a_{B,3}
	=
	a_{R,1}+a_{B,2}
	\leq
	a_{R,1}+a_{B,1}
	\leq n.
	\]
	
	First compare processes~1 and~3.
	Their red seed counts are equal, whereas
	$a_{B,1}\geq a_{B,3}$.
	The result proved above therefore gives
	\begin{equation}\label{eq:21092021I}
		A_{R,1}^*
		\leq_{\mathrm{st}}
		A_{R,3}^*
		\qquad\text{and}\qquad
		A_{B,3}^*
		\leq_{\mathrm{st}}
		A_{B,1}^*.
	\end{equation}
	
	Next compare processes~3 and~2.
	Their black seed counts are equal, whereas
	$a_{R,3}\leq a_{R,2}$.
	Interchanging the roles of red and black in the
	preceding argument yields
	\begin{equation}\label{eq:21092021II}
		A_{R,3}^*
		\leq_{\mathrm{st}}
		A_{R,2}^*
		\qquad\text{and}\qquad
		A_{B,2}^*
		\leq_{\mathrm{st}}
		A_{B,3}^*.
	\end{equation}
	
	Combining \eqref{eq:21092021I} and
	\eqref{eq:21092021II} and using transitivity
	of stochastic order proves the claim.
 }


\section{Proof of Lemmas  \ref{le:200523},  \ref{lemma:stop-CB} and \ref{lemma-140223}}\label{appendix-stop}

\begin{proof} (Lemma \ref{le:200523})
	By Lemma~\ref{lemma:deterministic-marks},
	\[
	\mathcal S_R[k]
	=
	\widehat{\mathcal S}_R
	(N_R[k],N_B[k])
	\qquad\text{a.s.}
	\]
	On the event
	$\mathcal N_{\ell_R,\ell_B}^{(k)}$ we have
	\[
	N_R[k]\ge\ell_R,
	\qquad
	N_B[k]\le\ell_B.
	\]
	Therefore, by the pathwise monotonicity property
	\eqref{static-monotonicity},
	\[
	\widehat{\mathcal S}_R
	(N_R[k],N_B[k])
	\supseteq
	\widehat{\mathcal S}_R(\ell_R,\ell_B),
	\]
		which proves the first inclusion in  \eqref{pathwise-RB-bound}.
	Similarly,
	\[
	\mathcal S_B[k]
	=
	\widehat{\mathcal S}_B
	(N_R[k],N_B[k]).
	\]
	For the \Bsupra{} condition, increasing the number of
	$R$-marks and decreasing the number of $B$-marks can only decrease
	the set of \Bsupra nodes. Hence, on
	$\mathcal N_{\ell_R,\ell_B}^{(k)}$,
	\[
	\widehat{\mathcal S}_B
	(N_R[k],N_B[k])
	\subseteq
	\widehat{\mathcal S}_B(\ell_R,\ell_B),
	\]	
	It follows that
	\[
	\left\{
	|\mathcal S_R[k]|\le x
	\right\}
	\cap
	\mathcal N_{\ell_R,\ell_B}^{(k)}
	\subseteq
	\left\{
	|\widehat{\mathcal S}_R(\ell_R,\ell_B)|
	\le x
	\right\},
	\]
	and
	\[
	\left\{
	|\mathcal S_B[k]|\ge x
	\right\}
	\cap
	\mathcal N_{\ell_R,\ell_B}^{(k)}
	\subseteq
	\left\{
	|\widehat{\mathcal S}_B(\ell_R,\ell_B)|
	\ge x
	\right\}.
	\]
	Taking probabilities and using
	\eqref{static-binomial} gives
	\eqref{prob-R-bound} and \eqref{prob-B-bound}.
\end{proof}

\begin{proof} (Lemma \ref{lemma:stop-CB}).
	{
	For notational simplicity,	throughout this proof, all variables without the superscript
		$\circ$ refer to the original physical process, stopped at
		$T_{K^*}$, rather than to its artificial extension. Thus, solely
		within this proof, we write $N_S$, $\mathcal V_S$, $\mathcal S_S$,
		$D_S^{(v)}$, and $T_k^S$ in place of
		$N_S^\circ$, $\mathcal V_S^\circ$, $\mathcal S_S^\circ$,
		$D_S^{(v),\circ}$, and $T_k^{S,\circ}$, respectively.
			Likewise, throughout this proof, all variables associated with the
		stopped process refer to its natural, non-prolonged version.}
	
Note that
\[
\{(T'^{\text{stop}}_k, V'^{\text{stop}}_k)\}_k := 
\{(T'_k, V'_k)\}_k
\]
and 
\[
\{E_i^{(v),\text{stop}}\}_{i\in\mathbb N}= \{E_i^{(v)}\}_{i\in\mathbb N},\qquad 
\{{E}_i^{'(v),\text{stop}}\}_{i\in\mathbb N}=\{E_i^{'(v)}\}_{i\in\mathbb N}. 
\]
Therefore,  for $S\in \{R,B\}$, 
\[
\mathcal{V}_S^{\text{stop}}(t)=\mathcal{V}_S (t), \quad 
\text{on the event } \{t\le  Z_{\text{stop}}\}.  
\]
On the event  $\{t > Z_{\text{stop}}\}$,  we have
\[
\mathcal{V}_R^{\text{stop}}  (t) = \mathcal{V}_R^{\text{stop}}  (Z_{\text{stop} }) = 
\mathcal{V}_R  (Z_{\text{stop} })\subseteq  \mathcal{V}_R  (t).
\]
Therefore
\begin{equation}\label{stopped-DR}
	D^{(v), \text{stop}}_R (T'_k)\le   D^{(v)}_R (T'_k),  \qquad \forall\; k \in \mathbb{N}, v\in \mathcal{V}_W.
\end{equation}
We proceed proving by induction that
\begin{equation} \label{stop-indu}
	\mathcal{V}_B(T'_k)\subseteq  \mathcal{V}^{\text{stop}}_B(T'_k), \qquad \forall k\in \mathbb{N}.
\end{equation}
The relation \eqref{stop-indu} is clearly true for $k=0$,  indeed $ \mathcal{V}_B(T'_0)=
\mathcal{V}^{\text{stop}}_B(T'_0)=\emptyset$ a.s.\footnote{We recall that conventionally $T'_0=0$.}  Assume that \eqref{stop-indu} is true for any $k\le k_0$.
Then
\begin{equation}\label{stopped-DB}
	D^{(v), \text{stop}}_B ({T'}_{k_0})\ge   D^{(v)}_B ({T'}_{k_0})  \qquad \forall\; v\in \mathcal{V}_W.
\end{equation}
Combining \eqref{stopped-DR} and \eqref{stopped-DB} we have
\[
\mathcal{S}_B(T'_{k_0})\subseteq
\mathcal{S}_B^{\text{stop}}(T'_{k_0}),
\] 
{
	We now consider the node $V'_{k_0+1}$ waking up at time
	$T'_{k_0+1}$. 
	Suppose that $V'_{k_0+1}$ becomes $B$-active in the
	original process at this time. Thus,
$	V'_{k_0+1}\in
	\mathcal V_W(T'_{k_0})\cap\mathcal S_B(T'_{k_0}).
$
	 If $V'_{k_0+1}$ is already $B$-active in the stopped process, there
	 is nothing to prove. Otherwise, since it is white in the original
	 process and
	 \[
	 \mathcal V_R^{\mathrm{stop}}(T'_{k_0})
	 \subseteq \mathcal V_R(T'_{k_0}),
	 \]
	 it cannot be $R$-active in the stopped process and is therefore still
	 white there. Moreover, since
	 $V'_{k_0+1}\in\mathcal S_B(T'_{k_0})$ and
	 \[
	 \mathcal S_B(T'_{k_0})
	 \subseteq\mathcal S_B^{\mathrm{stop}}(T'_{k_0}),
	 \]
	 it is $B$-suprathreshold also in the stopped process. Hence it becomes
	 $B$-active at time $T'_{k_0+1}$.

Therefore,
	\[
	\mathcal V_B(T'_{k_0+1})
	\subseteq
	\mathcal V_B^{\mathrm{stop}}(T'_{k_0+1}),
	\]
	which completes the induction step.
}

Then \eqref{stop-CB}  immediately follows noticing that
\[
A^{*}_B= \Big|\bigcup_k  \mathcal{V}_B(T'_{k})\Big|+a_B\le \Big |\bigcup_k  \mathcal{V}^{\text{stop}}_B(T'_{k})\Big|+a_B=
A^{*,\text{stop}}_B. 
\]
Moreover, since both processes can change only at wake-up times,
\eqref{stop-indu} implies
\[
N_B^{\mathrm{stop}}(t)\ge N_B(t),\qquad t\ge 0.
\]
Consequently,
\[
T_k^{B,\mathrm{stop}}\le T_k^B,\qquad k\in\mathbb N,
\]
which completes the proof.
\end{proof}

\begin{proof} (Lemma \ref{lemma-140223}).
We prove the lemma  by contradiction.
Assume that there exists $\alpha>0$ such that $\mathbb{P}(\limsup \{X_n>\alpha\})=
\mathbb{P}(\bigcap_n\bigcup_{m\ge n}\{X_m>\alpha\} )=\beta>0$.  Then
\[
\liminf_{n\to \infty} \sum_{m\ge n} \mathbb{P}(X_m>\alpha )\ge \lim_{n\to \infty}\mathbb{P}(\bigcup_{m\ge n}\{X_m>\alpha\} )
=	 \mathbb{P}(\bigcap_n\bigcup_{m\ge n}\{X_m>\alpha \} )=\beta.
\]
Therefore
\[
\sum_{n=0}^{\infty} \mathbb{P}(X_n>\alpha )= \infty.
\]
By the assumed  stochastic ordering,  it follows
\[
\sum_{n=0}^{\infty} \mathbb{P}(Y_n>\alpha )\ge  \sum_{n=0}^{\infty} \mathbb{P}(X_n>\alpha ) = \infty.
\]
Applying the  Borel--Cantelli lemma,  this  relation implies $\mathbb{P}(\limsup\{Y_n>\alpha\})=1$,  which contradicts the hypothesis that $Y_n\to 0$ as $n\to\infty$, a.s.
\end{proof}

\section{ Discussion of the case $\alpha_B=1$ and Proof of the inequality $\psi>\tau$ \label{proof-psi-tau}}

\subsection{The case $\alpha_B=1$ }
For $\alpha_B=1$, recall that
$z_B$ is the unique positive zero of $\beta_B$. Since, for
$q\ll p^{-1}$, the $B$-component satisfies the autonomous
equation
\[
g_B'(y)=\beta_B(g_B(y)),\qquad g_B(0)=0,
\]
its scalar solution is well defined for every finite $y$ and,
by uniqueness, cannot reach the equilibrium $z_B$ at any finite
time. Since $\kappa_{\bm g}<\infty$, continuity of the scalar
solution therefore yields
$g_B(\kappa_{\bm g})<z_B.$
\subsection{ Proof of the inequality $\psi>\tau$  }

Recall that, in both regimes considered in Theorem
\ref{lemma-nuovo}, we have $q\ll p^{-1}$ and, consequently,
\[
\beta_S(x_R,x_B)=\beta_S(x_S),
\qquad S\in\{R,B\}.
\]
Therefore, by Lemma~\ref{le:transfer},
\[
\bm f(x)=\bm g(z^{-1}(x)),
\qquad
z(y):=g_R(y)+g_B(y).
\]
Moreover, from \eqref{simplified:CPcoupled},
\[
z'(y)
=
g_R'(y)+g_B'(y)
=
\beta_R(g_R(y))+\beta_B(g_B(y)).
\]

Recalling the definition of $\tau$ in \eqref{170424}, we thus have
\begin{equation}
	\tau
	=
	\int_0^{\kappa}	\frac{1}{\beta_R(\bm f(x))+\beta_B(\bm f(x))}\,\dd x
	=
	\int_0^{\kappa}	\frac{1}{z'(z^{-1}(x))}\,\dd x.
\end{equation}
Applying the change of variables
\[
x=z(y),
\qquad
\dd x=z'(y)\,\dd y,
\]
we obtain
$
\tau=z^{-1}(\kappa).
$
In both cases considered in Theorem~\ref{lemma-nuovo},
Proposition~\ref{le:coupledphi} gives
$\kappa_{\bm f}=+\infty$. Furthermore, by
Lemma~\ref{le:transfer},
\[
\kappa_{\bm f}=z(\kappa_{\bm g}),
\]
where the latter expression is understood as the limit as
$y\uparrow\kappa_{\bm g}$. Hence, for every fixed
$\kappa<\infty$,
\begin{equation}
	\tau=z^{-1}(\kappa)<\kappa_{\bm g}.
	\label{eq:tau-less-kappag}
\end{equation}

We now consider $\psi$. From \eqref{simplified:CPcoupled},
\[
g_B'(y)=\beta_B(g_B(y)),
\qquad 0\le y<\kappa_{\bm g}.
\]
Therefore,
\[
\frac{g_B'(y)}{\beta_B(g_B(y))}=1,
\]
and, for every $y<\kappa_{\bm g}$,
\[
\int_0^{g_B(y)}
\frac{1}{\beta_B(v)}\,\dd v
=
y.
\]
Letting $y\uparrow\kappa_{\bm g}$ yields
\begin{equation}
	\int_0^{g_B(\kappa_{\bm g})}
	\frac{1}{\beta_B(v)}\,\dd v
	=
	\kappa_{\bm g}.
	\label{eq:kappag-B}
\end{equation}

As specified in Step~1 of the proof of
Theorem~\ref{lemma-nuovo}, $\varepsilon>0$ is chosen sufficiently
small so that $\beta_B(v)>0$ on
$[0,g_B(\kappa_{\bm g})+\varepsilon].$
Consequently, using \eqref{eq:kappag-B},
\begin{equation}
	\psi
	=
	\int_0^{g_B(\kappa_{\bm g})+\varepsilon}
	\frac{1}{\beta_B(v)}\,\dd v
	>
	\int_0^{g_B(\kappa_{\bm g})}
	\frac{1}{\beta_B(v)}\,\dd v
	=
	\kappa_{\bm g}.
	\label{eq:psi-greater-kappag}
\end{equation}

Combining \eqref{eq:tau-less-kappag} and
\eqref{eq:psi-greater-kappag}, we finally obtain
\[
\tau<\kappa_{\bm g}<\psi,
\]
which proves the claim.
{

\end{document}